\documentclass[pdflatex,sn-mathphys-num]{sn-jnl}
\usepackage{graphicx}%
\usepackage{multirow}%
\usepackage{amsmath,amssymb,amsfonts}%
\usepackage{amsthm}%
\usepackage{mathrsfs}%
\usepackage[title]{appendix}%
\usepackage{xcolor}%
\usepackage{textcomp}%
\usepackage{manyfoot}%
\usepackage{booktabs}%
\usepackage{algorithm}%
\usepackage{algorithmicx}%
\usepackage{algpseudocode}%
\usepackage{listings}%
\usepackage{tikz}
\usetikzlibrary{shapes,arrows}
\usepackage{subfigure}
\usepackage{float}
\usepackage{caption}

\theoremstyle{thmstyleone}%
\theoremstyle{thmstyletwo}%
\theoremstyle{thmstylethree}%
\begin{document}

\title[Article Title]{Arbitrary high order shaped stencils for time domain finite difference schemes in seismic wave propagation}

\author[1]{\fnm{Peixoto} \sur{Pedro S.}}\email{ppeixoto@usp.br; ORCID: 0000-0003-2358-3221}
\author*[1]{\fnm{Silva} \sur{Felipe A. G.}}\email{felipe.augusto.guedes@usp.br; ORCID: 0000-0003-0299-8153}
\affil*[1]{\orgdiv{Instituto de Matemática e Estatística}, \orgname{Universidade de São Paulo}, \orgaddress{\street{Rua do Matão, 1010 - Cidade Universitária}, \city{São Paulo}, \postcode{05508-090}, \state{São Paulo}, \country{Brazil}}}

\abstract{Finite Difference Schemes are widely used in the approximation of different hyperbolic (wave-like) differential equations, and are particularly important for seismic wave modelling and its applications. Classical methods based on Taylor Series are dominant in the literature; however, it is known that these methods can suffer from excessive numerical dispersion. In this paper, we review and extend existing high-order in space finite difference schemes for acoustic wave propagation, featuring different stencil geometries ranging from classical cross stencils to stencils with rhombus or square-like shapes, and propose a general mathematical framework for their derivation. The numerical implementation is performed in a symbolic, high-level framework (Devito), which compiles and runs highly optimized, stencil-based computations, allowing for a low-level interpretation of the methods' efficiency. We demonstrate that non-cross stencil shapes, such as rhombus and square-based stencils, do not necessarily provide added accuracy or dispersion reduction in general, despite their increased computational cost. However, results on both idealized and realistic velocity models confirm the benefits of using dispersion-optimized cross-stencils, indicating adequate accuracy with reduced computational cost on more compact stencils compared to classic approaches. Finally, our implementation of the methods provides ease of use for full-scale acoustic seismic inversion problems using Devito.}



\maketitle

\section{Introduction}\label{intro}

Seismic imaging relies heavily on the numerical solution of wave equations. Accurate inference of geophysical subsurface properties necessitates solutions of high fidelity—meaning they must effectively suppress spurious numerical waves and accurately preserve the physical dispersion relations. While a variety of numerical schemes are available, Finite Difference (FD) schemes have become the primary choice in seismic imaging because their efficiency helps manage the extremely high computational demands of related techniques, such as full waveform inversion.\citep{tarantola1984inversion,virieux2009overview,tago2012modelling}.

Finite difference schemes have been a cornerstone of numerical methods since their inception, leading to numerous approaches for approximating the differential operators in wave equations. Although the methods of construction may vary, their fundamental result is the same: they define a set of weights applied to variable values at pre-defined grid points. In essence, every finite difference scheme is uniquely characterized by its stencil \citep{fornberg1988generation}. Accordingly, the number of grid points in the stencil, their spatial configuration, and the definition of the weights together fundamentally determine the accuracy and stability of the resulting finite difference scheme \citep{strikwerda2004finite,collatz2012numerical}. While classical finite difference schemes are typically derived using Taylor series expansions truncated to achieve a desired order of accuracy, extensive research has also focused on developing purpose-specific schemes. These specialized methods have been widely studied, particularly for the numerical solution of acoustic and elastic wave equations (e.g., \cite{virieux1984sh,virieux1986p,levander1988fourth,schroeder1999finite, boyd2006efficient, tago2012modelling,moczo20113,moczo2014finite, wang2016effective,li2017time,chen2021framework}).

Due to the local nature of finite difference schemes and their simplicity in terms of stencil application, requiring local vector multiplication only, these are known to have lower computational costs per time step compared to other schemes, such as finite elements, finite volumes, or spectral schemes. Additionally, they are generally easier to implement computationally \citep{tago2012modelling}. However, and of relevance to seismic modelling, some implementations may yield an inadequate representation of wave dispersion and phase velocity \citep{wang2014seismic}. Numerical dispersion is a common problem inherent to finite-difference schemes, causing the phase velocity to become dependent on both the grid spacing and the discretization direction (anisotropy). This numerical error typically manifests as a slowing down of the high-frequency components of the calculated wavefield when compared to the analytical phase speed. Additionally, numerical direction biases can create directional dispersion (numerical anisotropy \citep{sescu2015numerical}). Together, these numerical artifacts can deteriorate overall wave representation and should be avoided as much as possible in practical applications.

One of the goals of this paper is to review formulations of time-domain finite difference schemes for acoustic wave equations that target high-order with arbitrary-shaped stencils, with implications in dispersion reduction. While the literature on the matter is quite vast, it seems that most seismic imaging software packages still rely on classic schemes. Therefore, the main scientific question addressed in this paper is whether effective potential gains, in terms of computational savings and better representation of realistic velocity models in seismic imaging, are possible with optimized FD stencils of different stencil shapes. Existing and novel stencil approaches will be investigated with a proposed general mathematical framework, focusing on the potential accuracy gains in realistic velocity models. A secondary goal, not of less importance, is to show how such schemes can be implemented and used within a high-level symbolic computational framework, Devito \citep{devito-api,devito-compiler,kukreja2016devito}, that yields efficient computations for realistic seismic imaging problems. 

Our main focus is on the isotropic acoustic case, which is a well-established case in seismic imaging and well-developed in Devito (e.g. \cite{devito-api, dolci2022effectiveness}). Nevertheless, Devito is a powerful tool that can be extended and applied in other geophysical/models problems, like Acoustic VTI/TTI (Vertical/Tilted Transverse Isotropy) wave equations \citep{alkhalifah2000acoustic,zhang2009one,etgen2009pseudo,zhan2011acoustic} and Acoustic RTM (Reverse Time Migration) \citep{baysal1983reverse,mcmechan2008migration,whitmore1983iterative,zhang2011stable,zhang2015stable} and even in other important models such as Anisotropic Viscoelastic model \citep{yang2016review}. The methods discussed in this paper can be adapted to Acoustic VTI/TTI and RTM and can be solved numerically using Devito. In particular, for problems where the Laplacian is explicitly present, the methods can be directly used. In cases where more general second-order elliptic operators appear, the stencil optimization would have to be adapted to the particular tensors.

Our paper will be organized as follows. In Section 2, we will provide a historical overview of time-domain finite-difference methods for wave propagation and present the primary references to be considered for implementation. In Section 3, we introduce the essential concepts for a general discretization framework for the acoustic wave equation. In Section 4, we will present descriptions of the selected FD schemes with dispersion reduction properties, starting with the Classical FD Scheme and a proposed addition of weights over this scheme. In Section 5, we will discuss the framework in which we implement the proposed schemes. In Section 6, we will conduct numerical experiments comparing the proposed methods by numerically solving the Acoustic Wave Equation with homogeneous and heterogeneous velocity profiles, thereby examining the properties of the proposed methods. 

\section{Literature Review and Methods of Interest}
\label{literaturereview}

The development of finite difference methods for seismic wave propagation goes back many decades and is mainly derived from work developed in the 1960s by Alterman and collaborators \cite{abramovici1965computations, alterman1968, alterman1968propagation}, but soon started attracting more attention \cite{alford1974accuracy, kelly1976synthetic}. For 2D waves, discretizations typically evaluate operators separately in each physical dimension, and early developments (e.g. \citep{virieux1984sh, levander1988fourth}) already presented fourth-order accurate schemes for elastic waves. It was quickly noticed that high-order spatial discrete operators were of significant interest for seismic wave modelling, as described in \cite{dablain1986application, etgen1986high, etgen1988evaluating}. A very comprehensive historical description of the literature of finite difference method developments for seismic modeling is presented in \cite{moczo2007finite}.

Naturally, many other applications, aside from seismic, require solutions to wave equations and have similar requirements. Historically, purpose-specific finite difference stencils have been frequently studied in many areas, for instance, in fluid dynamics (e.g., \citep{lele1992compact}), aeroacoustics (e.g., \citep{tam1993dispersion, kim1996optimized}), and electromagnetic problems \citep{cole1995high, zygiridis2004low, sun2005optimized, finkelstein2006dispersion,finkelstein2008comprehensive}. Interestingly, even though many applications rely on the same partial differential equations, here considering acoustic or elastic wave equations, many similar ``optimized'' schemes are proposed independently in different areas. This is particularly observed comparing the literature on Finite Difference Time Domain Method (FDTD) on electromagnetic waves (e.g. \cite{cole1995high}, \cite{nehrbass1998reducing}, \cite{lan1999higher}, \cite{yang2005isotropy}, \cite{yang2006least}, \cite{finkelstein2009spectral}), and related work on seismic wave modelling (e.g. \cite{dablain1986application}, \cite{fornberg1987pseudospectral}, \cite{moczo2007finite}, \cite{liu2009new}, \cite{liu2013globally}). Independent of the application, the overall goal is usually the derivation of high-order schemes that aim at reducing numerical dispersion effects or simply improving accuracy.%

In seismic modelling, the increasing demand for better representation of wave propagation led to the development of several numerical methods. The seismic wave problem becomes particularly challenging for time-domain finite difference schemes as the wave propagation usually occurs on heterogeneous media, in which sharp velocity contrasts may exist and can create spurious numerical artefacts \citep{fornberg1987pseudospectral,yang2002finite}. \citet{holberg1987computational,jastram1993accurate} discuss the adequate representation of high-frequency waves on discrete grids for seismic modelling, and propose the construction of optimized schemes focused on minimizing the peak relative error in group velocity for a spatial frequency range. The goal here, as in other areas, is to ensure adequate representation of physically relevant seismic waves at reasonable computational cost (fewer arithmetic operations and memory accesses). Better accuracy can be achieved by exploiting higher order schemes (e.g. \cite{dong2000staggered}), or optimizing low order schemes (e.g. \cite{takeuchi2000optimally, etgen2007tutorial}).

Optimized schemes can be derived focusing on specific spatial/temporal operators or optimized for the full wave operator. The resulting properties of the numerical scheme will naturally depend on both temporal and spatial stencils, for instance, with respect to dispersion and numerical stability properties. Interestingly, finite difference schemes for the acoustic wave equation allow an interpretation in which time derivatives can be ``traded'' with spatial derivatives (e.g. \cite{bilbao2004parameterized}), which facilitates some derivations. 

Some authors focus on optimizing the temporal discretization scheme, sometimes at the cost of relying on implicit time-stepping schemes (e.g. \cite{kim2007high, bilbao2004parameterized}). Implicit time stepping schemes can indeed provide better stability and allow further degrees of freedom for better dispersion properties, but then are generally more expensive per time step compared to explicit in time schemes, since they usually require the solution of a fully coupled linear system. \citet{geller1998optimally,takeuchi2000optimally} derived optimally accurate time-space finite difference operators, which rely on an implicit time-stepping procedure. They, however, convert the method to an explicit scheme, using a predictor-corrector approach, to avoid the computational burden of the linear solvers required in the implicit method. The drawbacks of temporal implicit schemes in seismic wave modelling are two-fold:  from one viewpoint, while stability of implicit schemes provide means of using larger time-step sizes, accuracy may limit the potential of having a time-step large enough as to compensate for the added cost of the linear solver; additionally, implicit solvers usually require global grid communication, which can impose limitations on spatial parallelism on modern supercomputers. Mainly for these reasons, it seems that realistic FD full-waveform imaging tools tend to adopt explicit temporal schemes with local spatial stencils, due to high-performance computation considerations  (e.g. \cite{devito-api, fabien2017time, de2022tt}).

Optimization of explicit temporal schemes is also possible, for example, using free parameters of Runge-Kutta schemes, as in \cite{berland2006low}, aiming for dispersion or dissipation reduction, depending on the problem at hand. High-order temporal schemes, with non-FD spatial discretization, have also been pursued (e.g. \cite{chen2007high}). Due to the aforementioned possibility of trading interpretation of temporal and spatial errors in time-domain finite differences for the acoustic wave equation, optimized explicit temporal and spatial discretizations are tightly connected. 

A significant problem, primarily present in 3D, is the memory requirements of the wave solver. As discussed in \cite{wu1996analysis}, for a fixed grid resolution, with the increase in the order of accuracy of the differential discrete operator, the scheme requires fewer samples per wavelength to represent wave propagation (with reduced dispersion) correctly. In practice, higher wavenumbers can be described on coarser grids using high-order schemes. So there is a natural trade-off between reducing memory footprint, by using a high order scheme and a coarser grid, or reducing operation count (computational processing time) by using a low order scheme on a fine grid (e.g. \citep{cheng2019robust}). With the increasing possibilities of massively parallel computer systems, allowing simultaneous processing of stencil operations of the wave solver (for instance, using Field-Programmable Gate Arrays - FPGAs, Graphics Processing Units - GPUs, Advanced RISC Machines - ARM processors), the main bottleneck of 3D seismic modelling tends to be the memory consumption, making high-order methods attractive on modern supercomputers.

Following the above discussion, we will target in this work results in the literature of seismic modelling aiming for dispersion reduction that allows arbitrary order accuracy and stencil shape, and that are explicit with respect to time discretization. 

Along these lines, there are different possible strategies. For instance, \cite{yang2008finite} discusses the importance of using even-order higher spatial discretization, whereas \cite{liu2009new} ensures the order of accuracy of Taylor-based schemes is attained not only in the dimension split (cross-line stencil) points but also in other directions (e.g, diagonal). However, most authors work on parameter space optimization (maximization or minimization) that minimizes, for instance, the phase velocity error \citep{liao20092, liu2013globally,liu2014optimal,ren2015acoustic,yang2017optimal,idesman2017optimal}. 

A relevant point in considering optimized stencils is the geometry of the points used in the stencils; different dispositions can lead to varying effects in wave representation \citep{wang2016effective}. The stencils' geometry presented in the previous references can be generalized to a cross-rhombus geometry \citep{wang2016effective}. However, there are limitations in the distribution of weights in the cross-rhombus geometry, as we will see in this paper, and we will discuss more general stencils. 

Of interest for the goal of this paper, we highlight the work of \cite{liu2009practical,liu2009numerical,zhang2013optimized,wang2014seismic,ren2015acoustic}. The methods developed there are general in terms of accuracy order and focus on reducing numerical dispersion. As discussed previously, there are several types of FD schemes with different strategies to reduce the numerical dispersion. In our case, we will focus on FD schemes of \cite{liu2009new,liu2013globally,wang2016effective}, that have the desired property, reduce the numerical dispersion, allied with the capability to increase the order of approximation, keeping the good dispersion properties. Moreover, these schemes have a historical evolution: the FD Scheme proposed in \citep{liu2009new} is based on its development over the Taylor Expansion and covers more directions of propagation. The FD Scheme proposed in \cite{liu2013globally} uses the phase velocity error shown in \citep{liu2009new} and uses this function to introduce a Least Squares problem to find the optimal weights of the FD scheme. At last, the FD schemes proposed in \citep{wang2016effective} are based on \citep{liu2009new,liu2013globally}, including the possibility to introduce more weights with respect to the stencils proposed by \citep{liu2009new,liu2013globally}, aiming to increase the accuracy of the solution. 

\section{Finite Difference Formulations}
\label{fundamentals}
First, we will introduce some fundamental notation and definitions regarding the Acoustic Wave problem and discretizations. Next, we will describe general finite difference formulations that encompass all methods to be investigated.

\subsection{Acoustic Wave Equation and Discrete Domain}
The Acoustic Wave Equation represents the displacement of a wave along the domain, under the influence of a certain velocity field. Here we will represent the displacement as $u(\mathbf{x},t)$, with $(\mathbf{x},t) \in \Omega \times I$, where $\Omega\subset\mathbb{R}^{2}$ is the spatial domain, with $\mathbf{x}=(x,z)$, and $I\subset\mathbb{R}$ is the time domain. Therefore, we define $u(\mathbf{x},t):\Omega \rightarrow \mathbb{R}$ and the Acoustic Wave Equation by
\begin{equation}
\label{eqn:wave_eq}
u_{tt}(\mathbf{x},t)-(c(\mathbf{x}))^{2}\nabla^{2}u(\mathbf{x},t) = f(\mathbf{x},t),
\end{equation}
where $u_{tt}$ represents the second partial derivative with respect to time $t$ and $\nabla^{2}(\cdot)$ represents the Laplacian operator with respect to $\mathbf{x} \in \Omega$. The function $c(\mathbf{x}):\Omega\rightarrow \mathbb{R}$ is the velocity field, which is assumed to be at least piecewise-constant and positive. The external force term $f(\mathbf{x},t):\Omega_{0}\rightarrow \mathbb{R}$ models the source of waves and is usually described by a Ricker Wavelet, with a positive dominant frequency $f_{0}\in\mathbb{R}$. 
Here, we consider the acoustic wave equation obeying the homogeneous initial conditions, given by $u(\mathbf{x}, 0)  = 0 = u_t(\mathbf{x}, 0) = 0$. Furthermore, for computational simulations, it is necessary to bound the domain $\Omega$. A limited area domain is defined as $\Omega=\left[x_{I},x_{F}\right]\times\left[z_{I},z_{F}\right]$. 
To avoid wave reflections at the boundaries, absorbing boundary conditions are usually imposed. In this work, we used a classic damping boundary layer, as discussed in \cite{dolci2022effectiveness}. In this case, we add a term $\eta(\mathbf{x})u_t(\mathbf{x}, t)$ in \eqref{eqn:wave_eq} where $\eta(\mathbf{x})$ represents a damping profile. Following the damping strategy, we extend our domain $\Omega$ to $\Omega=\left[x_{I}-L_{x},x_{F}+L_{x}\right]\times\left[z_{I}-L_{z},z_{F}+L_{z}\right]$, where $L_{x}$ and $L_{z}$ are positive real values proportional to difference $x_{F}-x_{I}$ and $z_{F}-z_{I}$, respectively. The extension in $\Omega$ crates the absorbing region (damping layer) where the waves are damped, avoiding reflections to the interior domain. 

Our two-dimensional discrete domain will be defined as the set of Cartesian points on the plane,
%
$\Pi_{xz} = \displaystyle\lbrace (x_{i},z_{j})\in\mathbb{R}^{2} | \mbox{ } i=0,\cdots,n_x \mbox{ and } j=0,\cdots,n_{z} \rbrace,$
%
where $\Delta x_{i} = x_{i} - x_{i-1}$ and $x_{I} = x_{0} < x_{1} \cdots < x_{n_x} = x_{F}$ and  $\Delta z_{j} = z_{j} - z_{j-1}$ and $z_{I} = z_{0} < z_{1} \cdots < z_{n_z} = z_{F}$. Similarly, we have
$\Pi_{t} = \displaystyle\lbrace t^{k}\in\mathbb{R} | k=0,\cdots,n_t \rbrace,  $
%
where $\Delta t^{k} = t^{k} - t^{k-1}$ and $t_{I} = t^{0} < t^{1} \cdots < t^{n_t} = t_{F}$.
An important simplification that we will consider for the analyses is that the spatial parameters, $\Delta x_{i}$ and $\Delta z_{j}$, are the same for all $i\in\lbrace0,\cdots,n_x\rbrace$ and $j\in\lbrace0,\cdots,n_z\rbrace$, that is, $\Delta x_{i} = \Delta z_{j} = h > 0$. This is not a requirement of the methods to be investigated, but merely a simplification for the sake of brevity of presentation. Again, the temporal parameter $\Delta t^{k}$ is the same for all $k\in\lbrace0,\cdots,n_t\rbrace$, that is, $\Delta t^{k} = \tau > 0$. These considerations lead us to a regular partition in the space-time domains. Still, more general settings could be easily obtained with analogous derivations as the ones shown in what follows.
We will indicate $u_{i,j}^{k}$ as being the numerical approximation for $u(x_{i},z_{j},t^{k})$ at point $(x_{i},z_{j},t^{k})$, where $u$ is the solution of the wave equation. Briefly, $u_{i,j}^{k} \approx u(x_{i},z_{j},t^{k})$. 

Finally, we remark that while the mathematical derivations and software implementations are focusing on 2D domains, the extension to 3D is straightforward in both aspects, except for heavier notation.

\subsection{Finite Difference Schemes}

Finite Difference Schemes are widely used to solve numerically several types of PDEs. There are two main ways to use this type of scheme: to approximate spatial and/or temporal differential operators. In the following sections, we will discuss these types of approximations and present some relevant schemes.

\subsubsection{Spatial Discretization}
A general finite difference approximation for the 2D Laplacian operator is given by:
\begin{equation}
\label{eqn:spacexz_derivative2}
\begin{array}{lcc}
\left.\Delta u(x,z,t)\right|_{(x_{i},z_{j},t^{k})} = \left.\displaystyle\frac{\partial^{2}u(\mathbf{x},t)}{\partial x^{2}}\right|_{(x_{i},z_{j},t^{k})} +\left.\displaystyle\frac{\partial^{2}u(\mathbf{x},t)}{\partial z^{2}}\right|_{(x_{i},z_{j},t^{k})} &  & \\ \\  \approx \displaystyle\frac{1}{h^{2}}\left[a_{0,0}u_{i,j}^{k}+\sum_{{\substack{m=-M \\ m\neq 0}}}^{M}a_{m,0}u_{i+m,j}^{k}+\sum_{{\substack{m=-M \\ m\neq 0}}}^{M}a_{0,m}u_{i,j+m}^{k}\right]  & + \\  \displaystyle\frac{1}{h^{2}}\left[\displaystyle\sum_{{\substack{p=-P \\ p\neq 0}}}^{P}\sum_{{\substack{q=-Q \\ q\neq 0}}}^{Q}a_{p,q}u_{i+p,j+q}^{k}\right] & , & \\
\end{array}
\end{equation}
where we have explicitly separated line stencil coefficients (weights), obtained from dimension splitting ($a_{0,m}$ and $a_{m,0}$), from mixed coefficients of the form  $a_{p,q}, \, p\neq 0, \, q\neq 0$, to allow an easier presentation of the several methods investigated in what follows.

The range of the stencil is defined by $M>0$, along with $P \geq 0$ and $Q \geq 0$ for the mixed terms, and can be used to define the geometric shape of the stencil. We will here assume the stencil is always symmetric around the point $(x_i,z_j)$ within each axis, and also symmetric between axes. Therefore, the 2D symmetries impose a quadrant symmetry given by $a_{p,q}=a_{q,p}$, $a_{p,q}=a_{-p,q}$, $a_{p,q}=a_{p,-q}$ and $a_{p,q}=a_{-p,-q}$, resulting in a reduced form of the approximation,
\begin{equation}
\begin{array}{lcc}
\left.\Delta u(x,z,t)\right|_{(x_{i},y_{j},t^{k})} \approx \displaystyle\frac{1}{h^{2}}a_{0,0}u_{i,j}^{k}+\frac{1}{h^{2}}\displaystyle\sum_{m=1}^{M}a_{m,0}\left(u_{i+m,j}^{k}+u_{i-m,j}^{k}+u_{i,j+m}^{k}+u_{i,j-m}^{k}\right) \\ \\ + \displaystyle\frac{1}{h^{2}}\displaystyle\sum_{q=1}^{Q}\sum_{p=P_0(q)}^{P(q)}a_{p,q}\left(u_{i-p,j-q}^{k}+u_{i+p,j+q}^{k}+u_{i-p,j+q}^{k}+u_{i+p,j-q}^{k}\right).
\label{eqn:spacexz_derivative2_sym} 
\end{array}
\end{equation}
In \eqref{eqn:spacexz_derivative2_sym} we have also included bounds for the $p$ index, $P_0(q)\geq 0$ to $P(q) \geq 0$, with the possibility of dependence of $P$ and $P_0$ on $q$, which will be helpful later in the presentation of the schemes. For notation simplicity, we assume here that when $Q=0$ or $P(q)=0$, the summation of the last term is empty, therefore not accounted in the construction of the scheme.

\subsubsection{Temporal Scheme}
In this work, we will focus our effort on the investigation of the spatial approximation. Therefore, we will consider a usual centred second-order FD scheme \citep{fornberg1988generation} to approximate the time partial derivative. This is a frequent choice in seismic wave propagation \citep{bunks1995multiscale,liu2011finite}, which ensures acceptable temporal error for realistic time-step sizes.
The temporal approximation at the point $(x_{i},z_{j},t^{k})$ is given by
\begin{equation}
\label{eqn:time_derivative1}
\left.\displaystyle\frac{\partial^{2} u(\mathbf{x},t)}{\partial t^{2}}\right|_{(x_{i},z_{j},t^{k})} = \displaystyle\frac{u(x_{i},z_{j},t^{k-1})-2u(x_{i},z_{j},t^{k})+u(x_{i},z_{j},t^{k+1})}{\tau^{2}} + O(\tau^2).
\end{equation}
The error of this approximation has an order of $O(\tau^{2})$, that is, we have a second-order bound on the time derivative approximation. However, the second-order temporal errors can be traded for spatial errors and higher orders can be obtained, assuming the solution has an exact fulfillment of the wave equation, as frequently done in Lax-Wendroff schemes \citep{strikwerda2004finite}.

\subsection{Spatial schemes}
Here, we will discuss schemes that focus on attaining specific properties based only on the spatial operator approximation.
We will explore three main approaches used in the literature, but here presented in a general unified way, to define the discrete spatial stencils: 
\begin{itemize}
    \item \textit{SpatTE}, stencils derived via \textbf{Spat}ial \textbf{T}aylor \textbf{E}xpansions (e.g. \cite{fornberg1988generation, liu2009new, lele1992compact}); 
    \item \textit{SpecTE}, stencils derived via \textbf{Spec}tral \textbf{T}aylor \textbf{E}xpansions of the spatial operator (e.g. \cite{finkelstein2009spectral, liu2009new, wang2016effective, chen2019numerical}); 
    \item \textit{SpecLS}, stencils derived via \textbf{L}east \textbf{S}quares approximation of the \textbf{Spec}tral Taylor series of the spatial operator (e.g. \cite{liu2013globally, wang2016effective}).
\end{itemize}
In all approaches, the coefficients of the finite differences scheme are obtained by solely looking at the spatial errors, allowing them to be coupled with arbitrary temporal schemes.

\subsubsection{Classic spatial Taylor truncation (SpatTE)}
\label{sec:spatTE}
 For the usual Taylor expansion, we expand $u(x_{i+l},z_{j+m},t^{k})$ around $u(x_{i},z_{j},t^{k})$, to obtain, 
\begin{equation}
\label{eqn:taylor_space}
u(x_{i+l},z_{j+m}) = \sum_{r=0}^{2M}\sum_{s=0}^{r}\left[
\binom{r}{s}
\frac{(lh)^{r-s}(mh)^{s}}{s!}
\frac{\partial^{r}u(x_{i},z_{j})}{\partial x^{r-s}\partial z^{s}}\right] +
O(h^{2M+1}),
\end{equation}
where we have omitted the temporal dimension ($t^k$) for a cleaner notation. We note that the symmetry impositions lead to the cancellation of all odd-powered terms in equation \eqref{eqn:spacexz_derivative2_sym}, since,
\begin{multline*}
u(x_{i+l},z_{j+m}) +u(x_{i+l},z_{j-m})+u(x_{i-l},z_{j+m})+u(x_{i-l},z_{j-m})  =  \\ 
\sum_{r=0}^{M}\sum_{s=0}^{r}\left[
\binom{2r}{2s}
\frac{4\,l^{2r-2s} m^{2s} h^{2r}}{2r!}\frac{\partial^{2r}u(x_{i},z_{j})}{\partial x^{2r-2s}\partial z^{2s}}\right] + O(h^{2M+2}),
\end{multline*}
for $l>0$, $m>0$, and
\begin{multline*}
u(x_{i+m},z_{j}) +u(x_{i-m},z_{j})+u(x_{i},z_{j+m})+u(x_{i},z_{j-m})  =  \\
\sum_{r=0}^{M}
\frac{2\,m^{2r} h^{2r}}{(2r)!}\left[\frac{\partial^{2r}u(x_{i},z_{j})}{\partial x^{2r}}+\frac{\partial^{2r}u(x_{i},z_{j})}{\partial z^{2r}}\right] + O(h^{2M+2}),
\end{multline*}
for $l=0$, $m>0$. As a result, the approximation of the Laplacian becomes,
\begin{equation}
\label{eqn:taylor_laplacian}
\begin{array}{lcc}
\Delta u(x_{i},z_{j}) = \displaystyle\frac{1}{h^{2}}a_{0,0}u(x_{i},z_{j}) + 
\displaystyle\frac{1}{h^{2}}\displaystyle\sum_{m=1}^{M}a_{m,0}\left(\sum_{r=0}^{M}
\frac{2\,m^{2r} h^{2r}}{(2r)!}\Delta ^r u(x_{i},z_{j})\right)  +  \\
\displaystyle\frac{1}{h^{2}}\displaystyle\sum_{q=1}^{Q}\sum_{p=P_0(q)}^{P(q)}a_{p,q}\left(\sum_{r=0}^{M}\sum_{s=0}^{r}\left[
\binom{2r}{2s}
\frac{4\,p^{2r-2s} q^{2s} h^{2r}}{(2r)!}\frac{\partial^{2r}u(x_{i},z_{j})}{\partial x^{2r-2s}\partial z^{2s}}\right]\right)+ O(h^{2M}).
\end{array}
\end{equation}
where $\Delta ^r  = \partial ^{2r}/\partial x^{2r}+\partial ^{2r}/\partial z^{2r}$ and $\partial^0 u=u$.
Ensuring an agreement of the Laplacian up to terms of order $h^{2M}$ implies respecting the following linear conditions on the coefficients. For derivatives of order $r=0$, we have
\begin{equation}
\label{eqn:taylor_system_r0}
\begin{array}{lcc}
\displaystyle a_{0,0} + 4\displaystyle\sum_{m=1}^{M}a_{m,0} + \displaystyle 4 \sum_{q=1}^{Q}\sum_{p=P_0(q)}^{P(q)}a_{p,q} = 0, &\quad& (r=0).
\end{array}
\end{equation}

For derivatives of second order ($r=1$), assuming the stencil is symmetric ($a_{p,q}=a_{q,p}$), so that the equation with respect to derivatives in $x$ is the same as that in $z$, we have
\begin{equation}
\label{eqn:taylor_system_r1}
\begin{array}{lcc}
\displaystyle \sum_{m=1}^{M} m^{2} a_{m,0}
   + \displaystyle  2 \sum_{q=1}^{Q}\sum_{p=P_0(q)}^{P(q)}\,p^{2}a_{p,q} = 1 ,&\quad & (r=1).
\end{array}
\end{equation}

For $r>1$, the main equations, that do not have mixed $x$ and $z$ derivatives (which can be recovered looking at $s=0$ and $s=r$), are similar to the case of $r=1$,
\begin{equation}
\label{eqn:taylor_system_r_s0}
\begin{array}{lcc}
\displaystyle \sum_{m=1}^{M} m^{2r} a_{m,0}
   + \displaystyle  2 \sum_{q=1}^{Q}\sum_{p=P_0(q)}^{P(q)}\,p^{2r}a_{p,q} = 0 ,&\quad & (r>1, s=0, r).
\end{array}
\end{equation}
Noting that the term with $r>1$ and $s=r$ has the same form, due to the stencil symmetries. Finally, the terms with mixed derivatives ($1<r \leq M$, $0<s \leq r/2$),
\begin{equation}
\label{eqn:taylor_system_r_s}
\begin{array}{lcc}
\displaystyle\displaystyle\sum_{q=1}^{Q}\sum_{p=P_0(q)}^{P(q)} \left(
p^{2r-2s} q^{2s} \right) a_{p,q} = 0,
\end{array}
\end{equation}
We point out that the equations for $r/2 < s < r$ are identical to those due to the symmetry of the stencils. Overall, we have therefore $M+1$ main equations, which are associated with non-mixed derivatives ($s=0$), given by equations \eqref{eqn:taylor_system_r0},\eqref{eqn:taylor_system_r1},\eqref{eqn:taylor_system_r_s0}, and $ \lfloor M/2 \rfloor  \lceil M/2 \rceil$ equations relative to the mixed derivatives \eqref{eqn:taylor_system_r_s}. The symbols $ \lfloor \cdot \rfloor$ and $ \lceil \cdot \rceil$ here refer to integer floor and ceiling, respectively.
We note, however, that the sub-system given by the sub-system \eqref{eqn:taylor_system_r_s} is independent of that of the main equations. Depending on specific choices of $Q$, $P$ and $P_{0}$, this sub-system will be either:
\begin{itemize}
    \item Determined: same number of independent equations and coefficients, $ \lfloor M/2 \rfloor  \lceil M/2 \rceil $. In this case, the solution will naturally be to have all coefficients zero ($a_{p,q} = 0$).
    \item Overdetermined: more independent equations than coefficients. In this case, the overdetermined system still has a valid solution that satisfies all equations, with all coefficients being zero ($a_{p,q} = 0$).
    \item Underdetermined: fewer independent equations than coefficients. In this case, the underdetermined system will have an infinite number of solutions (non-trivial kernel). However, it still has a valid solution with all coefficients zero ($a_{p,q} = 0$).
\end{itemize}

As a consequence, for the classic spatial Taylor truncation (SpatTE) discretization of the Laplacian, it is usually required to only solve the main equations, not the mixed term equations, by assuming $a_{p,q} = 0$, for $q=1, ..., Q$, $p=P_0(q), ..., P(q)$, resulting in the system, 
\begin{equation}
\begin{array}{lcc}
\displaystyle a_{0,0} + 4\displaystyle\sum_{m=1}^{M}a_{m,0}  = 0, &\quad& (r=0), \\
\displaystyle \sum_{m=1}^{M} m^{2} a_{m,0} = 1 ,&\quad & (r=1),
\\
\displaystyle \sum_{m=1}^{M} m^{2r} a_{m,0}  = 0 ,&\quad & (1<r\leq M),
\label{eqn:taylor_system}
\end{array}
\end{equation}
with $M+1$ equations and $M+1$ unknowns (coefficients).

\subsubsection{Spectral Taylor truncation (SpecTE)}
Following plane-wave theory, let $u(x,z,t)$ be of a wave-like form in space. Omitting the temporal dimension, we have, 
$
u(x,z) = u_0 e^{i(k_{x}x+k_{z}z)}
$,
for some constant $u_0$, where $k_{x}, k_z$ are the spatial wavenumbers. This assumption is derived from spectral theory in terms of Fourier analysis, and is frequently done for wave equations due to the natural connections between linear waves and Fourier theory. Therefore, in this form, $u$ can be seen as one of the terms of a full space-time Fourier expansion of the solution of the wave equation.

The spectral truncation approach (e.g. \cite{finkelstein2009spectral, liu2013globally, wang2016effective, chen2019numerical}) looks at approximations of the Laplacian in spectral space (the space of $(k_x,k_z)$ instead of $(x,z)$). From equation \eqref{eqn:taylor_laplacian} we now have,
\small{
\begin{equation}
\label{eqn:taylor_laplacian_spec_kxkz}
\begin{array}{lcc}
\Delta u(x_{i},z_{j}) = \displaystyle\frac{1}{h^{2}}a_{0,0}u(x_{i},z_{j}) + 
\displaystyle\frac{1}{h^{2}}\displaystyle\sum_{m=1}^{M}a_{m,0}\left(\sum_{r=0}^{M}
\frac{2\,m^{2r} h^{2r}}{(2r)!} (-1)^r (k_x^{2r}+k_z^{2r}) u(x_{i},z_{j})\right)  \\ +
\displaystyle\frac{1}{h^{2}}\displaystyle\sum_{q=1}^{Q}\sum_{p=P_0(q)}^{P(q)}a_{p,q}\left(\sum_{r=0}^{M}\sum_{s=0}^{r}\left[
\binom{2r}{2s}
\frac{4\,p^{2r-2s} q^{2s} h^{2r}}{(2r)!} (-1)^r (k_x^{2r-2s}k_z^{2s})u(x_{i},z_{j}) \right]\right) \\ + O(h^{2M}).
\end{array}
\end{equation}
}
Knowing that the analytical Laplacian of the wave-like function is given by
\begin{equation}
\label{eqn:lap_spec_exact}
    \Delta u(x_{i},z_{j})=\left((\mathrm{i}k_x)^{2}+(\mathrm{i}k_z)^{2}\right) u(x_{i},z_{j}) = -(k_x^2+k_z^2) u(x_{i},z_{j}),
\end{equation}
we can impose that the Taylor expansion matches the analytical expression for each polynomial term with respect to $k_x$ and $k_z$. As a result, we have the same linear system as the one obtained in the previous subsection for the usual spatial Taylor expansion (eqs. \eqref{eqn:taylor_system_r0} - \eqref{eqn:taylor_system_r_s}).
Again, since the sub-system of equations involving only mixed $k_xk_z$ terms only depend on $a_{p,q}$ coefficients and have a null right-hand-side, a valid general solution leads to $a_{p,q}=0$, for $q=1, ..., Q$, $p=P_0(q), ..., P(q)$, reducing the problem to the system \eqref{eqn:taylor_system}.

Alternatively (e.g. \cite{liu2009new}), writing $(k_x,k_z)$ in polar coordinates, that is, $(k_x,k_z)=k e^{\mathrm{i}\theta}$, so that $k_x=k\cos{\theta}$ and $k_z=k\sin{\theta}$, with $k=\|(k_x,k_z)\|_ 2=\sqrt{k^2_x+k_z^2}$ and $\theta\in\left[0,2\pi\right]$, the Taylor approximation of the Laplacian from equation \eqref{eqn:taylor_laplacian} becomes,
\small{
\begin{equation}
\label{eqn:taylor_laplacian_spec}
\begin{array}{lcc}
h^{2}\Delta u(x_{i},z_{j}) = \displaystyle a_{0,0}u(x_{i},z_{j}) \\ \\ + 
\displaystyle\displaystyle\sum_{m=1}^{M}a_{m,0}\left(\sum_{r=0}^{M}
\frac{2\,m^{2r} h^{2r}}{(2r)!}(-1)^r(k)^{2r}\left( \cos^{2r}(\theta) + \sin^{2r}(\theta)\right) u(x_{i},z_{j})\right)  +  \\ \\
\displaystyle\sum_{q=1}^{Q}\sum_{p=P_0(q)}^{P(q)}a_{p,q}\left(\sum_{r=0}^{M}\sum_{s=0}^{r}
\binom{2r}{2s}
\frac{4\,p^{2r-2s} q^{2s} h^{2r}(-1)^r(k)^{2r}(\cos(\theta))^{2r-2s}(\sin(\theta))^{2s}u(x_{i},z_{j})}{(2r)!} \right) \\ \\ + O(h^{2M}).
\end{array}
\end{equation}
}
From the analytical Laplacian of the wave-like function, $ \Delta u= -k^2 u$, we can impose that the Taylor expansion matches the analytical for each polynomial term with respect to $k$. As a result, we have a linear system similar to the one obtained in the previous subsection for the usual spatial Taylor expansion, given by,
\begin{equation}
\label{eqn:taylor_system_r0_spec}
\begin{array}{lcc}
\displaystyle a_{0,0} + 4\displaystyle\sum_{m=1}^{M}a_{m,0} + \displaystyle 4 \sum_{q=1}^{Q}\sum_{p=P_0(q)}^{P(q)}a_{p,q} = 0, &\quad& (r=0),
\end{array}
\end{equation}

\begin{equation}
\label{eqn:taylor_system_r1_spec}
\begin{array}{lcc}
\displaystyle \sum_{m=1}^{M} m^{2} a_{m,0}
   + \displaystyle  2 \sum_{q=1}^{Q}\sum_{p=P_0(q)}^{P(q)}\,\psi_{1}(\theta) a_{p,q} = 1 ,&\quad & (r=1),
\end{array}
\end{equation}
\begin{equation}
\label{eqn:taylor_system_r_s0_spec}
\begin{array}{lcc}
\displaystyle \sum_{m=1}^{M} m^{2r} (\cos^{2r}(\theta)+\sin^{2r}(\theta) )a_{m,0} +  \displaystyle  2 \sum_{q=1}^{Q}\sum_{p=P_0(q)}^{P(q)}\,\psi_{r}(\theta) a_{p,q} = 0 ,& & (r>1),
\end{array}
\end{equation}
 where 
 \begin{equation}
     \psi_{r}(\theta) = \sum_{s=0}^{r}
\binom{2r}{2s}
(p^{2r-2s} q^{2s} )(\cos(\theta))^{2r-2s}(\sin(\theta))^{2s}.
 \end{equation}

Overall, this approach leads to $M+1$ equations, but, in principle, needs to be satisfied for all $\theta \in [0,2\pi]$ to ensure accuracy of order $2M$ in spectral (Fourier) space. 

The main difference of this spectral approach from the previous approach is that now the linear system depends on a choice of the wavenumber angle in two dimension space, that is, it depends on a choice of $\theta$, and, as we will see in the following sections, $\theta$ is frequently used a free parameter to ensure reduction in dispersion error.
 
Finally, the \textit{SpecTE} schemes will therefore be of two forms: (i) either it is identical to the \textit{SpatTE}, if full ($k_x$, $k_z$) wavenumebers are matched with Taylor series, or, (ii) it is dependent on the choice of angle ($\theta$), and will not necessarily achieve spatial order of $O(h^{2M})$. We will usually refer to the former as simply \textit{SpecTE} (or \textit{SpatTE}), and the latter as \textit{SpecTE}-$\theta$.

\subsubsection{Least Squares Spectral Approximation of Spatial Operator (SpecLS)}
\label{sec:spec_LS}
The previous approaches employ exact agreement of truncated series (up to a certain degree) and analytical representation of the Laplacian. Alternatively, one can impose a minimization approach focused on reducing errors of interest (e.g. \cite{liu2013globally}, \cite{yang2006least}). 
From equations \eqref{eqn:lap_spec_exact} and \eqref{eqn:taylor_laplacian_spec} , we have that
\begin{equation}
\label{eqn:taylor_laplacian_spec_disp}
\begin{array}{lcc}
-k^2 = \displaystyle\frac{1}{h^{2}}a_{0,0} + 
\displaystyle\frac{1}{h^{2}}\displaystyle\sum_{m=1}^{M}a_{m,0}\left(\sum_{r=0}^{M}
\frac{2\,m^{2r} (hk)^{2r}}{(2r)!}(-1)^r\left( \cos^{2r}(\theta) + \sin^{2r}(\theta)\right) \right)  +  \\ \\
\displaystyle\frac{1}{h^{2}}\displaystyle\sum_{q=1}^{Q}\sum_{p=P_0(q)}^{P(q)}a_{p,q}\left(\sum_{r=0}^{M}\sum_{s=0}^{r}
\binom{2r}{2s}
\frac{4\,p^{2r-2s} q^{2s} (kh)^{2r}}{(2r)!}(-1)^r (\cos(\theta))^{2r-2s}(\sin(\theta))^{2s}\right) \\ \\ + O(h^{2M}).
\end{array}
\end{equation}
Assuming that the $a_{0,0}$ is exactly represented, via equation \eqref{eqn:taylor_system_r0_spec}, we may write equation \eqref{eqn:taylor_laplacian_spec_disp} as,
\begin{equation}
\label{eqn:taylor_laplacian_spec_disp2}
\begin{array}{lcc}
(kh)^2 =
\displaystyle\sum_{m=1}^{M}a_{m,0}\left(4-\sum_{r=0}^{M}
\frac{2\,m^{2r} (kh)^{2r}}{(2r)!}(-1)^r\left( \cos^{2r}(\theta) + \sin^{2r}(\theta)\right) \right)  +  \\ \\
\displaystyle\sum_{q=1}^{Q}\sum_{p=P_0(q)}^{P(q)}a_{p,q}\left(4-\sum_{r=0}^{M}\sum_{s=0}^{r}
\binom{2r}{2s}
\frac{4\,p^{2r-2s} q^{2s} (kh)^{2r}}{(2r)!}(-1)^r (\cos(\theta))^{2r-2s}(\sin(\theta))^{2s}\right) \\ \\+ O(h^{2M}).
\end{array}
\end{equation}

Following \cite{liu2013globally}, but now on a more general framework, we define a set of basis functions for a Least Squares approximation by letting $\beta = kh$, $\beta \in (0, \pi)$, and defining
\begin{equation}
\label{eqn:taylor_laplacian_spec_disp_phim}
\displaystyle \varphi_m (\beta, \theta)=
\displaystyle 4-\sum_{r=0}^{M}
\frac{2\,m^{2r} \beta^{2r}}{(2r)!}(-1)^r\left( \cos^{2r}(\theta) + \sin^{2r}(\theta)\right),
\end{equation}
\begin{equation}
\label{eqn:taylor_laplacian_spec_disp_phipq1}
\displaystyle \varphi_{p,q} (\beta, \theta)= 4-\sum_{r=0}^{M}\sum_{s=0}^{r}
\binom{2r}{2s} \frac{4\,p^{2r-2s} q^{2s} \beta^{2r}}{(2r)!}(-1)^r (\cos(\theta))^{2r-2s}(\sin(\theta))^{2s},
\end{equation}
and
\begin{equation}
\label{eqn:phibeta1}
\phi_{1}(\beta) = \beta^2.
\end{equation}
Which can be shown to be linearly independent in $\mathcal{C}([0,\pi] \times [0, 2\pi], \mathbb{R})$, the space of continuous functions with respect to $\beta$ and $\theta$, due to linear independence of the polynomials and trigonometric polynomials that form the basis.
The Least Squares problem is then that of finding $a_{m,0}$, $m=1, ..., M$, and $a_{p,q}$, $q=1, ..., Q$, $p=P_0(q), ..., P(q)$, such that the quadratic error is minimum,
Where we have assumed that the minimization with respect to $\beta$ can be focused on a specific range of the spectrum, limited to $\beta \in [ 0, b ] \subset [0, \pi]$, better capturing lower wavenumbers.
The unique Least Squares solution of the problem is obtained by solving the linear system,
\begin{equation}
\label{eqn:taylor_laplacian_spec_disp_linsys}
\begin{array}{lcc}
\displaystyle
\sum_{m=1}^{M}\langle \varphi_m, \varphi_{m'} \rangle a_{m,0}+ \sum_{q=1}^{Q}\sum_{p=P_0(q)}^{P(q)} \langle \varphi_{p,q}, \varphi_{m'} \rangle a_{p,q} = \langle \phi_{1}(\beta), \varphi_{m'} \rangle, 
 \\
\displaystyle
\sum_{m=1}^{M}\langle \varphi_m, \varphi_{p', q'} \rangle a_{m,0}+ \sum_{q=1}^{Q}\sum_{p=P_0(q)}^{P(q)} \langle \varphi_{p,q}, \varphi_{p', q'} \rangle a_{p,q} = \langle \phi_{1}(\beta), \varphi_{p', q'} \rangle ,
\end{array}
\end{equation}
for $m' = 1, ..., M$,  $q'=1, ..., Q$, $p'=P_0(q'), ..., P(q')$ and 
\begin{equation}
    \label{eqn:cont_dotprod}
    \langle f(\beta, \theta), g(\beta, \theta) \rangle = \int_{0}^{b}\int_{0}^{2\pi} f(\beta, \theta)  g(\beta, \theta) d\theta \, d\beta \, .
\end{equation}
On the other hand, following \cite{liu2013globally}, if we start with the equation (\ref{eqn:spacexz_derivative2_sym}) and apply $u(x,z)=u_0 e^{i(k_{x}x+k_{z}z)}$ together with the previous considerations about spectral expansion, we will obtain:
\begin{equation}
\label{eqn:taylor_laplacian_spec_disp_cont}
\begin{array}{lcc}
\beta^2 \approx \displaystyle\sum_{m=1}^{M}a_{m,0}\left(4 - 2\cos(m\beta\cos(\theta))-2\cos(m\beta\sin(\theta))\right)  + \\ \\ \displaystyle\sum_{q=1}^{Q}\sum_{p=P_0(q)}^{P(q)}a_{p,q}\left(4-2\cos(q\beta\cos(\theta)+p\beta\sin(\theta))-2\cos(q\beta\cos(\theta)-p\beta\sin(\theta))\right) \\ \\+ O(h^{2M}).
\end{array}
\end{equation}
Now, to simplify the previous equations, we define:
\begin{equation}
\label{eqn:phim_continuous}
\displaystyle \Phi_{m} (\beta, \theta)=
\displaystyle \frac{4-2\cos(m\beta\cos(\theta))-2\cos(m\beta\sin(\theta))}{\Psi_{1}(\beta)},
\end{equation}
\begin{equation}
\label{eqn:phipq_continuous}
\displaystyle \Phi_{p,q} (\beta, \theta)= \frac{4-2\cos(q\beta\cos(\theta)+p\beta\sin(\theta))-2\cos(q\beta\cos(\theta)-p\beta\sin(\theta))}{\Psi_{1}(\beta)}.
\end{equation}
where
\begin{equation}
\label{eqn:psibeta1}
\Psi_{1}(\beta) = \beta^2,
\end{equation}
Using the previous definition, we have the following equation:
\begin{equation}
\label{eqn:taylor_laplacian_spec_disp_cont2}
\begin{array}{lcc}
1 \approx \displaystyle\sum_{m=1}^{M}a_{m,0}\Phi_m(\beta,\theta)  + \displaystyle\sum_{q=1}^{Q}\sum_{p=P_0(q)}^{P(q)}a_{p,q}\Phi_{p,q}(\beta,\theta)+ O(h^{2M}).
\end{array}
\end{equation}
We can define the relative error, using \eqref{eqn:taylor_laplacian_spec_disp_cont2}, as:
\begin{equation}
    \label{eqn:LS_spec_disp_error}
    \displaystyle E_{\text{disp}}(b) = \int_{0}^{b}\int_{0}^{2\pi} \left( \sum_{m=1}^{M}a_{m,0}\Phi_m(\beta, \theta) + \sum_{q=1}^{Q}\sum_{p=P_0(q)}^{P(q)}a_{p,q} \Phi_{p,q}(\beta,\theta) - 1 \right)^2 d \theta \, d\beta .
\end{equation}
We have assumed that the minimization with respect to $\beta$ can be focused on a specific range of the spectrum, limited to $\beta \in [0,b] \subset [0,\pi]$, which better captures lower wavenumbers.
The unique Least Squares solution of the problem is obtained by solving the linear system,
\begin{equation}
\label{eqn:specls_cont_linsys}
\begin{array}{lcc}
\displaystyle
\sum_{m=1}^{M}\langle \Phi_m, \Phi_{m'} \rangle a_{m,0}+ \sum_{q=1}^{Q}\sum_{p=P_0(q)}^{P(q)} \langle \Phi_{p,q}, \Phi_{m'} \rangle a_{p,q} = \langle 1, \Phi_{m'} \rangle, 
 \\
\displaystyle
\sum_{m=1}^{M}\langle \Phi_m, \Phi_{p', q'} \rangle a_{m,0}+ \sum_{q=1}^{Q}\sum_{p=P_0(q)}^{P(q)} \langle \Phi_{p,q}, \Phi_{p', q'} \rangle a_{p,q} = \langle 1, \Phi_{p', q'} \rangle ,
\end{array}
\end{equation}
for $m' = 1, ..., M$,  $q'=1, ..., Q$, $p'=P_0(q'), ..., P(q')$ and using the inner product $\langle\cdot,\cdot\rangle$ defined in \eqref{eqn:cont_dotprod}. 

The systems \eqref{eqn:taylor_laplacian_spec_disp_linsys} and \eqref{eqn:specls_cont_linsys} have a similar basic idea, where both use the spectral expansion. The main difference is the starter equation: for \eqref{eqn:taylor_laplacian_spec_disp_linsys} we use \eqref{eqn:taylor_laplacian}, which was expanded in Taylor's series, and after that we apply the approximation to obtain \eqref{eqn:taylor_laplacian} in the spectral space, adding some relations in the spectral space to simplify and reduce to \eqref{eqn:taylor_laplacian_spec_disp_linsys}. On the other hand, in \eqref{eqn:specls_cont_linsys} we started with \eqref{eqn:spacexz_derivative2_sym}, which isn't expanded in terms of a Taylor series, and we use the plane-wave theory to have a wave-like form in space, omitting the temporal dimension. There isn't a clear relationship between \eqref{eqn:taylor_laplacian_spec_disp_linsys} and \eqref{eqn:specls_cont_linsys}, but following \citep{liu2013globally,wang2016effective} we observe that the cosine expansion in the right places of \eqref{eqn:taylor_laplacian_spec_disp_cont}, adding binomial relations, can produce a similar result as \eqref{eqn:taylor_laplacian_spec_disp2}. In the direction of maintaining coherence with existing results in the literature, we will adopt \eqref{eqn:taylor_laplacian_spec_disp_linsys} to define the SpecLS scheme, but it is essential to highlight that \eqref{eqn:taylor_laplacian_spec_disp_linsys}was produced using a different approach based on spectral expansion and can produce significant results in some numerical problems.  

\subsection{Dispersion Reduction Schemes}
The schemes discussed previously do not account for errors originating from temporal discretizations. Now we will discuss schemes that estimate the spatial coefficients also taking into account the temporal effects in the numerical dispersion relations of the method. The schemes we will consider are based on the following:
\begin{itemize}
    \item \textit{DispTE}, stencils obtained via optimal choices for \textbf{Disp}ersion reduction of the spectral \textbf{T}aylor \textbf{E}xpansions of the overall time-space discrete operator (e.g. \cite{finkelstein2009spectral, liu2009new, wang2016effective, chen2019numerical}); 
    \item \textit{DispLS}, stencils obtained via \textbf{L}east \textbf{S}quares approximation of the spectral Taylor series of the time-space discrete operator, aiming at \textbf{Disp}ersion error reduction (e.g. \cite{liu2013globally, wang2016effective}).
\end{itemize}

\subsubsection{Truncated Dispersion relations}
Since one of the main goals is to obtain schemes that reduce dispersion error, we here briefly discuss the dispersion relations, in general, of the formulation we will use. To analyse the numerical dispersion relations of the schemes, we will use the wave-like solutions, assuming that $f(\mathbf{x},t)=0$ and a constant velocity ($c(x,z)=c$), so that $u(x,z,t)=u_0\,e^{i(k_xx+k_zz-\omega t)}$. Using the continuous equation \eqref{eqn:wave_eq}, the analytical dispersions relations result in $\omega = c k$, 
where $k=\|(k_x,k_z)\|_2=\sqrt{k_{x}^{2}+k_{z}^{2}}$.
For the numerical schemes, the imposition of a wave-like solution results in slightly more complicated expressions for the dispersion relations, now dependent on the temporal and spatial steps. The temporal discretization (see Equation \eqref{eqn:time_derivative1}) leads to
\begin{equation}
\label{eqn:time_derivative2}
\left.\displaystyle\frac{\partial^{2} u(x, z,t)}{\partial t^{2}}\right|_{(x_{i},z_{j},t^{k})} = \frac{2u(x_{i},z_{j},t^{k})}{\tau^2} 
(\displaystyle \cos(\omega \tau)-1) + O(\tau^2),
\end{equation}
whereas the spatial discretization imposes what is presented in Equation \eqref{eqn:taylor_laplacian_spec_kxkz}, so that the numeric $\omega$ relation to ($k_x$, $k_z$) is given by
\begin{equation}
\label{eqn:full_disp}
\begin{array}{lcc}
\displaystyle \cos(\omega \tau)  =1+ \displaystyle\frac{C^2}{2}a_{0,0} + 
\displaystyle\frac{C^2}{2}\displaystyle\sum_{m=1}^{M}a_{m,0}\left(\sum_{r=0}^{M}
\frac{2\,m^{2r} h^{2r}}{(2r)!} (-1)^r (k_x^{2r}+k_z^{2r}) \right)  +  \\
\displaystyle\frac{C^2}{2}\displaystyle\sum_{q=1}^{Q}\sum_{p=P_0(q)}^{P(q)}a_{p,q}\left(\sum_{r=0}^{M}\sum_{s=0}^{r}\left[
\binom{2r}{2s}
\frac{4\,p^{2r-2s} q^{2s} h^{2r}}{(2r)!} (-1)^r (k_x^{2r-2s}k_z^{2s}) \right]\right) \\ \\ + C^2 O(h^{2M+2})+O(\tau^4),
\end{array}
\end{equation}
where we have denoted the Courant number as
\begin{equation}
\label{eqn:courant}
    C=\frac{c \tau}{h},
\end{equation}
which leads to numeric dispersion relation, here denoted as $\omega_{FD}=\omega_{FD}(k_x, k_z)$,
\begin{equation}
\label{eqn:full_disp_omegaFD}
\begin{array}{lcc}
\omega_{FD}  =\displaystyle\frac{1}{\tau} \arccos  \left\{1+ \displaystyle\frac{C^2}{2}a_{0,0} + 
\displaystyle\frac{C^2}{2}\displaystyle\sum_{m=1}^{M}a_{m,0}\left(\sum_{r=0}^{M}
\frac{2\,m^{2r} h^{2r}}{(2r)!} (-1)^r (k_x^{2r}+k_z^{2r}) \right)  +  \right.\\
 \left.\displaystyle\frac{C^2}{2}\displaystyle\sum_{q=1}^{Q}\sum_{p=P_0(q)}^{P(q)}a_{p,q}\left(\sum_{r=0}^{M}\sum_{s=0}^{r}\left[
\binom{2r}{2s}
\frac{4\,p^{2r-2s} q^{2s} h^{2r}}{(2r)!} (-1)^r (k_x^{2r-2s}k_z^{2s})\right]\right)\right\}  .
\end{array}
\end{equation}
Equivalently, in polar coordinates of $(k_x, k_z)$, leading to $(k, \theta)$, we have
\begin{equation}
\omega_{FD}^{\theta}=\omega_{FD}(k, \theta),    
\end{equation}
\scriptsize{
\begin{equation}
\label{eqn:full_disp_omegaFDtheta}
\begin{array}{lcc}
\omega^{\theta}_{FD}  =\displaystyle\frac{1}{\tau} \arccos \left\{ 1+ \displaystyle\frac{C^2 }{2 }a_{0,0} + \displaystyle C^2 \displaystyle\sum_{m=1}^{M}a_{m,0}\left(\sum_{r=0}^{M}
\frac{\,m^{2r} h^{2r}}{(2r)!}(-1)^r(k)^{2r}\left( \cos^{2r}(\theta) + \sin^{2r}(\theta)\right) \right)  + \right. \\ \\
\left. \displaystyle C^2 \displaystyle\sum_{q=1}^{Q}\sum_{p=P_0(q)}^{P(q)}a_{p,q}\left(\sum_{r=0}^{M}\sum_{s=0}^{r}
\binom{2r}{2s}
\frac{2\,p^{2r-2s} q^{2s} h^{2r}}{(2r)!}(-1)^r(k)^{2r}(\cos(\theta))^{2r-2s}(\sin(\theta))^{2s}\right) \right\} .
\end{array}
\end{equation}
}
The numerical phase velocity, here expressed as $v_{FD}$,  is obtained as  $v_{FD}=\omega_{FD}/k$. This can be compared to the analytic phase velocity, for which $v=\omega /k = c$, so that
\begin{equation}
\label{eqn:disp_rel}
\delta = \displaystyle\frac{v_{FD}}{v} = \frac{\omega_{FD}}{c\,k},
\end{equation}
provides an error metric of phase velocity error. As discussed in \cite{chen2021framework}, if $\delta<1$, spatial numerical dispersion error appears in the numerical solution, whereas if $\delta>1$, temporal dispersion error appears in the numerical solution. This happens due to the $\delta$ value depending, in general, on the value $\displaystyle\frac{1}{C}$, as we will see in the following sections. 

\subsubsection{Dispersion-based Taylor truncation schemes (DispTE)}
Equation \eqref{eqn:full_disp}, or its version in terms of polar coordinates (relative to equation \eqref{eqn:full_disp_omegaFDtheta}), is usually the basis for the formulation of schemes that aim at reducing dispersion errors. \citet{finkelstein2009spectral} provides details on a general formulation for such schemes, termed the spectral order of accuracy. Alongside, and apparently independently, \citet{liu2009new} proposes a similar, yet more specific, approach. We will here highlight some main characteristics of the methodology, and provide the general form for arbitrary stencil shapes, along similar lines as \cite{liu2013time} and \cite{wang2016effective}.

Assuming $\tau$ and $h$ are small, and using the analytic dispersion relation, $\omega=c k $, along with the Courant number, $C=c\tau/h$, it is reasonable to expand the $\cos(\omega \tau)=\cos(c k \tau)=\cos(C h k)$ in equation \eqref{eqn:full_disp} in its Taylor series around  (for $\tau$ and $h$ sufficiently small). As in \cite{wang2016effective}, we can additionally expand the binomial $k^{2j}=(k_x^2+k_z^2)^j$, that leads to
\begin{equation}
\label{eqn:full_disp_spec_taylor_cos_omega}
\begin{array}{lcc}
\displaystyle \sum_{r=1}^{M}(-1)^r\frac{C^{2r}(h)^{2r}}{(2r)!}\sum_{s=0}^{r} \binom{r}{s}k_x^{2r-2s}k_z^{2s} =  \displaystyle\frac{C^2}{2}a_{0,0}  \\ \\ + 
\displaystyle\frac{C^2}{2}\displaystyle\sum_{m=1}^{M}a_{m,0}\left(\sum_{r=0}^{M}
\frac{2\,m^{2r} h^{2r}}{(2r)!} (-1)^r (k_x^{2r}+k_z^{2r}) \right)  +  \\ \\ 
\displaystyle\frac{C^2}{2}\displaystyle\sum_{q=1}^{Q}\sum_{p=P_0(q)}^{P(q)}a_{p,q}\left(\sum_{r=0}^{M}\sum_{s=0}^{r}\left[
\binom{2r}{2s}
\frac{4\,p^{2r-2s} q^{2s} h^{2r}}{(2r)!} (-1)^r (k_x^{2r-2s}k_z^{2s}) \right]\right) \\ \\ + C^2 O(h^{2M+2})+O(\tau^4),
\end{array}
\end{equation}
where we note that the Taylor expansion of the $\cos(C h k)$ up to  $M$ terms ensures an overall order agreement up to $O(h^{2M})$, as also shown in the appendix of \cite{wang2016effective}. However, we highlight that an $O(\tau^4)$ is still present due to the temporal discretization. We can only avoid this by choosing coefficients ($a_{m,0}$, $a_{p,q}$) such that temporal errors are exactly cancelled out under the assumption of exact dispersion relation representation (trading of temporal derivatives by spatial derivatives).
Comparing the terms of the same polynomial degrees on $hk_x$ and $hk_z$, and using the symmetry of coefficients again, we arrive at the following linear system of equations,
\begin{equation}
\label{eqn:taylor_system_r0_disp}
\begin{array}{lcc}
\displaystyle a_{0,0} + 4\displaystyle\sum_{m=1}^{M}a_{m,0} + \displaystyle 4 \sum_{q=1}^{Q}\sum_{p=P_0(q)}^{P(q)}a_{p,q} = 0, &\quad& (r=0),
\end{array}
\end{equation}
\begin{equation}
\label{eqn:taylor_system_r1_disp}
\begin{array}{lcc}
\displaystyle \sum_{m=1}^{M} m^{2} a_{m,0}
   + \displaystyle  2 \sum_{q=1}^{Q}\sum_{p=P_0(q)}^{P(q)}\,p^{2}a_{p,q} = 1 ,&\quad & (r=1),
\end{array}
\end{equation}
\begin{equation}
\label{eqn:taylor_system_r_s0_disp}
\begin{array}{lcc}
\displaystyle \sum_{m=1}^{M} m^{2r} a_{m,0}
   + \displaystyle  2 \sum_{q=1}^{Q}\sum_{p=P_0(q)}^{P(q)}\,p^{2r}a_{p,q} = C^{2r-2} ,&\quad & (r>1, s=0, r),
\end{array}
\end{equation}
\begin{equation}
\label{eqn:taylor_system_r_s_disp}
\begin{array}{lcc}
\displaystyle\displaystyle  \sum_{q=1}^{Q}\sum_{p=P_0(q)}^{P(q)} \left(
p^{2r-2s} q^{2s} \right) a_{p,q} = \frac{1}{2}\frac{r!(2s)!(2r-2s)!}{(2r)!(r-s)!s!} C^{2r-2} , \quad (1<r \leq M, 0<s \leq r/2)
\end{array}
\end{equation}
where we note that this is the same matrix as shown in \textit{SpatTE} in equations \eqref{eqn:taylor_system_r0}-\eqref{eqn:taylor_system_r_s}, but now the system has a different right-hand-side, with non-zero terms in the mixed equations \eqref{eqn:taylor_system_r_s0_disp} and \eqref{eqn:taylor_system_r_s_disp}. Interestingly, in this case, the non-zero right-hand side of the mixed term sub-system of $a_{p,q}$ allows for non-zero coefficients for $p>0$ and $q>0$. Depending on the choices of $Q$, $P_0$, and $P$, the system may be determined (non-singular) and has a unique non-zero solution. 
Similarly, as with the \textit{SpecTE} scheme, here we can also use polar coordinates $(k, \theta)$ and compare the terms of the same polynomial degrees on $hk$, to arrive at the following linear system of equations,
\begin{equation}
\label{eqn:taylor_system_r0_disp_theta}
\begin{array}{lcc}
\displaystyle a_{0,0} + 4\displaystyle\sum_{m=1}^{M}a_{m,0} + \displaystyle 4 \sum_{q=1}^{Q}\sum_{p=P_0(q)}^{P(q)}a_{p,q} = 0, &\quad& (r=0),
\end{array}
\end{equation}
\begin{equation}
\label{eqn:taylor_system_r1_disp_theta}
\begin{array}{lcc}
\displaystyle \sum_{m=1}^{M} m^{2} a_{m,0}
   + \displaystyle  2 \sum_{q=1}^{Q}\sum_{p=P_0(q)}^{P(q)}\,\psi_{1}(\theta) a_{p,q} = 1 ,&\quad & (r=1),
\end{array}
\end{equation}
\begin{equation}
\label{eqn:taylor_system_r_s0_disp_theta}
\begin{array}{lcc}
\displaystyle \sum_{m=1}^{M} m^{2r} (\cos^{2r}(\theta)+\sin^{2r}(\theta) )a_{m,0} +  \displaystyle  2 \sum_{q=1}^{Q}\sum_{p=P_0(q)}^{P(q)}\,\psi_{r}(\theta) a_{p,q} = C^{2r-2} ,& & (r>1),
\end{array}
\end{equation}
where 
\begin{equation}
     \psi_{r}(\theta) = \sum_{s=0}^{r}
\binom{2r}{2s}
(p^{2r-2s} q^{2s} )(\cos(\theta))^{2r-2s}(\sin(\theta))^{2s}.
 \end{equation}
Here we remark that the system is essentially the same as the one obtained with the spatial-only spectral Taylor truncation (\eqref{eqn:taylor_system_r0_spec}-\eqref{eqn:taylor_system_r_s0_spec}), except that equation \eqref{eqn:taylor_system_r_s0_spec} now has a non-zero right-hand-side in the respective equation \eqref{eqn:taylor_system_r_s0_disp_theta}.
Finally, we note that there are two possibilities with \textit{DispTE}: (i) using complete expansion 2D wavenumbers ($k_x, k_z$) and solving the system \eqref{eqn:taylor_system_r0_disp}-\eqref{eqn:taylor_system_r_s_disp} and capturing high order accuracy in all directions, or, (ii) using polar coordinates and solving the system \eqref{eqn:taylor_system_r0_disp_theta}-\eqref{eqn:taylor_system_r_s0_disp_theta}, but being restricted to special choices of $\theta$ to control accuracy and dispersion errors. We will refer to the former as \textit{DispTE} and to the latter as \textit{DispTE}-$\theta$. Interestingly, both \textit{DispTE} and \textit{DispTE}-$\theta$ are almost identical to \textit{SpatTE}/\textit{SpecTE} and \textit{SpecTE}-$\theta$, respectively, apart from the right-hand side of the linear systems, so the \textit{SpecTE} schemes can be computationally written as particular cases of the \textit{DispTE} schemes (choosing $C=0$).

\subsubsection{Dispersion-based Least-Squares schemes (DispLS)}
\citet{holberg1987computational} proposes that the group velocity error $\left(\displaystyle\frac{\partial \omega}{\partial k}\right)$ should be minimized for the elastic wave equations to obtain a finite difference scheme with adequate reduction of dispersion error. There, a Least Squares approach is suggested, minimizing the dispersion error in a specific bandwidth of wavenumbers. For exact acoustic waves, group and phase velocities are identical, and a similar Least Squares minimization of dispersion errors is possible, as suggested in \cite{liu2013globally}. We will describe here a general formulation for such a case.
The formulation is very similar to the one proposed in Section \ref{sec:spec_LS} for the \textit{SpecLS} method, except that now we consider the full time-space spectral representation of the wave equation from \eqref{eqn:full_disp} using polar coordinates in spectral space,
\begin{equation}
\label{eqn:full_disp_spec_LS}
\begin{array}{lcc}
1-\cos(C \beta)  =  
\displaystyle \frac{C^2} {2} \displaystyle\sum_{m=1}^{M}a_{m,0} \varphi_m(\beta, \theta) +  
\displaystyle \frac{C^2 }{2} \displaystyle\sum_{q=1}^{Q}\sum_{p=P_0(q)}^{P(q)}a_{p,q}\varphi_{p,q}(\beta, \theta) \\ \\ + O(h^{2M}) + O(\tau^4),
\end{array}
\end{equation}
and we have used the previously defined $\varphi_m$ and $\varphi_{p,q}$ from \eqref{eqn:taylor_laplacian_spec_disp_phim} and \eqref{eqn:taylor_laplacian_spec_disp_phipq1} and written the relation in terms of $\beta = k h$.
As in Section \ref{sec:spec_LS}, the Least Squares problem can now be formulated as minimizing the error,
\begin{equation}
    \label{eqn:LS_disp_error}
    \displaystyle E_{\text{disp}}(b) = \int_{0}^{b}\int_{0}^{2\pi} \left( \sum_{m=1}^{M}a_{m,0}\varphi_m(\beta, \theta) + \sum_{q=1}^{Q}\sum_{p=P_0(q)}^{P(q)}a_{p,q} \varphi_{p,q}(\beta, \theta) - \phi_{2}(\beta) \right)^2 d \theta \, d\beta ,
\end{equation}
where
\begin{equation}
\label{eqn:phibeta2}
\phi_{2}(\beta)= \frac{2-2\cos(C \beta)}{C^2}.
\end{equation}
We have again assumed that the minimization with respect to $\beta$ can be focused on a specific range of the spectrum, limited to $\beta \in [0, b] \subset [0, \pi]$, which better captures lower wavenumbers.
The unique Least Squares solution of the problem is obtained by solving the linear system,
\begin{equation}
\label{eqn:taylor_laplacian_spec_disp_linsys1}
\begin{array}{lcc}
\displaystyle
\sum_{m=1}^{M}\langle \varphi_m, \varphi_{m'} \rangle a_{m,0}+ \sum_{q=1}^{Q}\sum_{p=P_0(q)}^{P(q)} \langle \varphi_{p,q}, \varphi_{m'} \rangle a_{p,q} = \langle \phi_{2}(\beta), \varphi_{m'} \rangle, 
 \\
\displaystyle
\sum_{m=1}^{M}\langle \varphi_m, \varphi_{p', q'} \rangle a_{m,0}+ \sum_{q=1}^{Q}\sum_{p=P_0(q)}^{P(q)} \langle \varphi_{p,q}, \varphi_{p', q'} \rangle a_{p,q} = \langle \phi_{2}(\beta), \varphi_{p', q'} \rangle ,
\end{array}
\end{equation}
for $m' = 1, ..., M$,  $q'=1, ..., Q$, $p'=P_0(q'), ..., P(q')$ and using the inner product $\langle\cdot,\cdot\rangle$ defined in \eqref{eqn:cont_dotprod}. The generated system \eqref{eqn:taylor_laplacian_spec_disp_linsys1} is similar to \eqref{eqn:taylor_laplacian_spec_disp_linsys}, except for the right side, where now we have a different $\phi(\beta)$. Similarly to the SpecLS case, we can start with the equation (\ref{eqn:spacexz_derivative2_sym}) and apply $u(x,z,t)= e^{i(k_{x}x+k_{z}z-\omega\tau)}$, where $\omega$ is the angular frequency. In this way, we obtain:
\begin{equation}
\label{eqn:full_disp_spec_LS_cont}
\begin{array}{lcc}
1-\cos(C \beta) \approx  
\displaystyle \frac{C^2} {2}\displaystyle\sum_{m=1}^{M}a_{m,0} \Phi_m(\beta, \theta) +  
\displaystyle \frac{C^2} {2}\displaystyle\sum_{q=1}^{Q}\sum_{p=P_0(q)}^{P(q)}a_{p,q}\Phi_{p,q}(\beta, \theta) \\ \\ + O(h^{2M}) + O(\tau^4).
\end{array}
\end{equation}
When we use the relative dispersion error \eqref{eqn:LS_spec_disp_error} we obtain a very similar Least Square system as in \eqref{eqn:specls_cont_linsys}, given by:
\begin{equation}
\label{eqn:displs_cont_linsys}
\begin{array}{lcc}
\displaystyle
\sum_{m=1}^{M}\langle \Phi_m, \Phi_{m'} \rangle a_{m,0}+ \sum_{q=1}^{Q}\sum_{p=P_0(q)}^{P(q)} \langle \Phi_{p,q}, \Phi_{m'} \rangle a_{p,q} = \langle 1, \Phi_{m'} \rangle, 
 \\
\displaystyle
\sum_{m=1}^{M}\langle \Phi_m, \Phi_{p', q'} \rangle a_{m,0}+ \sum_{q=1}^{Q}\sum_{p=P_0(q)}^{P(q)} \langle \Phi_{p,q}, \Phi_{p', q'} \rangle a_{p,q} = \langle 1, \Phi_{p', q'} \rangle ,
\end{array}
\end{equation}
Where now we define the functions $\Phi_{m}(\beta,\theta)$ and $\Phi_{p,q}(\beta,\theta)$ as:
\begin{equation}
\label{eqn:phim_continuous2}
\displaystyle \Phi_{m} (\beta, \theta)=
\displaystyle \frac{4-2\cos(m\beta\cos(\theta))-2\cos(m\beta\sin(\theta))}{\Psi_{2}(\beta)},
\end{equation}
and
\begin{equation}
\label{eqn:phipq_continuous2}
\displaystyle \Phi_{p,q} (\beta, \theta)= \frac{4-2\cos(q\beta\cos(\theta)+p\beta\sin(\theta))-2\cos(q\beta\cos(\theta)-p\beta\sin(\theta))}{\Psi_{2}(\beta)},
\end{equation}
using the following denominator:
\begin{equation}
\label{eqn:psibeta2}
\Psi_{2}(\beta) = \displaystyle\frac{2-\cos(C\beta)}{C^2}.
\end{equation}
As we discussed in the SpecLS case, the systems \eqref{eqn:displs_cont_linsys} and \eqref{eqn:taylor_laplacian_spec_disp_linsys1} were obtained using full time-space spectral representation. However, \eqref{eqn:taylor_laplacian_spec_disp_linsys1} is obtained starting from \eqref{eqn:taylor_laplacian_spec_disp_linsys}, while \eqref{eqn:displs_cont_linsys} is obtained from \eqref{eqn:spacexz_derivative2_sym}. The explanation about the origin and developments to obtain \eqref{eqn:displs_cont_linsys} and \eqref{eqn:taylor_laplacian_spec_disp_linsys1} is very similar to the SpecLS case, observing now that for \eqref{eqn:displs_cont_linsys} we need to consider the time dependence. Again, for coherence with existing literature, we will adopt \eqref{eqn:displs_cont_linsys} to define the DispLS scheme, but one should be aware that \eqref{eqn:taylor_laplacian_spec_disp_linsys1} can produce different results in certain numerical problems.

\section{Geometrical shape of the spatial stencils}
We will focus on approximations of the Laplacian operator (equation (\ref{eqn:spacexz_derivative2_sym})) for different values of $M$, $P_0$, and $P$. Combinations of these values give the total number of grid points that we will use to approximate the Laplacian operator. 
In general, the integer $M$ controls the number of weights over the axes, while $P$ and $Q$ control the ``inner'' points, with weights not on the axes, which forms the ``extra'', or ``non-cross'', points that can we add in the stencil usually to increase the accuracy of the approximation and/or add other desired properties. 
The geometric shape is generally connected with the spatial order of approximation, but implies other properties of the scheme, including computational aspects, such as memory access/storage of data. Usually, with the implied symmetries, given a stencil with $M$ to-be-defined coefficients over an axis results in a stencil providing accuracy order of order $2M$, provided $0 \leq P, Q\leq M$.

\subsection{Basic shapes}
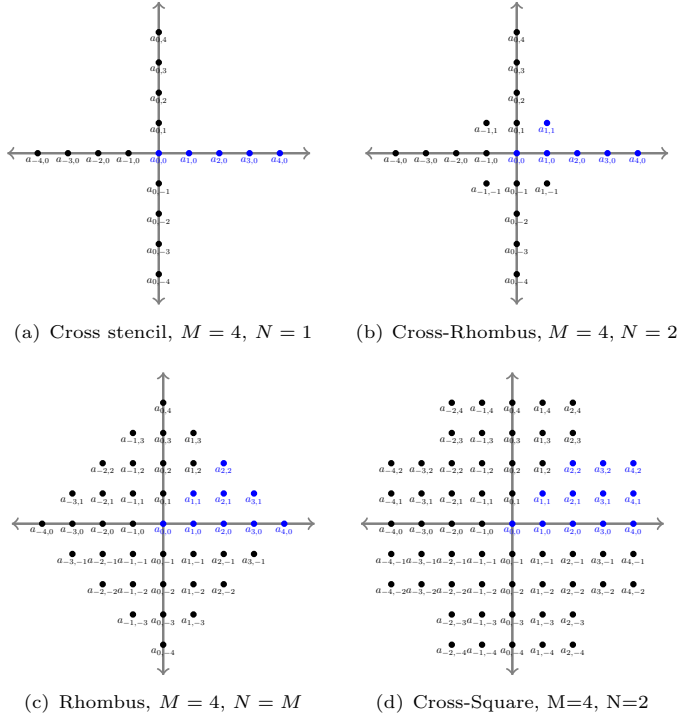
\begin{figure}[H]
\centering

\subfigure[Cross stencil, $M=4$, $N=1$]{
\begin{tikzpicture}[thick,scale=0.4, every node/.style={scale=0.4}]
\draw[gray, thick] (-5,0) -- (5,0);
\draw[gray, thick] (0,-5) -- (0,5);
\draw[gray, <->] (-5,0) -- (5,0);
\draw[gray, <->] (0,-5) -- (0,5);

\filldraw[black] (-4,0) circle (2pt) node[anchor=north]{$a_{-4,0}$};
\filldraw[black] (-3,0) circle (2pt) node[anchor=north]{$a_{-3,0}$};
\filldraw[black] (-2,0) circle (2pt) node[anchor=north]{$a_{-2,0}$};
\filldraw[black] (-1,0) circle (2pt) node[anchor=north]{$a_{-1,0}$};
\filldraw[blue] (0,0) circle (2pt) node[anchor=north]{$a_{0,0}$};
\filldraw[blue] (1,0) circle (2pt) node[anchor=north]{$a_{1,0}$};
\filldraw[blue] (2,0) circle (2pt) node[anchor=north]{$a_{2,0}$};
\filldraw[blue] (3,0) circle (2pt) node[anchor=north]{$a_{3,0}$};
\filldraw[blue] (4,0) circle (2pt) node[anchor=north]{$a_{4,0}$};

\filldraw[black] (0,1) circle (2pt) node[anchor=north]{$a_{0,1}$};
\filldraw[black] (0,2) circle (2pt) node[anchor=north]{$a_{0,2}$};
\filldraw[black] (0,3) circle (2pt) node[anchor=north]{$a_{0,3}$};
\filldraw[black] (0,4) circle (2pt) node[anchor=north]{$a_{0,4}$};

\filldraw[black] (0,-1) circle (2pt) node[anchor=north]{$a_{0,-1}$};
\filldraw[black] (0,-2) circle (2pt) node[anchor=north]{$a_{0,-2}$};
\filldraw[black] (0,-3) circle (2pt) node[anchor=north]{$a_{0,-3}$};
\filldraw[black] (0,-4) circle (2pt) node[anchor=north]{$a_{0,-4}$};

\end{tikzpicture}
\label{fig:cross_line_stencil}
}
~
\subfigure[Cross-Rhombus, $M=4$, $N=2$]{
\begin{tikzpicture}[thick,scale=0.4, every node/.style={scale=0.4}]
\draw[gray, thick] (-5,0) -- (5,0);
\draw[gray, thick] (0,-5) -- (0,5);
\draw[gray, <->] (-5,0) -- (5,0);
\draw[gray, <->] (0,-5) -- (0,5);

\filldraw[black] (-4,0) circle (2pt) node[anchor=north]{$a_{-4,0}$};
\filldraw[black] (-3,0) circle (2pt) node[anchor=north]{$a_{-3,0}$};
\filldraw[black] (-2,0) circle (2pt) node[anchor=north]{$a_{-2,0}$};
\filldraw[black] (-1,0) circle (2pt) node[anchor=north]{$a_{-1,0}$};
\filldraw[blue] (0,0) circle (2pt) node[anchor=north]{$a_{0,0}$};
\filldraw[blue] (1,0) circle (2pt) node[anchor=north]{$a_{1,0}$};
\filldraw[blue] (2,0) circle (2pt) node[anchor=north]{$a_{2,0}$};
\filldraw[blue] (3,0) circle (2pt) node[anchor=north]{$a_{3,0}$};
\filldraw[blue] (4,0) circle (2pt) node[anchor=north]{$a_{4,0}$};

\filldraw[black] (-1,1) circle (2pt) node[anchor=north]{$a_{-1,1}$};
\filldraw[black] (0,1) circle (2pt) node[anchor=north]{$a_{0,1}$};
\filldraw[blue] (1,1) circle (2pt) node[anchor=north]{$a_{1,1}$};

\filldraw[black] (0,2) circle (2pt) node[anchor=north]{$a_{0,2}$};
\filldraw[black] (0,3) circle (2pt) node[anchor=north]{$a_{0,3}$};
\filldraw[black] (0,4) circle (2pt) node[anchor=north]{$a_{0,4}$};

\filldraw[black] (-1,-1) circle (2pt) node[anchor=north]{$a_{-1,-1}$};
\filldraw[black] (0,-1) circle (2pt) node[anchor=north]{$a_{0,-1}$};
\filldraw[black] (1,-1) circle (2pt) node[anchor=north]{$a_{1,-1}$};

\filldraw[black] (0,-2) circle (2pt) node[anchor=north]{$a_{0,-2}$};
\filldraw[black] (0,-3) circle (2pt) node[anchor=north]{$a_{0,-3}$};
\filldraw[black] (0,-4) circle (2pt) node[anchor=north]{$a_{0,-4}$};

\end{tikzpicture}
\label{fig:cross_rhombus_stencil}
}

\subfigure[Rhombus, $M=4$, $N=M$]{
\begin{tikzpicture}[thick,scale=0.4, every node/.style={scale=0.4}]
\draw[gray, thick] (-5,0) -- (5,0);
\draw[gray, thick] (0,-5) -- (0,5);
\draw[gray, <->] (-5,0) -- (5,0);
\draw[gray, <->] (0,-5) -- (0,5);

\filldraw[black] (-4,0) circle (2pt) node[anchor=north]{$a_{-4,0}$};
\filldraw[black] (-3,0) circle (2pt) node[anchor=north]{$a_{-3,0}$};
\filldraw[black] (-2,0) circle (2pt) node[anchor=north]{$a_{-2,0}$};
\filldraw[black] (-1,0) circle (2pt) node[anchor=north]{$a_{-1,0}$};
\filldraw[blue] (0,0) circle (2pt) node[anchor=north]{$a_{0,0}$};
\filldraw[blue] (1,0) circle (2pt) node[anchor=north]{$a_{1,0}$};
\filldraw[blue] (2,0) circle (2pt) node[anchor=north]{$a_{2,0}$};
\filldraw[blue] (3,0) circle (2pt) node[anchor=north]{$a_{3,0}$};
\filldraw[blue] (4,0) circle (2pt) node[anchor=north]{$a_{4,0}$};

\filldraw[black] (-3,1) circle (2pt) node[anchor=north]{$a_{-3,1}$};
\filldraw[black] (-2,1) circle (2pt) node[anchor=north]{$a_{-2,1}$};
\filldraw[black] (-1,1) circle (2pt) node[anchor=north]{$a_{-1,1}$};
\filldraw[black] (0,1) circle (2pt) node[anchor=north]{$a_{0,1}$};
\filldraw[blue] (1,1) circle (2pt) node[anchor=north]{$a_{1,1}$};
\filldraw[blue] (2,1) circle (2pt) node[anchor=north]{$a_{2,1}$};
\filldraw[blue] (3,1) circle (2pt) node[anchor=north]{$a_{3,1}$};

\filldraw[black] (-2,2) circle (2pt) node[anchor=north]{$a_{-2,2}$};
\filldraw[black] (-1,2) circle (2pt) node[anchor=north]{$a_{-1,2}$};
\filldraw[black] (0,2) circle (2pt) node[anchor=north]{$a_{0,2}$};
\filldraw[black] (1,2) circle (2pt) node[anchor=north]{$a_{1,2}$};
\filldraw[blue] (2,2) circle (2pt) node[anchor=north]{$a_{2,2}$};

\filldraw[black] (-1,3) circle (2pt) node[anchor=north]{$a_{-1,3}$};
\filldraw[black] (0,3) circle (2pt) node[anchor=north]{$a_{0,3}$};
\filldraw[black] (1,3) circle (2pt) node[anchor=north]{$a_{1,3}$};

\filldraw[black] (0,4) circle (2pt) node[anchor=north]{$a_{0,4}$};

\filldraw[black] (-3,-1) circle (2pt) node[anchor=north]{$a_{-3,-1}$};
\filldraw[black] (-2,-1) circle (2pt) node[anchor=north]{$a_{-2,-1}$};
\filldraw[black] (-1,-1) circle (2pt) node[anchor=north]{$a_{-1,-1}$};
\filldraw[black] (0,-1) circle (2pt) node[anchor=north]{$a_{0,-1}$};
\filldraw[black] (1,-1) circle (2pt) node[anchor=north]{$a_{1,-1}$};
\filldraw[black] (2,-1) circle (2pt) node[anchor=north]{$a_{2,-1}$};
\filldraw[black] (3,-1) circle (2pt) node[anchor=north]{$a_{3,-1}$};

\filldraw[black] (-2,-2) circle (2pt) node[anchor=north]{$a_{-2,-2}$};
\filldraw[black] (-1,-2) circle (2pt) node[anchor=north]{$a_{-1,-2}$};
\filldraw[black] (0,-2) circle (2pt) node[anchor=north]{$a_{0,-2}$};
\filldraw[black] (1,-2) circle (2pt) node[anchor=north]{$a_{1,-2}$};
\filldraw[black] (2,-2) circle (2pt) node[anchor=north]{$a_{2,-2}$};

\filldraw[black] (-1,-3) circle (2pt) node[anchor=north]{$a_{-1,-3}$};
\filldraw[black] (0,-3) circle (2pt) node[anchor=north]{$a_{0,-3}$};
\filldraw[black] (1,-3) circle (2pt) node[anchor=north]{$a_{1,-3}$};

\filldraw[black] (0,-4) circle (2pt) node[anchor=north]{$a_{0,-4}$};

\end{tikzpicture}
\label{fig:rhombus_stencil}
}
~
\subfigure[Cross-Square, M=4, N=2]{
\begin{tikzpicture}[thick,scale=0.4, every node/.style={scale=0.4}]
\draw[gray, thick] (-5,0) -- (5,0);
\draw[gray, thick] (0,-5) -- (0,5);
\draw[gray, <->] (-5,0) -- (5,0);
\draw[gray, <->] (0,-5) -- (0,5);

\filldraw[black] (-4,0) circle (2pt) node[anchor=north]{$a_{-4,0}$};
\filldraw[black] (-3,0) circle (2pt) node[anchor=north]{$a_{-3,0}$};
\filldraw[black] (-2,0) circle (2pt) node[anchor=north]{$a_{-2,0}$};
\filldraw[black] (-1,0) circle (2pt) node[anchor=north]{$a_{-1,0}$};
\filldraw[blue] (0,0) circle (2pt) node[anchor=north]{$a_{0,0}$};
\filldraw[blue] (1,0) circle (2pt) node[anchor=north]{$a_{1,0}$};
\filldraw[blue] (2,0) circle (2pt) node[anchor=north]{$a_{2,0}$};
\filldraw[blue] (3,0) circle (2pt) node[anchor=north]{$a_{3,0}$};
\filldraw[blue] (4,0) circle (2pt) node[anchor=north]{$a_{4,0}$};

\filldraw[black] (-4,1) circle (2pt) node[anchor=north]{$a_{-4,1}$};
\filldraw[black] (-3,1) circle (2pt) node[anchor=north]{$a_{-3,1}$};
\filldraw[black] (-2,1) circle (2pt) node[anchor=north]{$a_{-2,1}$};
\filldraw[black] (-1,1) circle (2pt) node[anchor=north]{$a_{-1,1}$};
\filldraw[black] (0,1) circle (2pt) node[anchor=north]{$a_{0,1}$};
\filldraw[blue] (1,1) circle (2pt) node[anchor=north]{$a_{1,1}$};
\filldraw[blue] (2,1) circle (2pt) node[anchor=north]{$a_{2,1}$};
\filldraw[blue] (3,1) circle (2pt) node[anchor=north]{$a_{3,1}$};
\filldraw[blue] (4,1) circle (2pt) node[anchor=north]{$a_{4,1}$};

\filldraw[black] (-4,2) circle (2pt) node[anchor=north]{$a_{-4,2}$};
\filldraw[black] (-3,2) circle (2pt) node[anchor=north]{$a_{-3,2}$};
\filldraw[black] (-2,2) circle (2pt) node[anchor=north]{$a_{-2,2}$};
\filldraw[black] (-1,2) circle (2pt) node[anchor=north]{$a_{-1,2}$};
\filldraw[black] (0,2) circle (2pt) node[anchor=north]{$a_{0,2}$};
\filldraw[black] (1,2) circle (2pt) node[anchor=north]{$a_{1,2}$};
\filldraw[blue] (2,2) circle (2pt) node[anchor=north]{$a_{2,2}$};
\filldraw[blue] (3,2) circle (2pt) node[anchor=north]{$a_{3,2}$};
\filldraw[blue] (4,2) circle (2pt) node[anchor=north]{$a_{4,2}$};

\filldraw[black] (-2,3) circle (2pt) node[anchor=north]{$a_{-2,3}$};
\filldraw[black] (-1,3) circle (2pt) node[anchor=north]{$a_{-1,3}$};
\filldraw[black] (0,3) circle (2pt) node[anchor=north]{$a_{0,3}$};
\filldraw[black] (1,3) circle (2pt) node[anchor=north]{$a_{1,3}$};
\filldraw[black] (2,3) circle (2pt) node[anchor=north]{$a_{2,3}$};

\filldraw[black] (-2,4) circle (2pt) node[anchor=north]{$a_{-2,4}$};
\filldraw[black] (-1,4) circle (2pt) node[anchor=north]{$a_{-1,4}$};
\filldraw[black] (0,4) circle (2pt) node[anchor=north]{$a_{0,4}$};
\filldraw[black] (1,4) circle (2pt) node[anchor=north]{$a_{1,4}$};
\filldraw[black] (2,4) circle (2pt) node[anchor=north]{$a_{2,4}$};

\filldraw[black] (-4,-1) circle (2pt) node[anchor=north]{$a_{-4,-1}$};
\filldraw[black] (-3,-1) circle (2pt) node[anchor=north]{$a_{-3,-1}$};
\filldraw[black] (-2,-1) circle (2pt) node[anchor=north]{$a_{-2,-1}$};
\filldraw[black] (-1,-1) circle (2pt) node[anchor=north]{$a_{-1,-1}$};
\filldraw[black] (0,-1) circle (2pt) node[anchor=north]{$a_{0,-1}$};
\filldraw[black] (1,-1) circle (2pt) node[anchor=north]{$a_{1,-1}$};
\filldraw[black] (2,-1) circle (2pt) node[anchor=north]{$a_{2,-1}$};
\filldraw[black] (3,-1) circle (2pt) node[anchor=north]{$a_{3,-1}$};
\filldraw[black] (4,-1) circle (2pt) node[anchor=north]{$a_{4,-1}$};

\filldraw[black] (-4,-2) circle (2pt) node[anchor=north]{$a_{-4,-2}$};
\filldraw[black] (-3,-2) circle (2pt) node[anchor=north]{$a_{-3,-2}$};
\filldraw[black] (-2,-2) circle (2pt) node[anchor=north]{$a_{-2,-2}$};
\filldraw[black] (-1,-2) circle (2pt) node[anchor=north]{$a_{-1,-2}$};
\filldraw[black] (0,-2) circle (2pt) node[anchor=north]{$a_{0,-2}$};
\filldraw[black] (1,-2) circle (2pt) node[anchor=north]{$a_{1,-2}$};
\filldraw[black] (2,-2) circle (2pt) node[anchor=north]{$a_{2,-2}$};
\filldraw[black] (3,-2) circle (2pt) node[anchor=north]{$a_{3,-2}$};
\filldraw[black] (4,-2) circle (2pt) node[anchor=north]{$a_{4,-2}$};

\filldraw[black] (-2,-3) circle (2pt) node[anchor=north]{$a_{-2,-3}$};
\filldraw[black] (-1,-3) circle (2pt) node[anchor=north]{$a_{-1,-3}$};
\filldraw[black] (0,-3) circle (2pt) node[anchor=north]{$a_{0,-3}$};
\filldraw[black] (1,-3) circle (2pt) node[anchor=north]{$a_{1,-3}$};
\filldraw[black] (2,-3) circle (2pt) node[anchor=north]{$a_{2,-3}$};

\filldraw[black] (-2,-4) circle (2pt) node[anchor=north]{$a_{-2,-4}$};
\filldraw[black] (-1,-4) circle (2pt) node[anchor=north]{$a_{-1,-4}$};
\filldraw[black] (0,-4) circle (2pt) node[anchor=north]{$a_{0,-4}$};
\filldraw[black] (1,-4) circle (2pt) node[anchor=north]{$a_{1,-4}$};
\filldraw[black] (2,-4) circle (2pt) node[anchor=north]{$a_{2,-4}$};
\end{tikzpicture}
\label{fig:cross_extra_stencil}
}
\caption{Example of Geometric Representation of Stencils. The blue points indicate the effective weights required for the stencil construction, and all other weights (points) are obtained through symmetry of the blue points. }
\label{fig:stencils}
\end{figure}

Working with one-dimensional stencils to separately approximate the second partial derivatives in each direction is denoted as \emph{line stencil}.
Naturally, bi-dimensional stencils can have more elaborate shapes. Using two one-dimensional \emph{line} stencils gives rise to what is denoted as a \emph{cross-line} stencil, or simply \textit{cross} stencil, as illustrated in Fig. \ref{fig:cross_line_stencil}.
In this case, based on \eqref{eqn:spacexz_derivative2_sym}, we have that $P=Q=0$, that is, we don't have weights out of the axes (no ``extra'' or ``inner'' points).
Due to the symmetry impositions, for the \emph{cross} stencil we just need to calculate $M+1$ coefficients, $a_{m,0}$, for $m=0,1,...,M$. 

Other stencil shapes are also of interest. For example,  \emph{cross-rhombus} and \emph{rhombus} stencils, represented in  Figures \ref{fig:cross_rhombus_stencil} and \ref{fig:rhombus_stencil}, respectively (e.g. \cite{wang2016effective}). These stencils use the same quantities of grid points over the axis as in the cross-line stencil, but also add more points off the axes, with specific rules. 

The full \textit{rhombus} stencil is obtained by choosing $Q=\lfloor M/2 \rfloor$, $P_0(q)=q$ and $P(q)=M-q$. As can be noted in Figure \ref{fig:rhombus_stencil}, due to the symmetry impositions, the number of effective weights to be calculated is given by the ``line'' terms, $M+1$, plus the off-axis terms, which vary depending on $M$. The number of non-cross weights in this rhombus case is given by $\lfloor M/2 \rfloor \lceil M/2 \rceil $.

The \textit{cross-rhombus} stencil is an intermediate stencil between the \textit{rhombus} and the \textit{cross} ones, simply assuming that the rhombus shape is restricted to an inner portion of the stencil only. Effectively, this is obtained by setting $Q=\rfloor N/2 \lfloor $, $P_0(q)=q$ and $P(q)=N-q$, where $0 < N\leq M$ defines the rhombus size, ranging from the cross stencil ($N=1$), to the full rhombus case ($N=M$). Figure \ref{fig:cross_rhombus_stencil} illustrates the case where $M=4$ and $N=2$.
The number of weights that need to be calculated for the \textit{cross-rhombus} stencil depends on the parity of $M$, as with the rhombus stencil, but now we need to consider the value of $N$. If $M$ is even, then the number of weights is $\displaystyle\frac{N^{2}+4M+4}{4}$ and, if $M$ is odd then, the number of weights is $\displaystyle\frac{N^{2}+4M+3}{4}$. Note that in the cross-rhombus stencil, the integers $M$ and $N$ control the geometry and consequently the number of weights, similar to the rhombus stencil.

A \textit{square} stencil can be obtained by choosing $Q=M$, $P_0(q)=q$ and $P(q)=M$, which results in $\displaystyle\frac{M(M+1)}{2}$ weights to be calculated. As with the \textit{rhombus}, it is possible to define an intermediate stencil for the \textit{square} stencil, defining a \textit{cross-square} stencil. Here, we adopt the case with $Q=N \leq M$, $P_0(q)=q$ and $P(q)= M$, which leads to a ``thick'' cross like stencil, which we denote as \textit{cross-square} stencil and show an example in Figure \ref{fig:cross_extra_stencil} for the case with $M=4$ and $N=2$. For $N=0$, one recovers the \textit{cross} stencil, whereas for $N=M$, one recovers the full \textit{square} stencil. The overall number of weights required in this case is given by $\displaystyle\frac{(N+1)(2M-N+2)}{2}$.

Observe that in all the previous schemes we strongly use the symmetrical assumption, that is, $a_{p,q} = a_{q,p}$, when $p,q\geq1$ and $p\neq q$. This fact reduces the number of weights that need to be calculated and directly affects the values of $Q, P(q)$ and $P_{0}(q)$, giving a reduced structure in linear systems.

On the other hand, the system for coefficients $a_{q,p}$ needs to be considered alongside the system of coefficients $a_{p,q}$ to maintain the symmetrical assumption. Without loss of generality, one may define and solve for
%
the sum of these coefficients
$
c_{p,q} = a_{p,q} + a_{q,p}, \mbox { when }\, p,q\geq 1 \mbox{ and } p\neq q. 
$
This method updates the $a_{p,q}$ coefficient and maintains the number of variables and equations described above, following the symmetrical assumption.  

\subsection{Shape and scheme combinations of interest}
Different shapes of stencils (Cross, Cross-Rhombus, Rhombus, Cross-Square, Square) can be combined with other methods of calculations of the coefficients (SpatTE, SpecTE, SpecLS, DispTE, DispLS). We will discuss here the main possible connections between shapes and methods under investigation, and define the most relevant strategies to be investigated.
\begin{itemize}

    \item \textbf{SpatTE-Cross(M)}: For all stencil shapes, the SpatTE scheme always reduces to a \textit{Cross(M)} stencil, with $M+1$ main equations and coefficients given by the solution of the system \eqref{eqn:taylor_system}. As discussed in section \ref{sec:spatTE}, attempts to add non-cross coefficients tend to result in zero coefficients outside the central cross of the stencil. For practical purposes, the SpatTE-Cross(M) scheme/stencil is equivalent to the SpecTE-Cross(M). As with the SpatTE, the SpecTE schemes always reduce to cross-shaped non-zero coefficients.

    \item \textbf{SpecTE-$\theta$-Cross(M)}: The systems for the SpecTE-$\theta$ schemes have $M+1$ equations, so if $Q>0$ and $P>0$, just one valid solution could be obtained when the mixed coefficients $a_{p,q}=0$, that is, just a cross-stencil can be obtained. In this particular case, following \citep{liu2009new} we choose $\theta=\pi/8$. In particular, when $\theta=0$ we recover the classic \textit{SpatTE}.
    
    \item \textbf{DispTE-CrossRb(M, N)}, \textbf{DispLS-CrossRb(M, N)}, \textbf{DispLS-CrossSq(M, N)}, \textbf{SpecLS-CrossRb(M, N)} and \textbf{SpecLS-CrossSq(M, N)}: The systems for these Optimized FD Schemes have the same number of equations and weights, with a non-zero right-hand side. Even though we are working with a square system, the use of the LS solver for the obtained linear system can generate non-feasible solutions in some cases. In such cases, FD stencils that have quasi-zero non-cross weights or weights that decrease the quality of the numerical solution, even though they are acceptable as the LS minimum norm solution. 
    
    \item \textbf{DispTE-CrossSq(M, N)}: The system for these Optimized FD Schemes has, in general, fewer equations than weights, with a non-zero right side. For this reason, the problem is ill-posed, and attempts to obtain numerical solutions satisfying such equations can lead to numerical instabilities. In practice, we observe that for certain choices of $M$ and $N$ the stencils generated decrease the quality of the numerical solution. 
        
\end{itemize}

In the cases of \textit{DispTE-CrossRb(M, N)}, \textit{DispLS-CrossRb(M, N)}, and \textit{SpecLS-CrossRb (M, N)} choosing $N=1$, we recovered the \textit{Cross(M)} case for each spatial choice. Similarly, choosing $N=0$ for \textit{DispTE-CrossSq(M, N)}, \textit{DispLS-CrossSq(M, N)}, and \textit{SpecLS-CrossSq(M, N)}, we recovered the \textit{Cross(M)} case for each spatial choice. 

It should also be clear that \textit{DispTE-CrossRb(M,1)} and \textit{DispTE-CrossSq(M,0)} are the same as the \textit{DispTE-Cross(M)} scheme, and a similar conclusion can be reached for \textit{DispLS-CrossRb(M,1)} and \textit{DispLS-CrossSq(M,0)}, which are the same as \textit{DispLS-Cross(M)}, and \textit{SpecLS-CrossRb(M,1)} and \textit{SpecLS-CrossSq(M,0)} are the same as \textit{SpecLS-Cross(M)}.

\section{Computational Benchmarks}

\subsection{Software Framework}
The numerical simulations were carried out by using the Devito framework \cite{devito-api,devito-compiler,kukreja2016devito}. Devito is a domain-specific Language (DSL) and code generation framework for developing highly optimized finite difference kernels for use in inversion methods, which are widely used in seismic problems. There are other frameworks with similar functionalities as Devito, for example, DENISE \citep{kohn2012influence} and SAVA \citep{kohn2015waveform}. However, Devito utilizes \textit{SymPy} library to allow the definition of operators from high-level symbolic equations and generates optimized and automatically tuned code specific to a given target computational architecture. Symbolic computation is a powerful tool that allows users to build complex solvers from only a few lines of high-level code. It also helps to build automated performance optimizations for generated codes, and gives flexibility to adjust stencil discretization at run time, making efficient code development much less time-consuming. 

For the implementation of different stencils in Devito, we have replaced the default Finite Difference schemes \cite{fornberg1988generation} with a customized symbolic function based on the schemes defined previously.

In our case, we introduce the new stencils to approximate the Laplacian operator. As these new weights are pre-computed symbolically, they will be automatically incorporated by Devito into the optimized generated code.

Special boundary conditions can be easily added (see \cite{dolci2022effectiveness} and references therein).

As mentioned in the previous sections, we consider a damping layer strategy \citep{sochacki1987absorbing,cerjan1985nonreflecting}, which can reduce reflections at the boundaries, thereby avoiding most of the effects induced by reflected waves.

\subsection{Data}
We highlight three essential data objects frequently used in seismic imaging: the \textbf{snapshots}, the \textbf{seismic traces}, and the \textbf{receiver}. For this propose, let $u_{num}(x_{i},z_{j},t^{k})$ be a generic numerical solution obtained using one of the FD schemes proposed previously (for a suitable spatial and spatial scheme) in a spatial general position $(x_{i},z_{j})$ and a time step $t^{k}$. 

The \textbf{snapshot} is the numerical solution selected in the time step $t^{k}=t^{*}$ and for a restricted set of spatial points, where, in general, we consider the entire discretized spatial domain $\Pi_{xz}$. We will denote a given snapshot in the time step $t^{*}$ over all the points in $\Pi_{xz}$ by $u_{snap}(t^{*})$. More precisely, we have:
$
u_{snap}(t^{*}) = \left\{\left. u_{num}(x_{i},z_{j},t^{*})\right|_{\forall(x_{i},z_{j})\in\Pi_{xz}}\right\}.  
$

The \textbf{seismic trace} is the numerical solution selected for a fixed spatial position $(x_{*},z_{*})$ and a selected time range, where, in general, we consider all the step times in $\Pi_{t}$. We will denote a given seismic trace in the position $(x_{*},z_{*})$ over all the time-steps in $\Pi_{t}$ by $u_{st}(x_{*},z_{*})$. More precisely, we have:
$u_{st}(x_{*},y_{*}) = \left\{\left. u_{num}(x_{*},z_{*},t^{k})\right|_{\forall t^{k}\in\Pi_{t}}\right\}.  $
Finally, the \textbf{receiver} data can be interpreted as a set of seismic traces for a set of positions. In general, the receivers are placed using the same value for $z_{j}$ and varying the $x_{i}$ position, not necessarily with uniform spacing for the selected points $x_{i}$. In terms of time steps, we generally choose all the points in $\Pi_{t}$. Regarding spatial points, we will indicate an arbitrary set $A_{xz}\subset\Pi_{xz}$ of selected points. We will denote the receiver as $u_{rec}$ and define it as:
$u_{rec}(A_{xz}) = \left\{\left.u_{st}(x_{p},z_{p})\right|_{\forall (x_{p},z_{p})\in A_{xz}}\right\}$
Note that a snapshot is fixed in time, but varies in space; a seismic trace varies in time for a fixed spatial position; and a receiver data varies in space and time. This behaviour of the data objects leads us to choose specific tools to evaluate these objects, as we will see in the next section.

\subsection{Evaluation metrics}

The evaluation of seismic wave propagation in the literature employs a heavy use of visual inspection of snapshots, seismic traces, and receiver data (e.g. \citep{finkelstein2009spectral, liu2009new, wang2016effective, chen2019numerical, liu2013globally, wang2016effective}). Here, we will discuss quantitative measures that help to quantify what is usually noted via visual inspection.

We will use classic error norm evaluations and propose error metrics using the Fourier and Wavelet Transform. The main goal is to seek dispersion effects, usually captured by visual inspection, quantitatively in the numerical solutions.

\subsubsection{Classical Evaluation}
The classical evaluation is based on the usual norm of matrices and vectors in $\mathbb{R}^{p}\times\mathbb{R}^{q}$, with $p,q\in\mathbb{N}$. We can see the snapshots/receiver as bi-dimensional matrices and the seismic trace as a one-dimensional matrix, or simply, a vector. In this case, we can use the usual vector norms $\left\|\cdot\right\|_{1}$, $\left\|\cdot\right\|_{2}$, and $\|\cdot\|_{\infty}$.

Independent of the norm type, the error is constructed by the norm of the difference between the numerical and a reference solution. Following the previous notation, we will denote $u_{i,j,k}^{num} = u_{num}(x_{i},y_{j},t^{k})$ and $u_{i,j,k}^{ref} = u_{ref}(x_{i},y_{j},t^{k})$ for an arbitrary triple $(x_{i},y_{j},t^{k})$ in $\Pi_{xz}\times\Pi_{t}$. For convenience, we will use $u_{snap}^{num}(t^{*})$, $u_{st}^{num}(x_{*},z_{*})$ and $u_{rec}^{num}(A_{xz})$ for the data objects of the numerical solution and $u_{snap}^{ref}(t^{*})$, $u_{st}^{ref}(x_{*},z_{*})$ and $u_{rec}^{ref}(A_{xz})$ for the data objects of the numerical reference. Using classic notation for error norms, we define the following error norms for the snapshots:
\begin{equation}
\label{eqn:snapnorm1}
\left\|u_{snap}^{num}(t^{*})-u_{snap}^{ref}(t^{*})\right\|_{1} = \displaystyle\frac{1}{h^{2}}\left(\sum_{i=1}^{n_{x}}\sum_{j=1}^{n_{z}}\vert u_{i,j,*}^{num}-u_{i,j,*}^{ref}\vert\right),
\end{equation}
\begin{equation}
\label{eqn:snapnorm2}
\left\|u_{snap}^{num}(t^{*})-u_{snap}^{ref}(t^{*})\right\|_{2} = \displaystyle\frac{1}{h}\left(\sqrt{\sum_{i=1}^{n_{x}}\sum_{j=1}^{n_{z}}\vert u_{i,j,*}^{num}-u_{i,j,*}^{ref}\vert^{2}}\right),
\end{equation}
\begin{equation}
\label{eqn:snapnormmax}
\left\|u_{snap}^{num}(t^{*})-u_{snap}^{ref}(t^{*})\right\|_{\infty} = \displaystyle\max\limits_{1\leq i\leq n_{x},1\leq j\leq n_{z}}\left(\vert u_{i,j,*}^{num}-u_{i,j,*}^{ref}\vert\right).
\end{equation}
For the seismic trace and receivers, we have analogous relations, changing the space dimension for the time dimension where appropriate, and $h$ for $\tau$ similarly. 

\subsubsection{Fourier and Wavelet Evaluation}
Fourier Transform and Wavelet Transform are essential tools in several contexts, from theoretical to practical results in some areas, particularly, in the seismic application \citep{broughton2018discrete,debnath2015wavelet,grossmann1984decomposition,mallat1999wavelet,meyer1993wavelets,morlet1982wave}. In our case, we will use the Discrete Fourier Transform (DFT) \citep{cooley1965algorithm} and the Continuous Wavelet Transform (CWT) \citep{torrence1998practical}. These transformations have a huge mathematical background, and the main results can be applied from different perspectives.  
These two discrete transformations are applied in functions defined in the time domain; however, the final domain depends on the type of transformation. The DFT results are functions in the frequency domain, while the CWT results are functions in the time-frequency domain.

The connection between DFT/CWT and the FD schemes studied in this work is that the dispersion effect can be better understood by examining numerical and reference solutions in the frequency and/or time-frequency domains, which facilitate the easier determination of wave amplitudes and phases.
These two quantities, amplitude and phase, are numerical indicators of the dissipative and dispersion properties of the numerical and reference solutions. In particby examining the phase of the numerical solution, we can demonstrate how the wave is delayed or advanced relative to se of the reference solution. The amplitude indicates how the energy in the trough and crest of the numerical solution deviates from the reference solution.
In practical terms, we will denote the DFT by $\mathscr{F}(\cdot)$ and CWT by $\mathscr{W}(\cdot)$, which will be used for seismic traces and receiver data, which contain time evolutions. We will describe the process in more detail in the Appendix for obtaining two essential analysis components: $\gamma_{k}^{\mathscr{F}}$ and $\gamma_{n,k}^{\mathscr{W}}$. These values measure the weighted phase difference for each $k$ wavenumber (frequency) in the case of DFT transformation, and the weighted phase difference for each $k$ wavenumber  (frequency) and $n$ time step in the case of CWT transformation. In both cases, these values are calculated at the point $(x_{*},z_{*})$ over a simulated time range, that is, choosing a seismic trace at $(x_{*},z_{*})$ we calculate the respective DFT and CWT transforms and use the respective phase and amplitude to obtain the previous values based on numerical and reference solutions. These values are given by:   
\begin{equation}
\label{eqn:gammadftloc}
\gamma_{k}^{\mathscr{F}} = \displaystyle\left(\frac{\hat{\phi}_{dif}^{\mathscr{F}}(k)}{\pi}\right)\left(\frac{|\Gamma^{num}_{A_{k}}|}{\Gamma_{A}^{num}}\right),
\end{equation}
and
\begin{equation}
\label{eqn:gammacwtloc}
\gamma_{n,k}^{\mathscr{W}} = \displaystyle\left(\frac{\hat{\phi}_{dif}^{\mathscr{W}}(n,k)}{\pi}\right)\left(\frac{|\Gamma^{num}_{S_{n,k}}|}{\Gamma_{S}^{num}}\right),
\end{equation}
where the elements of the previous quantities involve the phase difference and amplitude of numerical and reference solutions, in direction to quantify the dispersion effects for a given seismic trace.   
A complete and detailed description of the previous formulas is in Appendix A.
We can then compute the norms using these quantities to obtain unique error estimates associated with dispersion effects.
In summary, for snapshots/receivers, we will use classical norms to better understand the error in the spatial dimension of the numerical solution and for seismic traces we will employ DFT- and CWT-based norms to understand dispersion effects across spatial and temporal dimensions.

\section{Numerical Experiments}

\subsection{General Considerations}
The goal of the following numerical experiments is to highlight the pros and cons of the proposed stencils relative to classical stencils. Ultimately, we aim to find optimal settings for the use of such stencils.

To ensure reproducibility, we discuss all parameters used in the investigations here. We evaluate $M$ from 1 to 6 for all FD Schemes proposed in this work. That is, we use spatial orders 2-12 for the numerical simulations we present (recall that M is half the accuracy order). For the cross-rhombus and cross-square stencils, we consider the ranges $ 1\leq N\leq M$ and $ 0\leq N\leq M$, respectively, for each $N$. We present results up to order 12 because, in practical applications, this is an admissible choice for the spatial operator. Spatial operators with orders higher than 12 generally do not significantly improve accuracy in realistic settings compared to 8th-12th order methods; however, they are computationally more expensive when applied to scenarios that require solving the Acoustic Wave Equation multiple times, such as in the FWI procedure.  

For simplicity, we assume $\Delta x = \Delta z$, and choose $\tau$ such that we achieve stability considering the maximum medium velocity.

Moreover, we will adopt the same order for each FD stencil used in each test; that is, we calculate the maximum velocity in the velocity field and construct the FD scheme based on this maximum velocity and the discrete parameters of the problem at hand, for the FD schemes that use the maximum velocity in the calculation of the weights. This strategy will reveal possible limitations in the dispersion-reduction effects in velocity fields with high velocity contrast.  

For all tests, a reference solution will be used to compare and analyze the errors in the numerical solutions, adopting spatial/temporal mesh resolutions 16 times smaller than the numerical solution.

\subsection{Homogeneous Velocity Model}
As a first benchmark, we adopt the classical Homogeneous Velocity Model. Table \ref{table1:homogeneo} presents the set of parameters we will use in the test, following a similar homogeneous approach in \citep{wang2016effective}.
\begin{table}
\centering
\caption{Set of Parameters for Homogeneous Velocity Model}
\label{table1:homogeneo}
\begin{tabular}{|c|c|}
\hline
\textbf{Parameters} & \textbf{Values}    \\
\hline
[$x_{I}$, $x_{F}$]     & [$0m$, $6000m$]   \\
\hline
[$z_{I}$, $z_{F}$]    & [$0m$, $6000m$]   \\
\hline
$t_{F}$    & $1500ms$   \\
\hline
$nx=nz$       & $401$   \\
\hline
$\Delta x=\Delta z$ & $15m$   \\
\hline
$\Delta t$ & $2ms$   \\
\hline
$nt$       & $1500$   \\
\hline
$f_{0}$    & $30Hz$   \\
\hline
$c$   & $3$ km/s \\
\hline
\end{tabular}
\end{table}
The constant velocity is assumed to be $c=3$ km/s, and the source position is at the middle of the domain. In this experiment, a Damping strategy was applied to minimize boundary reflections and isolate dispersion phenomena. By selecting a displacement snapshot prior to boundary interaction and employing a homogeneous-velocity model with a single propagation velocity, the analysis focuses primarily on dispersion effects and the performance of optimized finite-difference (FD) schemes.

\begin{figure}[H]
\centering
\includegraphics[scale=.30]{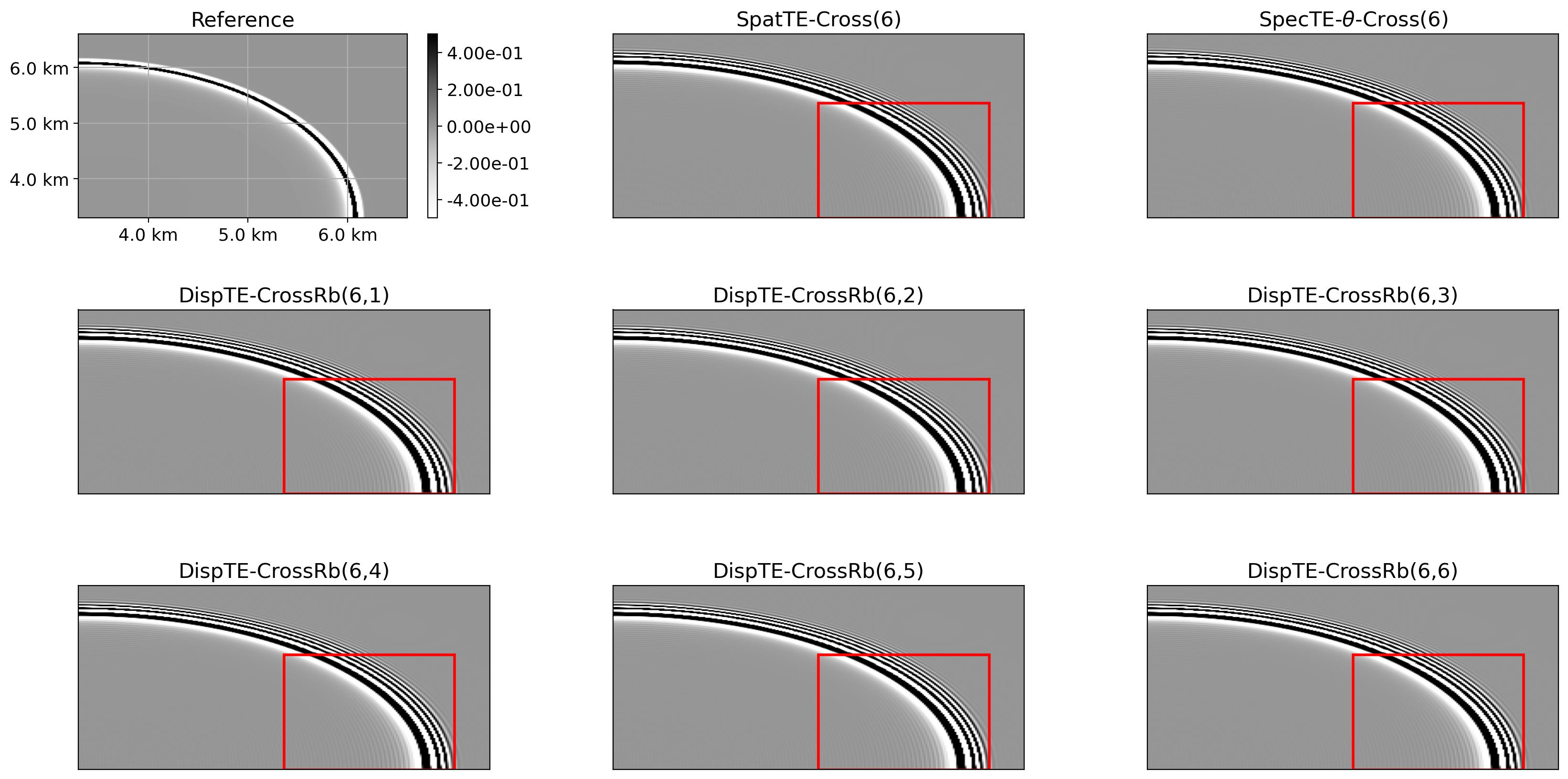} 
\caption{Homogeneous Velocity Model - Displacement - $DispTE-CrossRb(6,N)$ at $T=1.05s$.}
\label{fig:homo1}
\end{figure}

\begin{figure}[H]
\centering
\includegraphics[scale=.30]{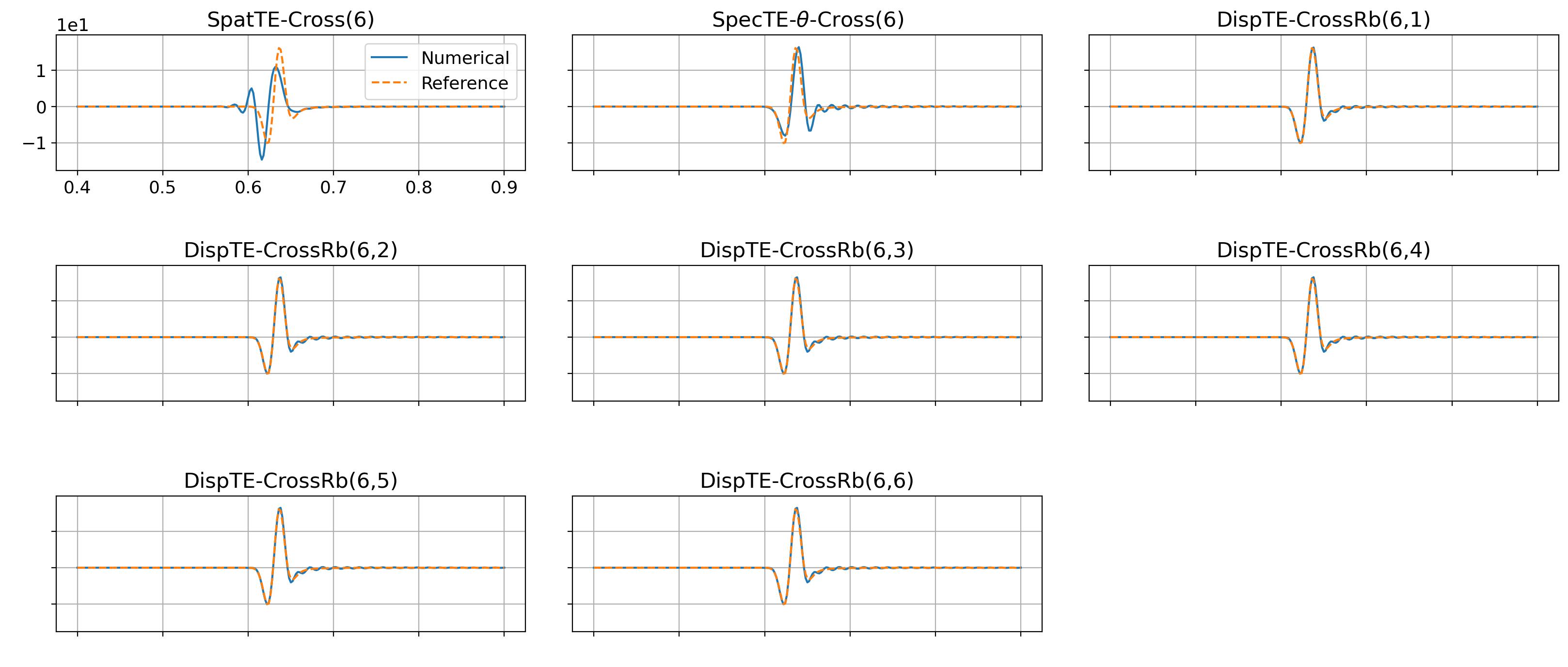} 
\caption{Homogeneous Velocity Model - Seismic Trace - $DispTE-CrossRb(6,N)$ at $x=1200m$ and $z=3000m$.}
\label{fig:homo3}
\end{figure}

\begin{figure}[H]
\centering
\includegraphics[scale=.30]{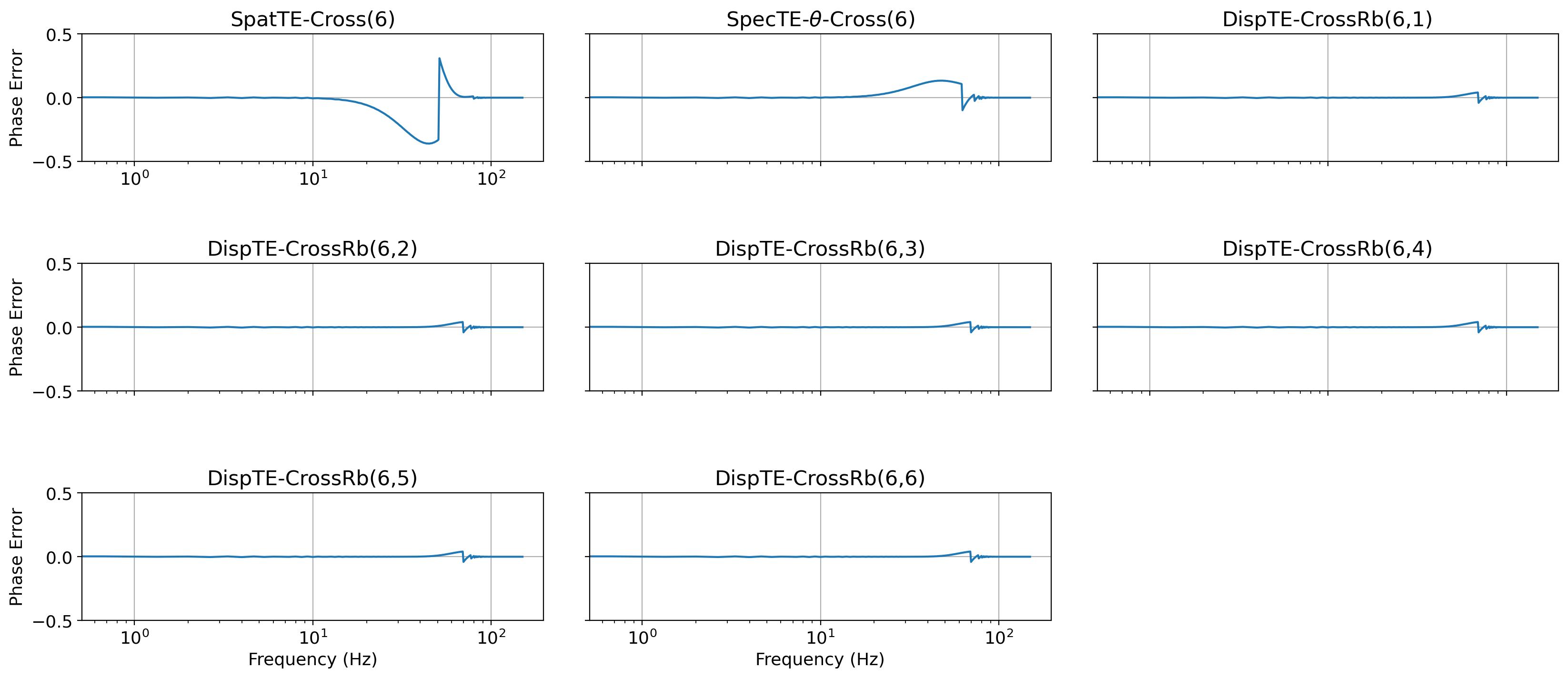} 
\caption{Homogeneous Velocity Model - FFT Analysis of Seismic Trace - $DispTE-CrossRb(6,N)$ at $x=1200m$ and $z=3000m$.}
\label{fig:homo5}
\end{figure}

\begin{figure}[H]
\centering
\includegraphics[scale=.30]{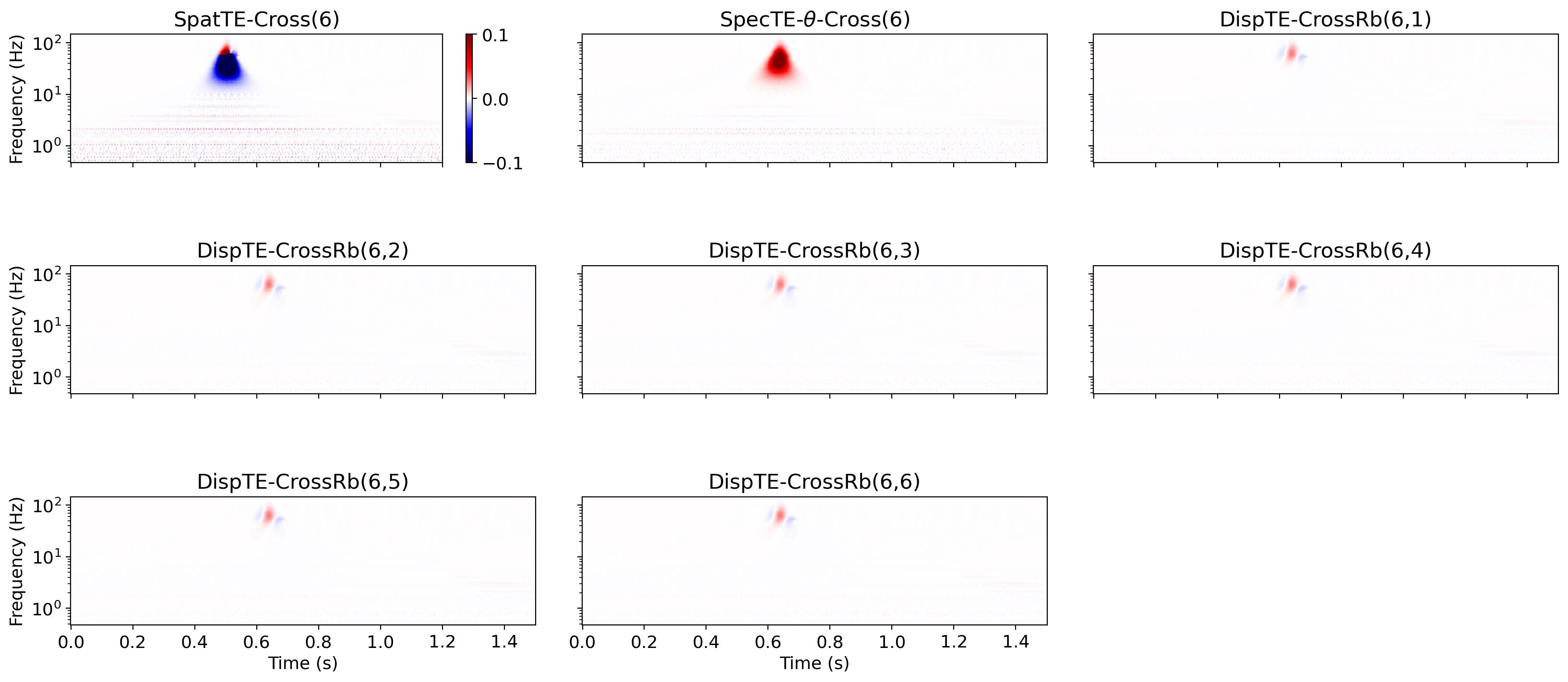} 
\caption{Homogeneous Velocity Model - CWT Analysis of Seismic Trace - $DispTE-CrossRb(6,N)$ at $x=1200m$ and $z=3000m$.}
\label{fig:homo7}
\end{figure}

Figures~\ref{fig:homo1},~\ref{fig:homo3},~\ref{fig:homo5}, and~\ref{fig:homo7} present comparative results for the $SpatTE$-, $SpecTE$-, and $DispTE$-based Cross schemes, evaluating displacement, seismic traces, FFT, and CWT dispersion curves, respectively. From Figure~\ref{fig:homo1}, it is evident that $DispTE$-$CrossRb(6,N)$ exhibits the least dispersion, followed by $SpecTE$-$\theta$-$Cross(6)$ and $SpatTE$-$Cross(6)$. A similar pattern is observed in the seismic traces (Figure~\ref{fig:homo3}), where $DispTE$-$CrossRb(6,N)$ and $SpecTE$-$\theta$-$Cross(6)$ closely match the reference solution, while $SpatTE$-$Cross(6)$ shows noticeable dispersion. The frequency- and time–frequency-domain analyses in Figures~\ref{fig:homo5} and~\ref{fig:homo7} confirm these observations quantitatively, with $SpatTE$-$Cross(6)$ being the most dispersive scheme.

\begin{figure}[H]
\centering
\includegraphics[scale=.30]{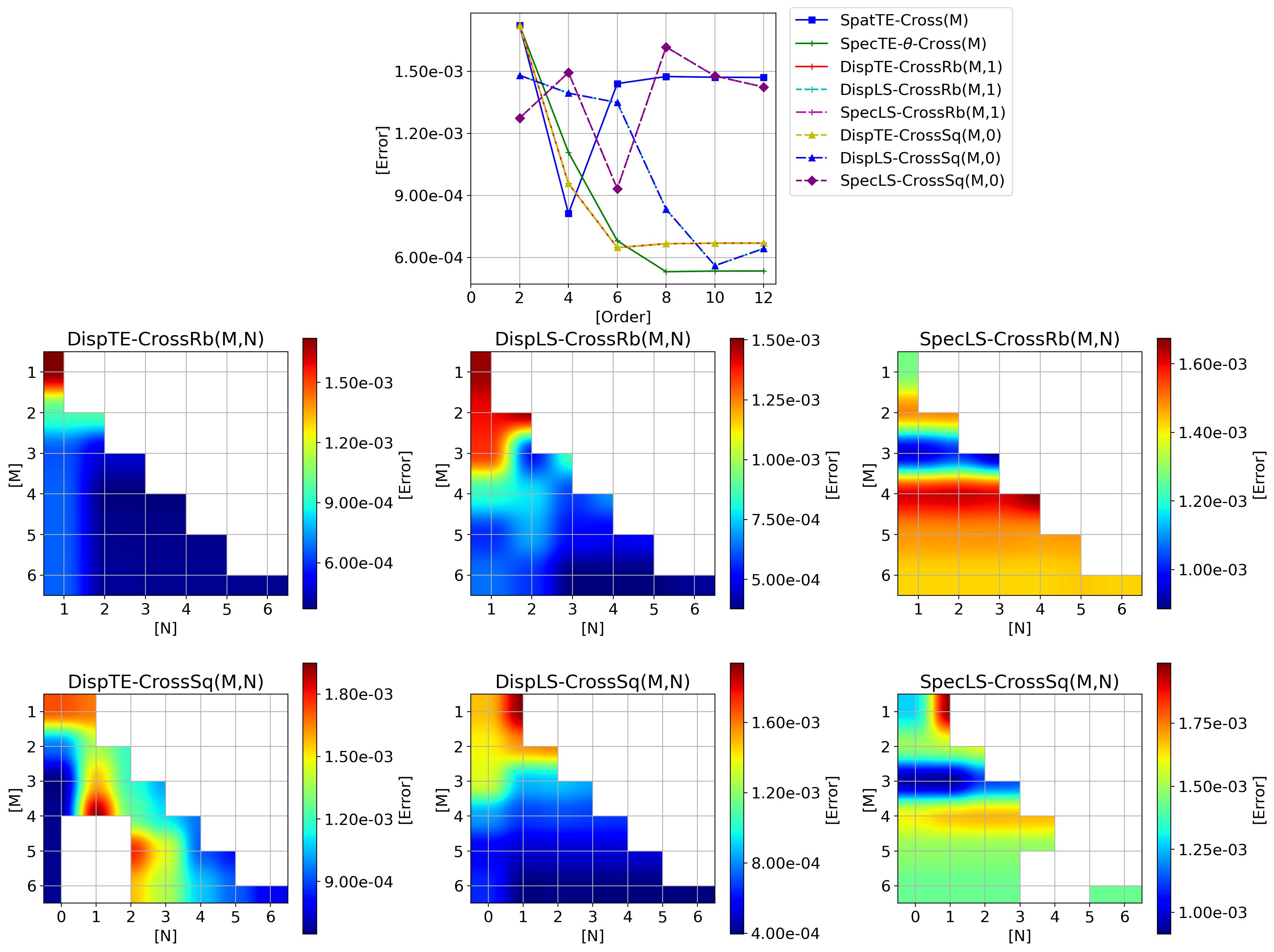} 
\caption{Homogeneous Velocity Model - Displacement - $\|\cdot\|_{2}$-norm at $T=1.05s$.}
\label{fig:homo2}
\end{figure}

\begin{figure}[H]
\centering
\includegraphics[scale=.30]{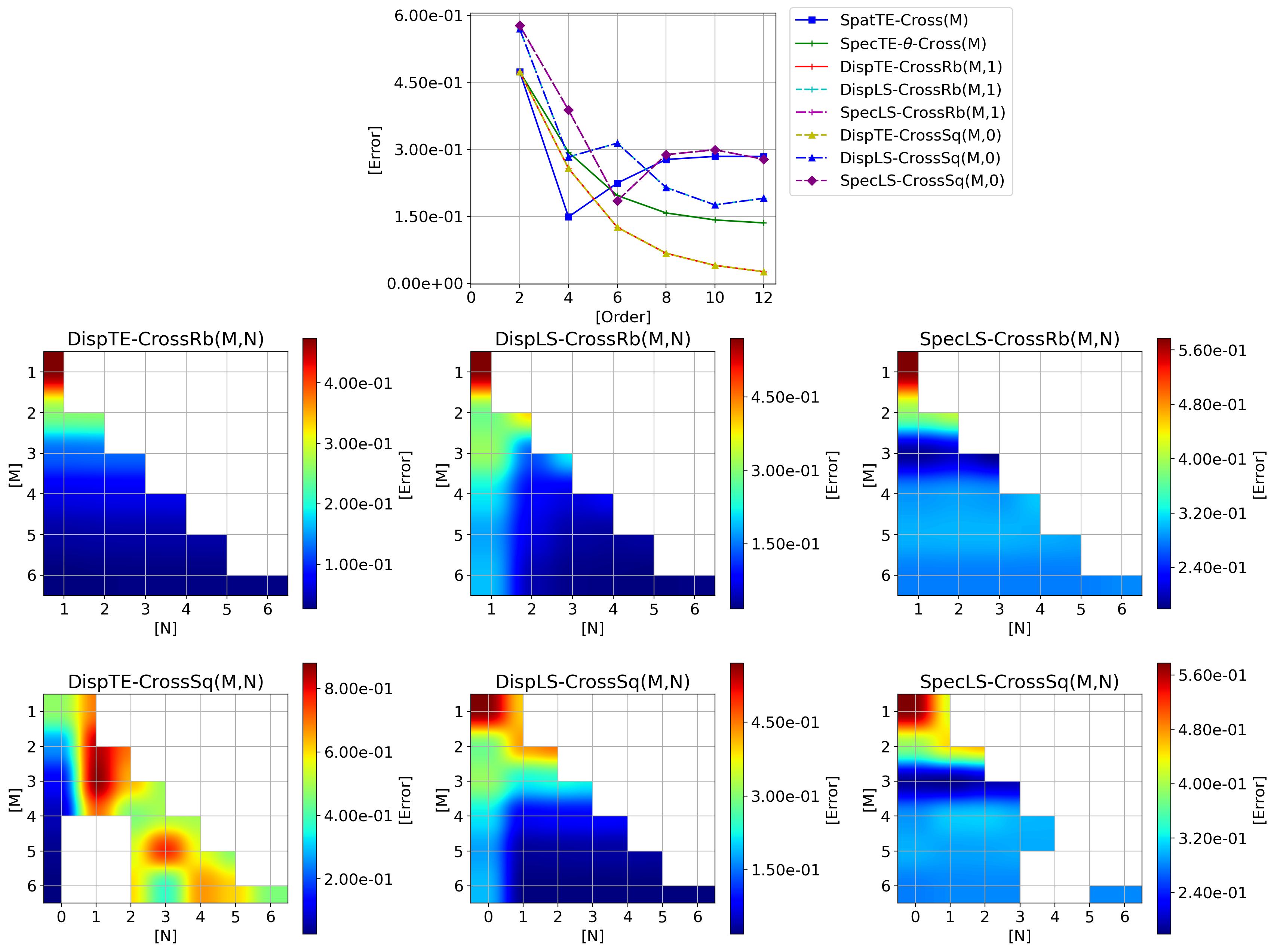} 
\caption{Homogeneous Velocity Model - Seismic Trace - $\|\cdot\|_{2}$-norm at $x=1200m$ and $z=3000m$.}
\label{fig:homo4}
\end{figure}

\begin{figure}[H]
\centering
\includegraphics[scale=.30]{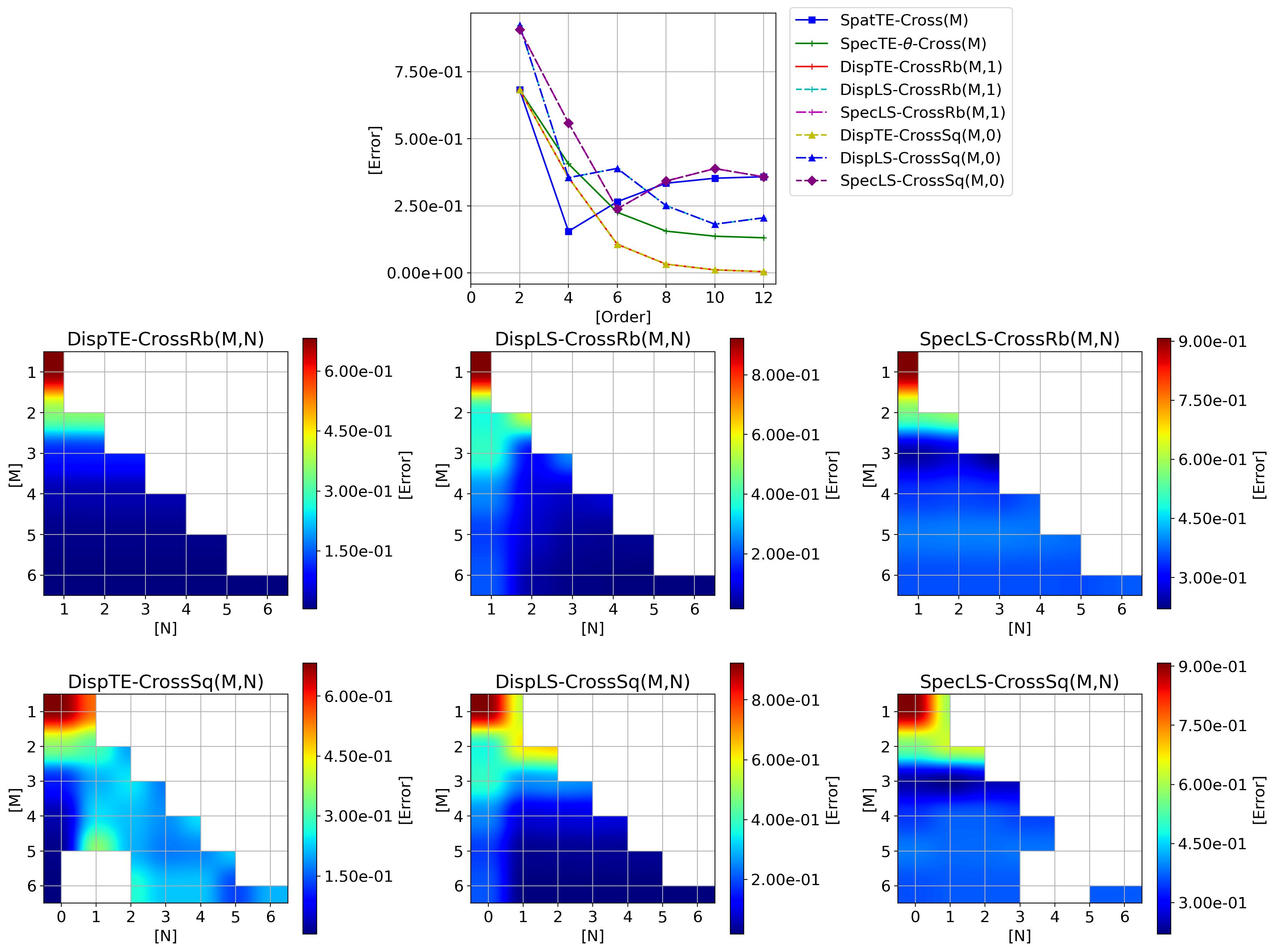} 
\caption{Homogeneous Velocity Model - Seismic Trace - $\|\gamma^{\mathscr{F}}_{u_{st}(x,z)}\|_{2} $-norm at $x=1200m$ and $z=3000m$.}
\label{fig:homo6}
\end{figure}

\begin{figure}[H]
\centering
\includegraphics[scale=.30]{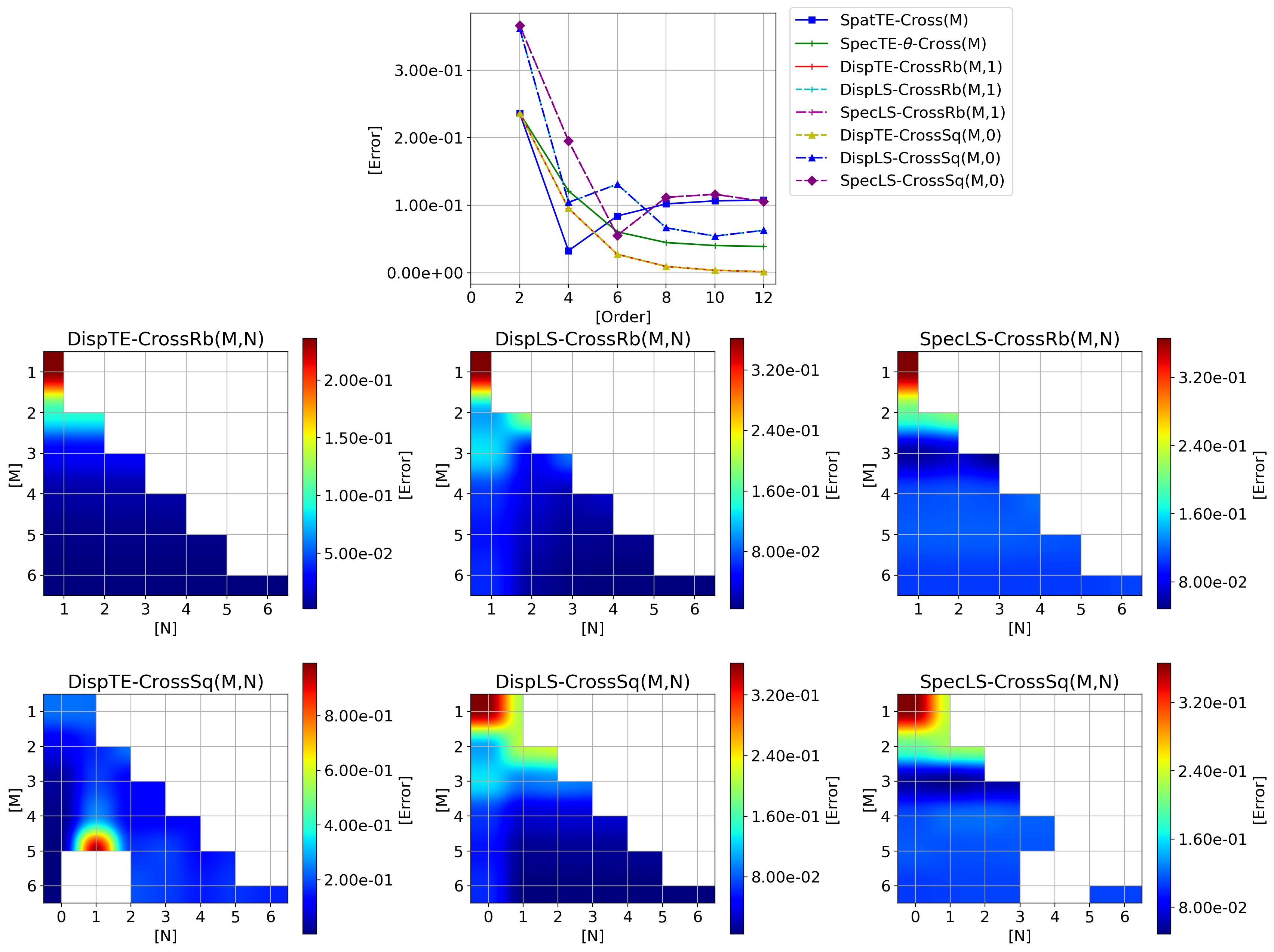} 
\caption{Homogeneous Velocity Model - Seismic Trace - $\|\gamma^{\mathscr{W}}_{u_{st}(x,z)}\|_{2} $-norm at $x=1200m$ and $z=3000m$.}
\label{fig:homo8}
\end{figure}

Norm evaluations in Figures~\ref{fig:homo1}–\ref{fig:homo7} reinforce these trends: $SpatTE$-$Cross(6)$ consistently yields the highest norm values, while optimized FD schemes achieve lower errors. The analysis suggests organizing the schemes into two groups—$Cross$ and $Cross$-$Rb$/$Cross$-$Sq$—to better assess their convergence behavior. This extended comparison is shown in Figures~\ref{fig:homo2},~\ref{fig:homo4},~\ref{fig:homo6}, and~\ref{fig:homo8}, where all FD-optimized configurations are contrasted.

For the $Cross$ configurations ($N=1$ for $Cross$-$Rb$ and $N=0$ for $Cross$-$Sq$), the optimized FD schemes outperform the classical $SpatTE$-$Cross(M)$, with $DispTE$-$CrossRb(M,1)$ generally showing the best results, except in Figure~\ref{fig:homo2}, where $SpecTE$-$\theta$-$Cross(M)$ performs slightly better. For higher-order $Cross$-$Rb$ and $Cross$-$Sq$ cases ($N\geq1$ and $N\geq0$), $DispTE$-$CrossRb(M,N)$, $DispLS$-$CrossRb(M,N)$, and $DispLS$-$CrossSq(M,N)$ display similar trends: increasing $M$ reduces the error, while varying $N$ has a limited effect. Hence, the $Cross$ schemes appear more efficient than their $Cross$-$Rb$ and $Cross$-$Sq$ counterparts.

The $Cross$-$Sq$ schemes present additional challenges. For $DispTE$-$CrossSq(M,N)$, the number of weights can exceed the number of equations for $N\geq1$, producing either unstable or negligible weights when solved using least squares (LS) methods. Consequently, some schemes yield nearly equivalent or even degraded performance compared to the $Cross$ variants. High condition numbers in the system matrices of $DispTE$-$CrossSq(M,N)$, $DispLS$-$CrossSq(M,N)$, and $SpecLS$-$CrossSq(M,N)$ can further compromise numerical stability. Conversely, the $Cross$-$Rb$ configurations remain stable for moderate orders ($M,N \leq 10$), as they form well-posed square systems.

Norm analyses in Figures~\ref{fig:homo2}–\ref{fig:homo8} confirm that the trends observed for $M=6$ extend to lower orders ($1\leq M\leq5$), which are relevant to practical applications. Overall, the homogeneous case demonstrates that optimized FD stencils outperform classical schemes, especially when suitable spatial–temporal parameters are chosen. For small parameter values, however, all stable schemes exhibit comparable low-dispersion behavior.

An important observation is that optimized schemes depend strongly on the velocity used in the stencil optimization. In the homogeneous scenario, a uniform velocity field enables the consistent application of a single optimized stencil, thereby enhancing performance. In heterogeneous media, by contrast, this uniform optimization is impractical and may diminish accuracy. 

In summary, for the homogeneous case, appropriate tuning of spatial–temporal parameters and velocity values enables clear differentiation between optimized and classical FD schemes across displacement, seismic trace, FFT, and CWT analyses, as evidenced in Figures~\ref{fig:homo1}–\ref{fig:homo8}.

\subsection{Two Layers Velocity Model}
To address the variable velocity case, we first adopt a classic Two-Layer Velocity Model test.
The set of parameters for the Heterogeneous Velocity Model test is similar to the Homogeneous Velocity Model test, as we can see in Table \ref{table1:hetero}. Moreover, we adopt the same parameters as those used in \citep{wang2016effective}, where the velocity is $1.5Km/s$ from $z=0$ to $z=1200m$ and increases to $3.0Km/s$, after $z=1200m$, where the source is positioned at $ x=2000 m$ and $ z=1600 m$.

\begin{figure}[H]
\centering
\includegraphics[scale=.30]{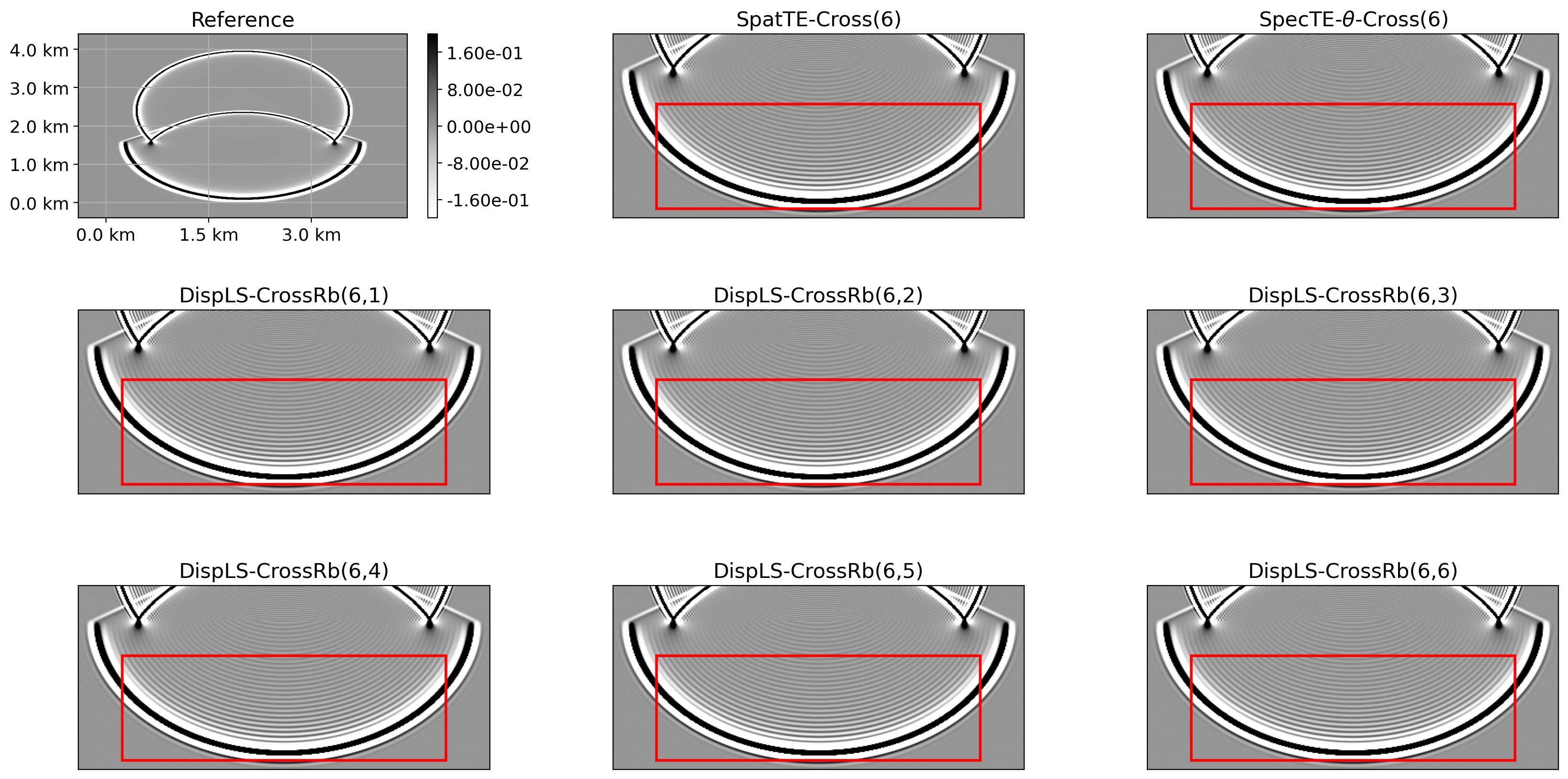} 
\caption{Heterogeneous Velocity Model - Displacement - $DispLS-CrossRb(6,N)$ at $T=1.08s$.}
\label{fig:het1}
\end{figure}

\begin{figure}[H]
\centering
\includegraphics[scale=.30]{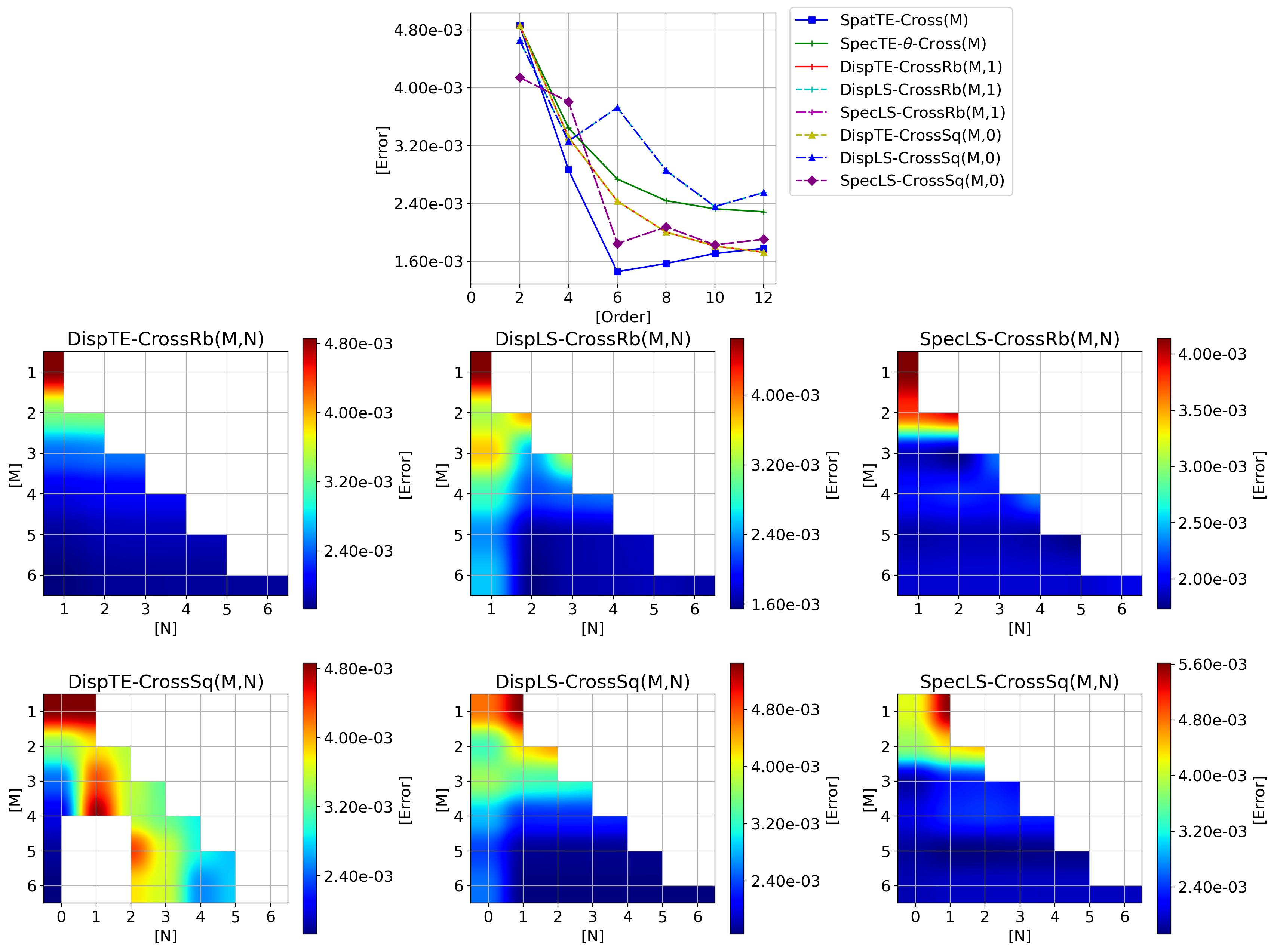} 
\caption{Heterogeneous Velocity Model - Displacement - $\|\cdot\|_{2}$-norm at $T=1.08s$.}
\label{fig:het2}
\end{figure}

\begin{figure}[H]
\centering
\includegraphics[scale=.30]{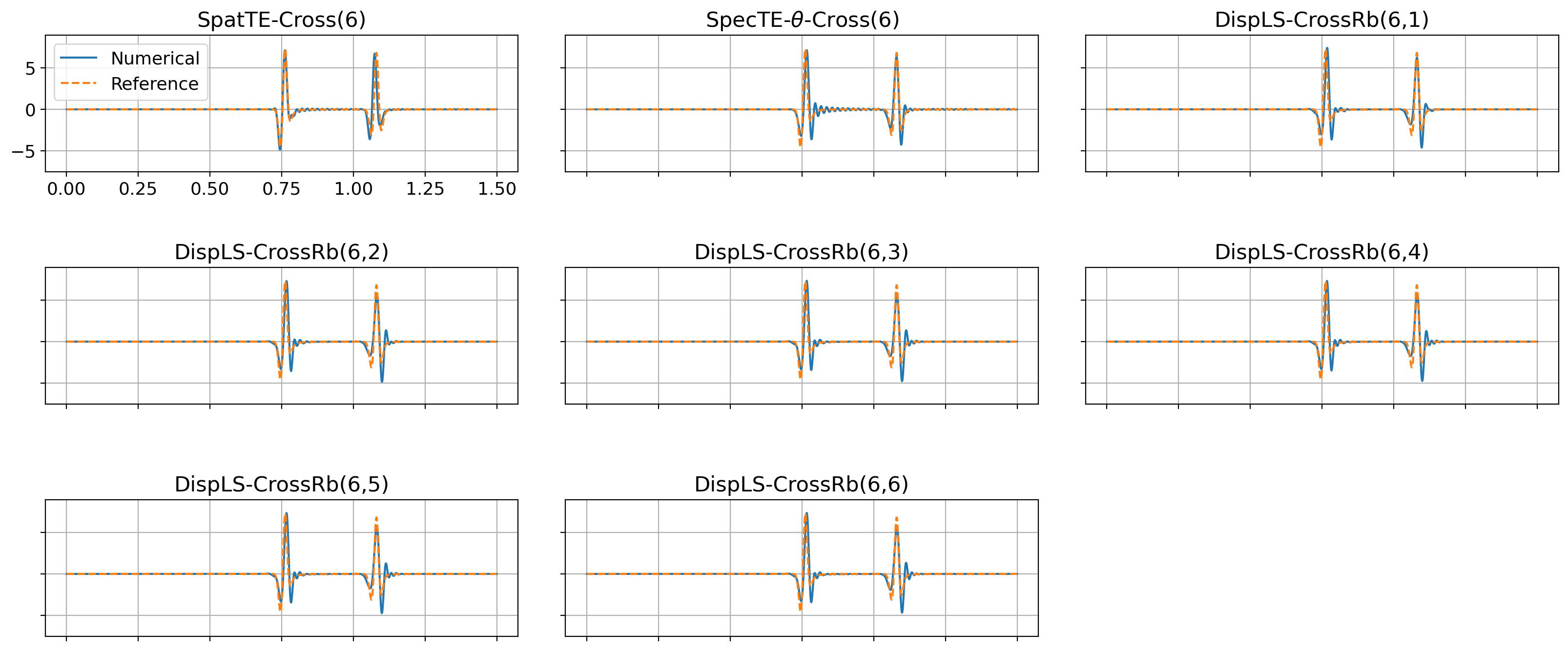} 
\caption{Heterogeneous Velocity Model - Seismic Trace - $DispLS-CrossRb(6,N)$ at $x=1000m$ and $z=2000m$.}
\label{fig:het3}
\end{figure}

\begin{figure}[H]
\centering
\includegraphics[scale=.30]{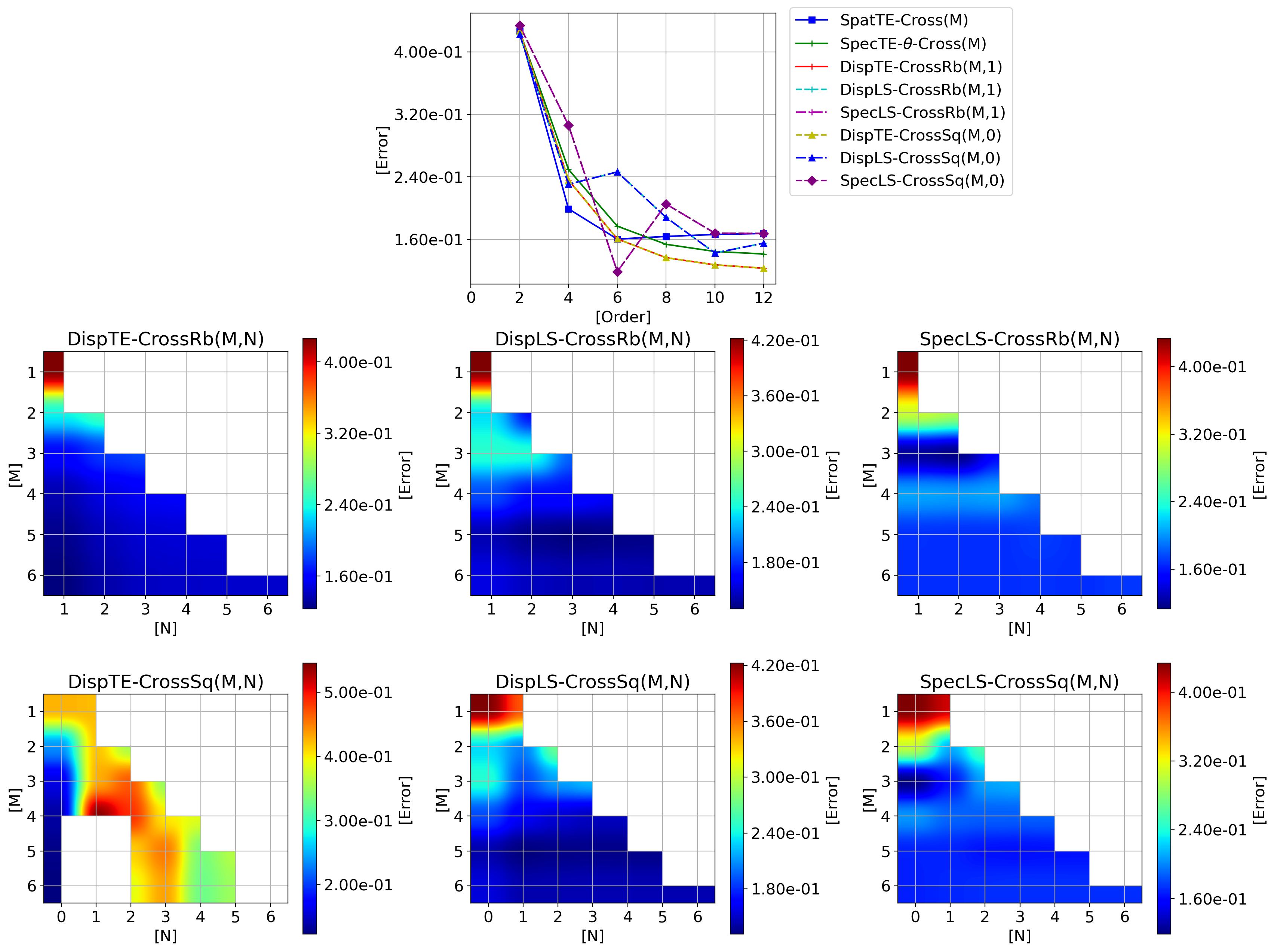} 
\caption{Heterogeneous Velocity Model - Seismic Trace - $\|\cdot\|_{2}$-norm at $x=1000m$ and $z=2000m$.}
\label{fig:het4}
\end{figure}

\begin{figure}[H]
\centering
\includegraphics[scale=.30]{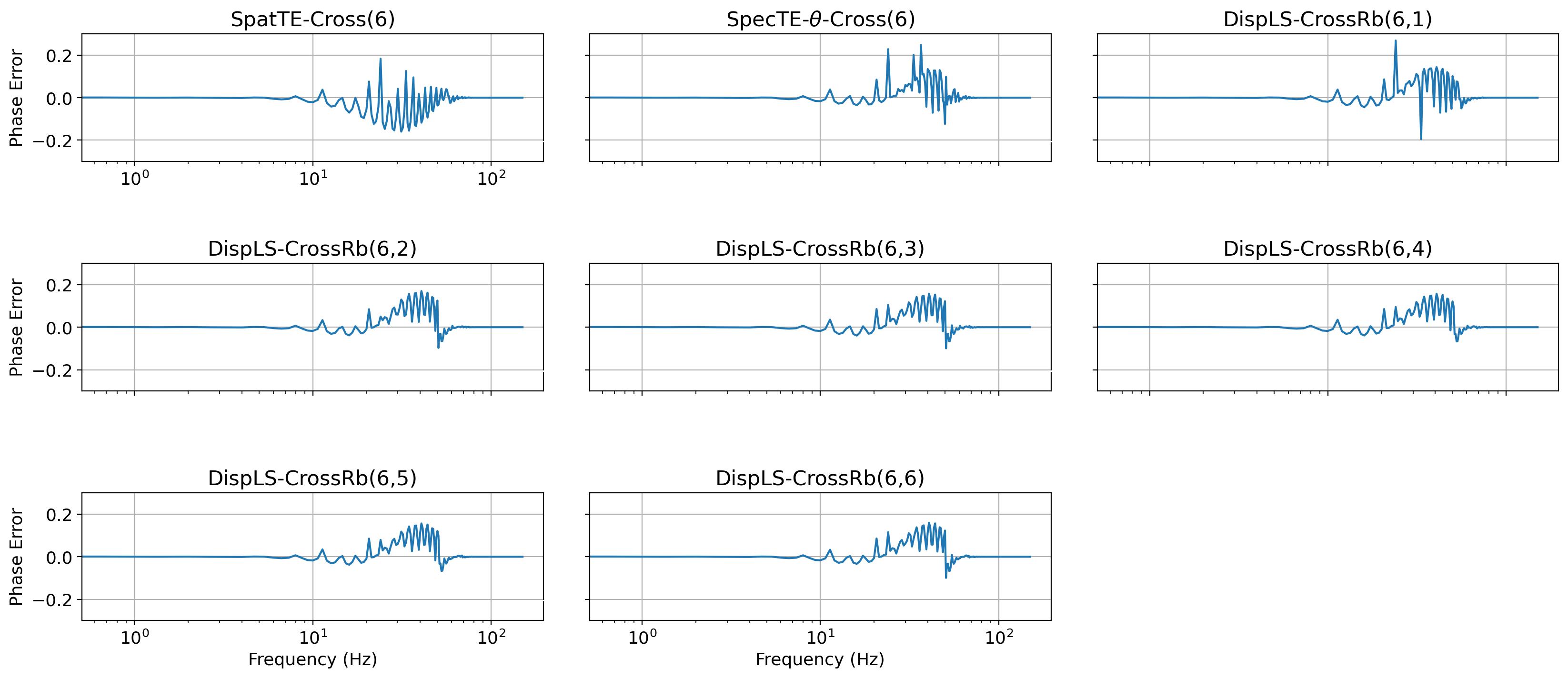} 
\caption{Heterogeneous Velocity Model - FFT Analysis of Seismic Trace - $DispLS-CrossRb(6,N)$ at $x=1000m$ and $z=2000m$.}
\label{fig:het5}
\end{figure}

\begin{figure}[H]
\centering
\includegraphics[scale=.30]{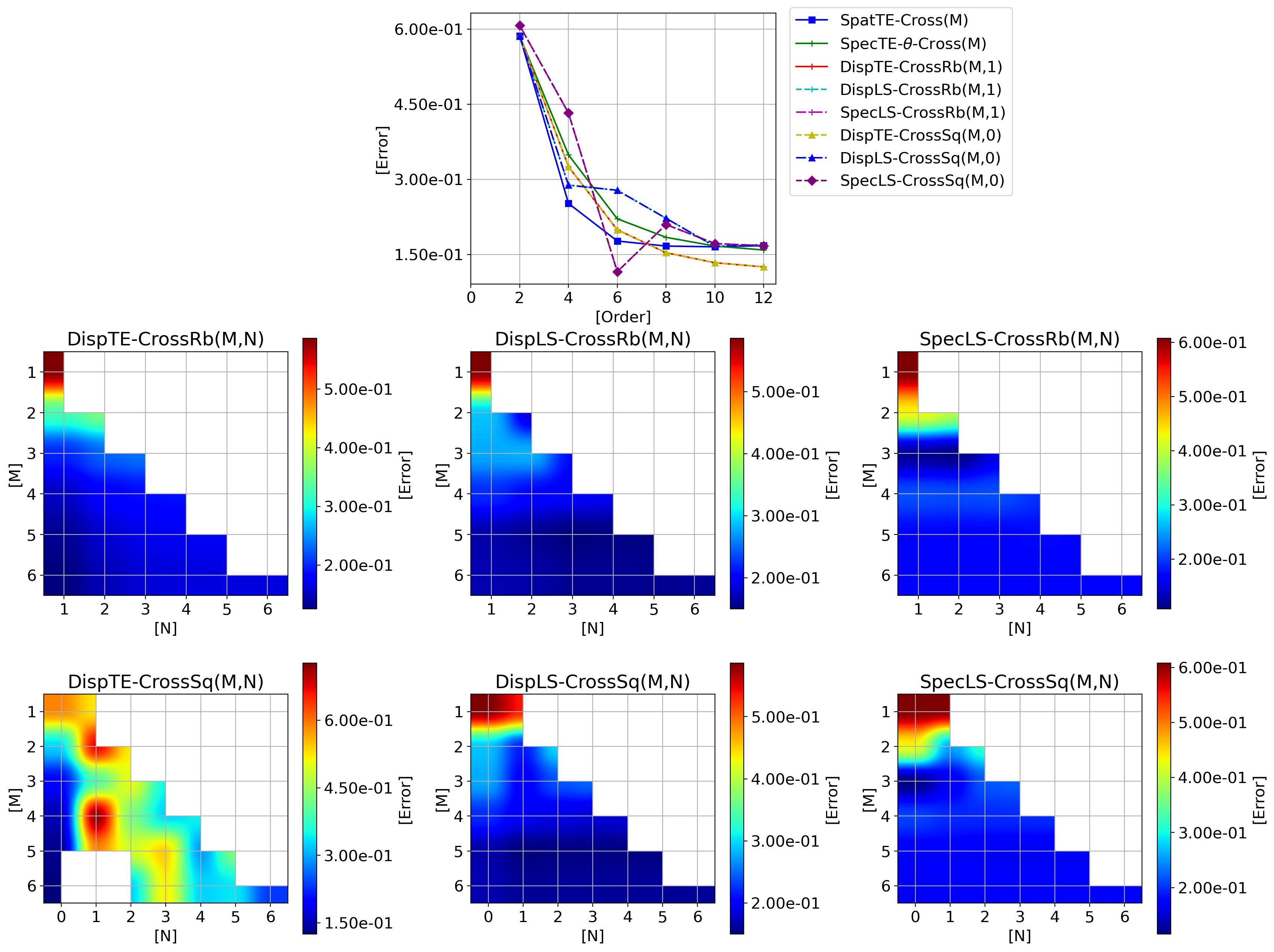} 
\caption{Heterogeneous Velocity Model - Seismic Trace - $\|\gamma^{\mathscr{F}}_{u_{st}(x,z)}\|_{2} $-norm at $x=1000m$ and $z=2000m$.}
\label{fig:het6}
\end{figure}

\begin{table}
\centering
\caption{Set of Parameters for Heterogeneous Velocity Model}
\label{table1:hetero}
\begin{tabular}{|c|c|}
\hline
\textbf{Parameters} & \textbf{Values}    \\
\hline
[$x_{I}$,$x_{F}$] & [$0m$,$4000m$]   \\
\hline
[$z_{I}$,$z_{F}$] & [$0m$,$4000m$]   \\
\hline
$t_{F}$    & $1500s$   \\
\hline
$nx=nz$       & $401$   \\
\hline
$\Delta x$ & $10mm$   \\
\hline
$\Delta x =\Delta z$ & $10m$   \\
\hline
$\Delta t$ & $1.0ms$   \\
\hline
$nt$       & $1500$   \\
\hline
$f_{0}$    & $25Hz$   \\
\hline
\end{tabular}
\end{table}
Figures~\ref{fig:het1},~\ref{fig:het3} and ~\ref{fig:het5}, present comparative results for the $SpatTE$-$Cross(6)$, $SpecTE$-$\theta$-$Cross(6)$, and $DispLS$-$CrossRb(6,N)$ schemes, with $1 \leq N \leq 6$, analyzing the displacement, seismic trace, and FFT dispersion curves, respectively. Figures~\ref{fig:het2},~\ref{fig:het4},~\ref{fig:het6}, display the corresponding norm comparisons for the same configurations.

As in the homogeneous case, a Damping strategy is employed to suppress boundary reflections, and displacement plots are chosen such that waves do not reach the domain boundaries. However, the presence of two propagation velocities introduces additional complexity due to wave interactions across media, potentially influencing numerical results. The optimized FD stencils were constructed using the model's maximum velocity, with the spatial–temporal parameters derived from this value. Ideally, separate stencils should be designed for each velocity region, but for this analysis, a single FD scheme was fixed per configuration, which may affect the accuracy of the results.

Figure~\ref{fig:het1} shows that $SpatTE$-$Cross(6)$, $SpecTE$-$\theta$-$Cross(6)$, and $DispLS$-$CrossRb(6,1)$ exhibit pronounced dispersion artifacts (highlighted in the red bound-box), while increasing $N$ in $DispLS$-$CrossRb(6,N)$ ($N \geq 2$) progressively reduces these effects. The norm analysis in Figure~\ref{fig:het2} reveals distinct trends among the schemes: while the errors of $SpecTE$-$\theta$-$Cross(M)$ and $DispLS$-$CrossRb(M,1)$ decrease with increasing $M$, the $SpatTE$-$Cross(M)$ error grows for $M \geq 3$. Moreover, $DispTE$-$CrossSq(M,N)$ exhibits an instability region, consistent with observations from the homogeneous case.

Although $SpatTE$-$Cross(M)$ errors increase with $M$, the scheme remains competitive—sometimes even outperforming the optimized FD schemes—owing to its independence from the velocity model. Nonetheless, for fixed $M$, varying $N$ yields mixed results: in $DispTE$-$CrossRb(M,N)$ and $SpecLS$-$CrossRb(M,N)$, higher $N$ does not significantly improve the norm, while for $DispLS$-$CrossRb(M,N)$, $DispLS$-$CrossSq(M,N)$, and $SpecLS$-$CrossSq(M,N)$, increases in $N$ can either have negligible or adverse effects, particularly for $M \geq 4$. These deteriorations are attributed to the increasing condition number of the system matrices, in which the least-squares (LS) solver produces non-feasible yet numerically stable solutions.

In Figure~\ref{fig:het3}, the seismic trace reveals that $SpatTE$-$Cross(6)$ lags slightly behind the reference solution, while $SpecTE$-$\theta$-$Cross(6)$ and $DispLS$-$CrossRb(6,N)$ are slightly ahead. Amplitude discrepancies and oscillations in all solutions indicate residual dispersion, though these effects are notably reduced in $DispLS$-$CrossRb(6,N)$. The norm comparison in Figure~\ref{fig:het4} further demonstrates that, for the $Cross$ configurations, $SpatTE$-$Cross(M)$ performs worse than both $SpecTE$-$\theta$-$Cross(M)$ and $DispLS$-$CrossRb(M,0)$. Except for the unstable $DispTE$-$CrossSq(M,N)$ case, the general trend across all FD schemes is that increasing $M$ reduces dispersion, while varying $N$ has little influence, aside from specific anomalies in $DispLS$-$CrossSq(M,N)$ and $SpecLS$-$CrossSq(M,N)$.

The FFT dispersion curves in Figure~\ref{fig:het5} exhibit mild noise across all methods, evidencing residual dispersion effects. Figure~\ref{fig:het6} confirms that $DispLS$-$CrossRb(M,1)$ achieves the best results among the $Cross$ schemes, and that increasing $M$ generally reduces dispersion without strong dependence on $N$, except in the unstable $DispTE$-$CrossSq(M,N)$ configuration. The results using CWT dispersion curves show similar behavior for the test, highlighting that, even with small differences, optimized FD schemes consistently exhibit less dispersion than $SpatTE$-$Cross(M)$. 

Nevertheless, in this two-layer model, demonstrating the benefits of optimized FD schemes requires detailed analysis using FFT dispersion curves and their norm-based quantification, unlike the homogeneous case, where improvements are more visually apparent.

Overall, while the optimized FD schemes outperform the classical $SpatTE$-$Cross(M)$ in dispersion control, their advantages are less pronounced in the heterogeneous model because they require constructing stencils based on the maximum velocity. This choice renders the classical scheme comparatively competitive. Therefore, as in the homogeneous-velocity model, appropriate parameter selection enables optimized FD schemes to produce accurate and stable results, even when using a single stencil based on the maximum velocity.

\subsection{2D SEG EAGE Salt Velocity Model}
\begin{figure}[H]
\centering
\includegraphics[scale=.15]{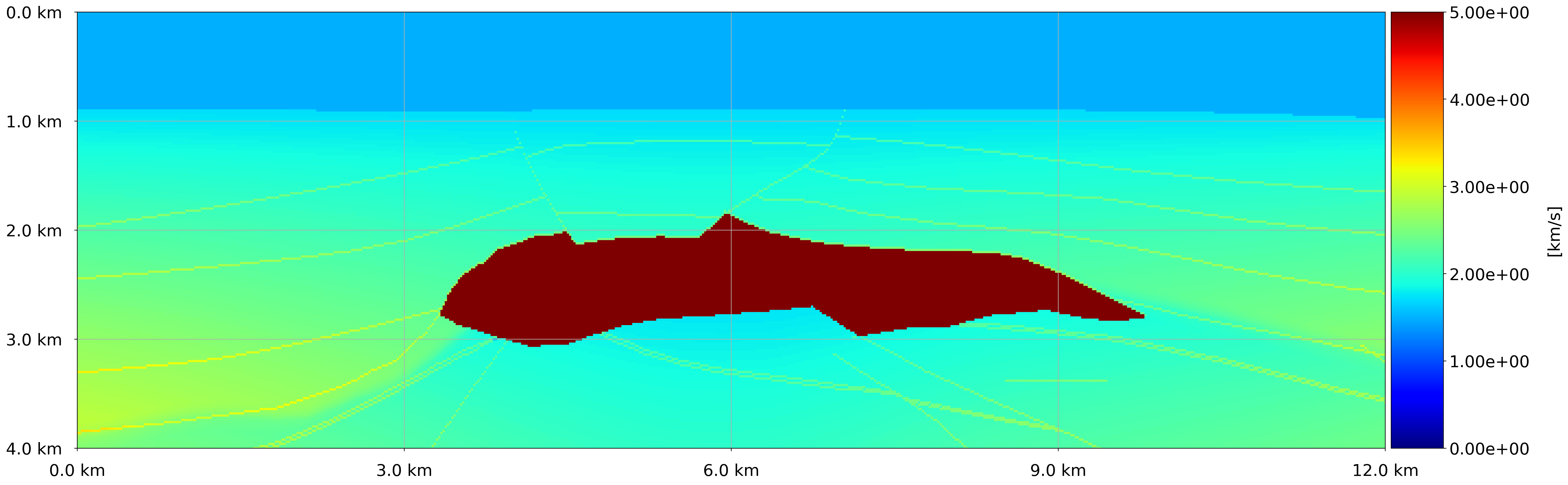} \caption{2D SEG/EAGE Salt Velocity Model.}
\label{fig:vel_seg}
\end{figure}

A well-known benchmark is performed using a 2D SEG/EAGE Salt Velocity Model \citep{aminzadeh19973,aminzadeh1996three,aminzadeh19953,aminzadeh1994seg,lecomte1994building}. This velocity model is widely used in several papers to evaluate the behavior of FD Schemes in a more realistic case, as it features a significant velocity contrast. In this model, the velocity varies across all fields, including regions that reproduce salt-dome profiles. We use a slice of the original velocity model, given in Figure \ref{fig:vel_seg}.
The set of parameters for the 2D SEG/EAGE Salt Velocity Model reproduces similar properties that we found in Seismic Imaging Problems, as we can see in Table \ref{table1:segeage2dsalt}, where we adopt the parameters of a similar test in \citep{wang2016effective}.
\begin{table}
\centering
\caption{Set of Parameters for 2D SEG/EAGE Salt Velocity Model}
\label{table1:segeage2dsalt}
\begin{tabular}{|c|c|}
\hline
\textbf{Parameters} & \textbf{Values}    \\
\hline
[$x_{I}$,$x_{F}$]    & [$0m$,$12000m$]   \\
\hline
[$z_{I}$,$z_{F}$]    & [$0m$,$4000m$]   \\
\hline
$t_{F}$    & $3000ms$   \\
\hline
$nx$ | $nz$       & $601$ | $201$   \\
\hline
$\Delta x = \Delta x$ & $20m$   \\
\hline
$\Delta t$ & $1.0ms$   \\
\hline
$nt$       & $3000$   \\
\hline
$f_{0}$    & $15Hz$   \\
\hline
\end{tabular}
\end{table}
Some information about the 2D SEG/EAGE Velocity Model is provided in Table \ref{table3:segeage2dsalt}, highlighting the source/receiver positions, as well as the number of receivers. The positions of the source/receivers replicate practical cases in which we seek to obtain information along the top of the oil basin.
\begin{table}
\centering
\caption{2D SEG/EAGE Salt Velocity Model Properties}
\label{table3:segeage2dsalt}
\begin{tabular}{|c|c|}
\hline
\textbf{Maximum Velocity} & $1.5Km/s$ \\
\hline
\textbf{Minimum Velocity} & $5.0Km/s$ \\
\hline
$\mathbf{x}$ \textbf{Source Position} & $6000m$ \\
\hline
$\mathbf{z}$ \textbf{Source Position} & $100m$ \\
\hline
$\mathbf{x}$ \textbf{Receiver Position} & from $0m$ to $12000m$ spaced by $20m$ \\
\hline
$\mathbf{z}$ \textbf{Receiver Position} & $20m$ \\
\hline
\end{tabular}
\end{table}

\noindent

In this case, a Damping strategy is again applied to minimize boundary reflections; however, the high velocity field heterogeneity introduces substantial wave interactions, complicating the analysis. Despite pronounced spatial variations in velocity, we construct optimized FD stencils using the maximum velocity, which is also used to determine the spatial–temporal parameters for each method and stencil configuration.

\begin{figure}[H]
\centering
\includegraphics[scale=.3]{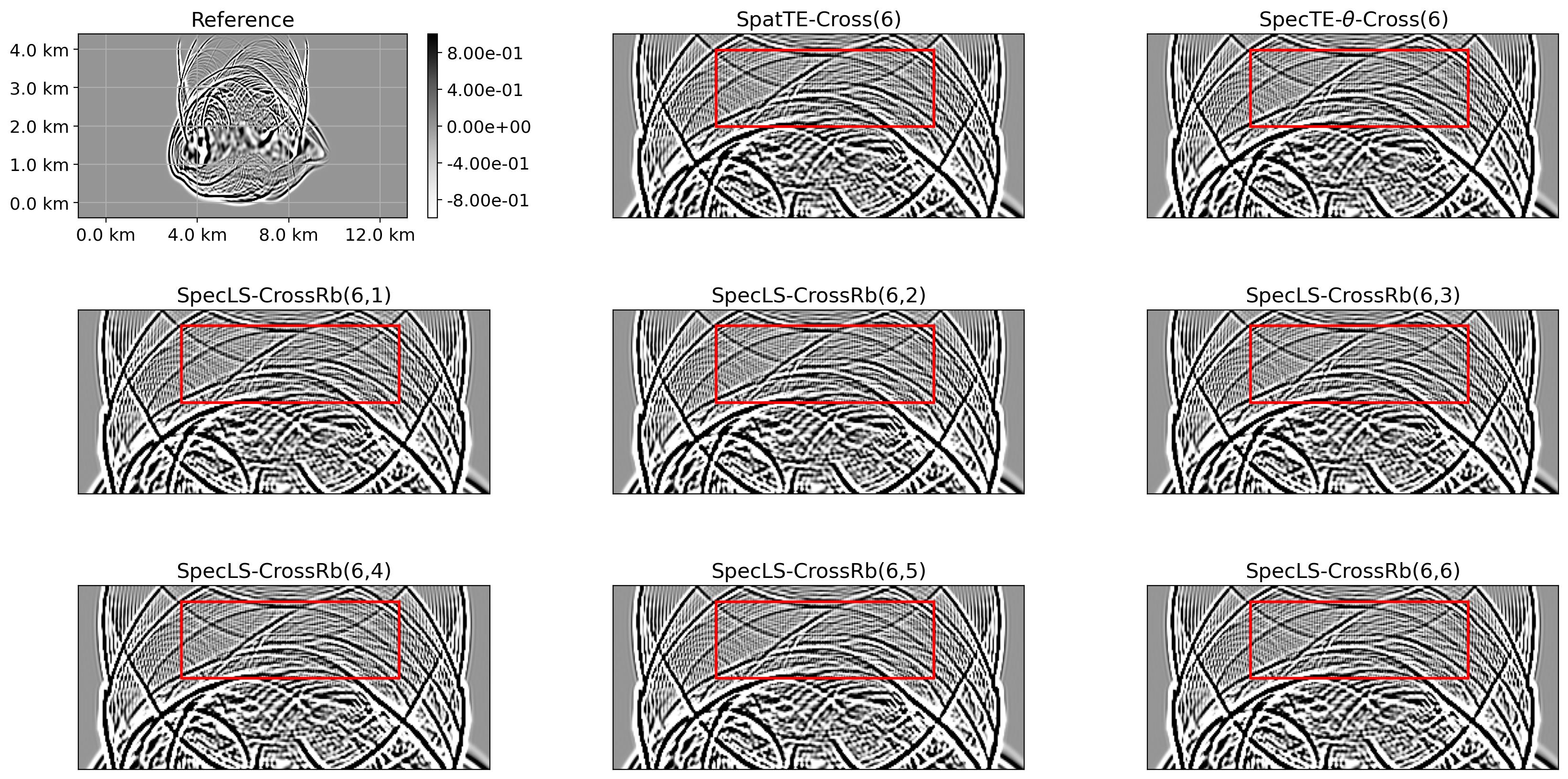} 
\caption{2D SEG EAGE Salt Velocity Model - Displacement - $SpecLS-CrossRb(6,N)$ at $T=1.95s$.}
\label{fig:segeage1}
\end{figure}

\begin{figure}[H]
\centering
\includegraphics[scale=.3]{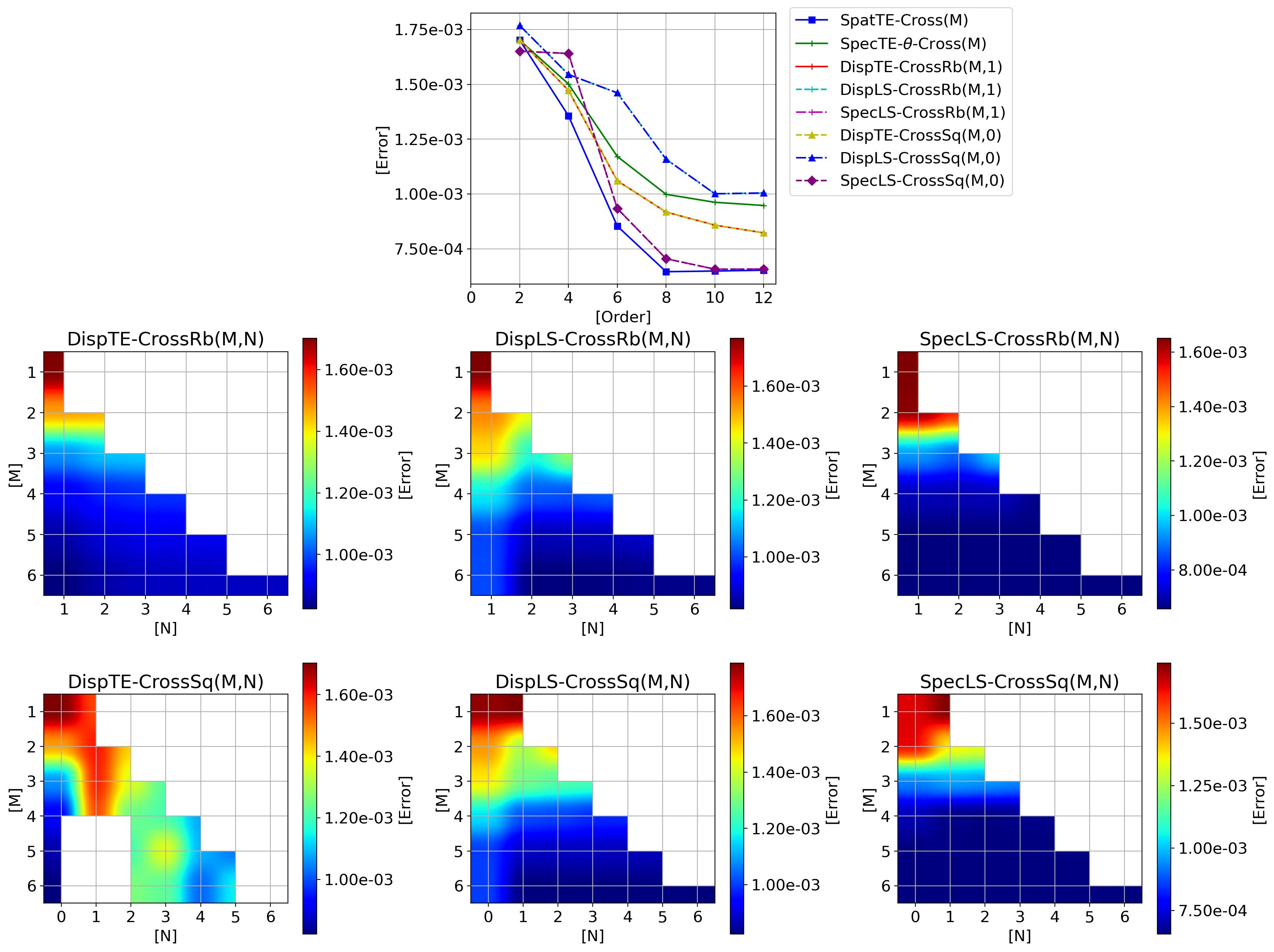} 
\caption{2D SEG EAGE Salt Velocity Model - Displacement - $\|\cdot\|_{2}$-norm at $T=1.95s$.}
\label{fig:segeage2}
\end{figure}

Following the same procedure as in the previous tests, both visual and numerical comparisons (using norm evaluations) are carried out, including analysis of the seismogram. Figure \ref{fig:segeage1} presents displacement results for $SpatTE$-$Cross(6)$, $SpecTE$-$\theta$-$Cross(6)$ and $SpecLS$-$CrossRb(6,N)$, while Figure~\ref{fig:segeage2} provides their corresponding norm quantification. Due to the strong wave interactions arising from the complex velocity field, the displacement plots in Figure \ref{fig:segeage1} are visually challenging to interpret. Red-bound boxes highlight regions of significant dispersion effects.

In Figure~\ref{fig:segeage1}, the dispersion effects are slightly reduced for $SpatTE$-$Cross(6)$, and the $SpecLS$-$CrossRb(6,N)$ solution shows noticeably weaker dispersion artifacts. The norm comparison in Figure~\ref{fig:segeage2} confirms that, for the $Cross$ configuration, $SpecLS$-$CrossRb(M,N)$ and $SpatTE$-$Cross(M)$ achieve the best performance among all tested schemes. As observed in previous tests, $DispTE$-$CrossSq(M,N)$ exhibits instability in certain regions, while the norm generally decreases with increasing $M$ and shows little variation with $N$.

A central insight from these results concerns why $SpecLS$-$CrossRb(M,N)$ and $SpatTE$-$Cross(M)$ outperform other optimized schemes. These two families of methods share a key characteristic: the velocity model is not directly used in the computation of their FD weights. Specifically, in $SpecLS$-$CrossRb(M,N)$ and $SpecLS$-$CrossSq(M,N)$, only the spatial parameter is considered in the optimization. This distinction becomes critical in highly heterogeneous environments. As shown in \citep{wang2016effective}, low-order optimized FD schemes tend to perform poorly in the \textit{SEG/EAGE 2D Salt} model. The study proposes two remedies: the use of very high-order stencils (e.g., $M=28$) or adaptive schemes in which $M$ decreases as the velocity increases, while $N$ remains fixed and small. Both strategies improve accuracy and reduce dispersion but incur high computational costs, particularly in full-waveform inversion (FWI) applications.

\begin{figure}[H]
\centering
\includegraphics[scale=.3]{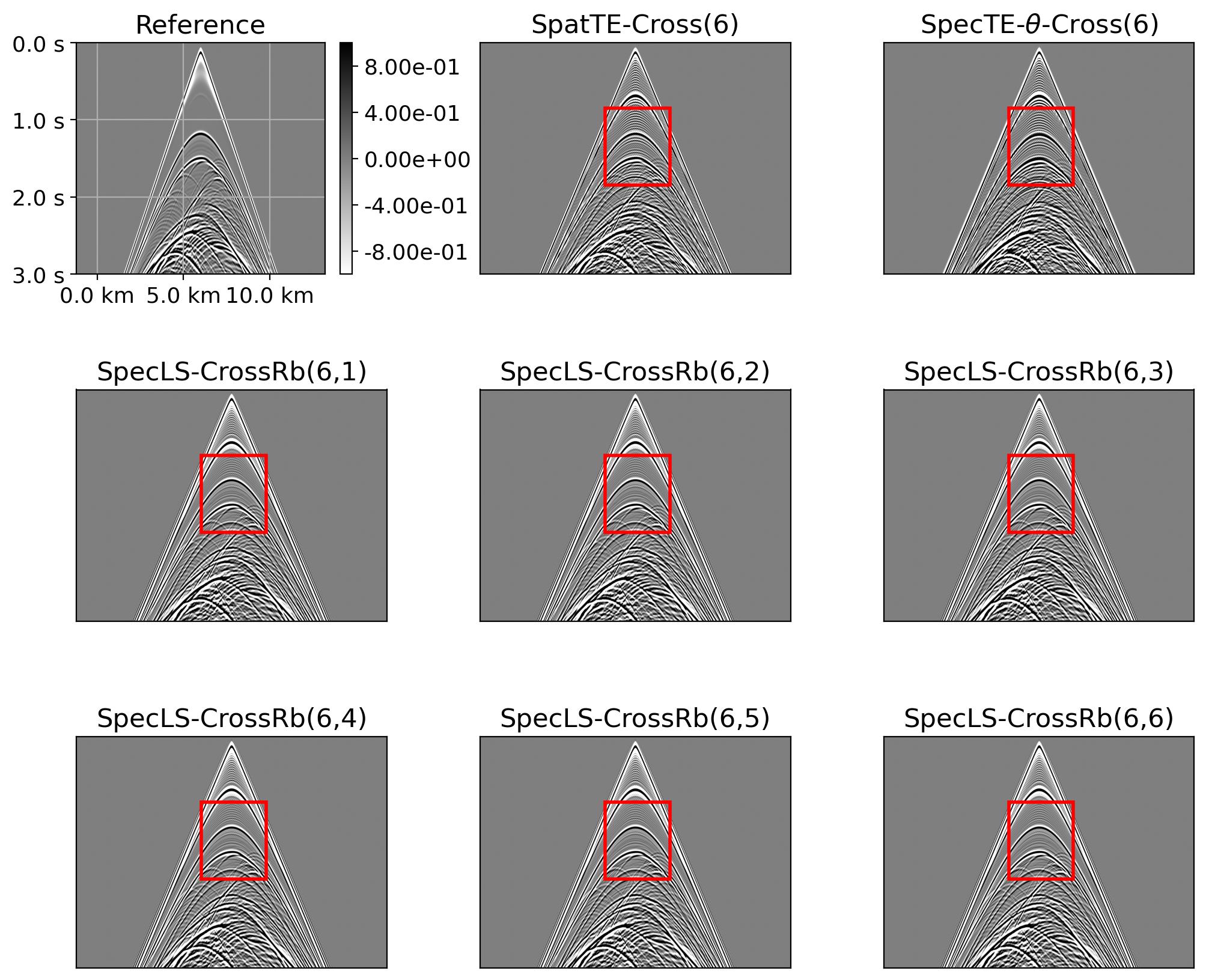} 
\caption{2D SEG EAGE Salt Velocity Model - Receiver - $SpecLS-CrossRb(6,N)$ at $T=1.95s$.}
\label{fig:segeage3}
\end{figure}

\begin{figure}[H]
\centering
\includegraphics[scale=.3]{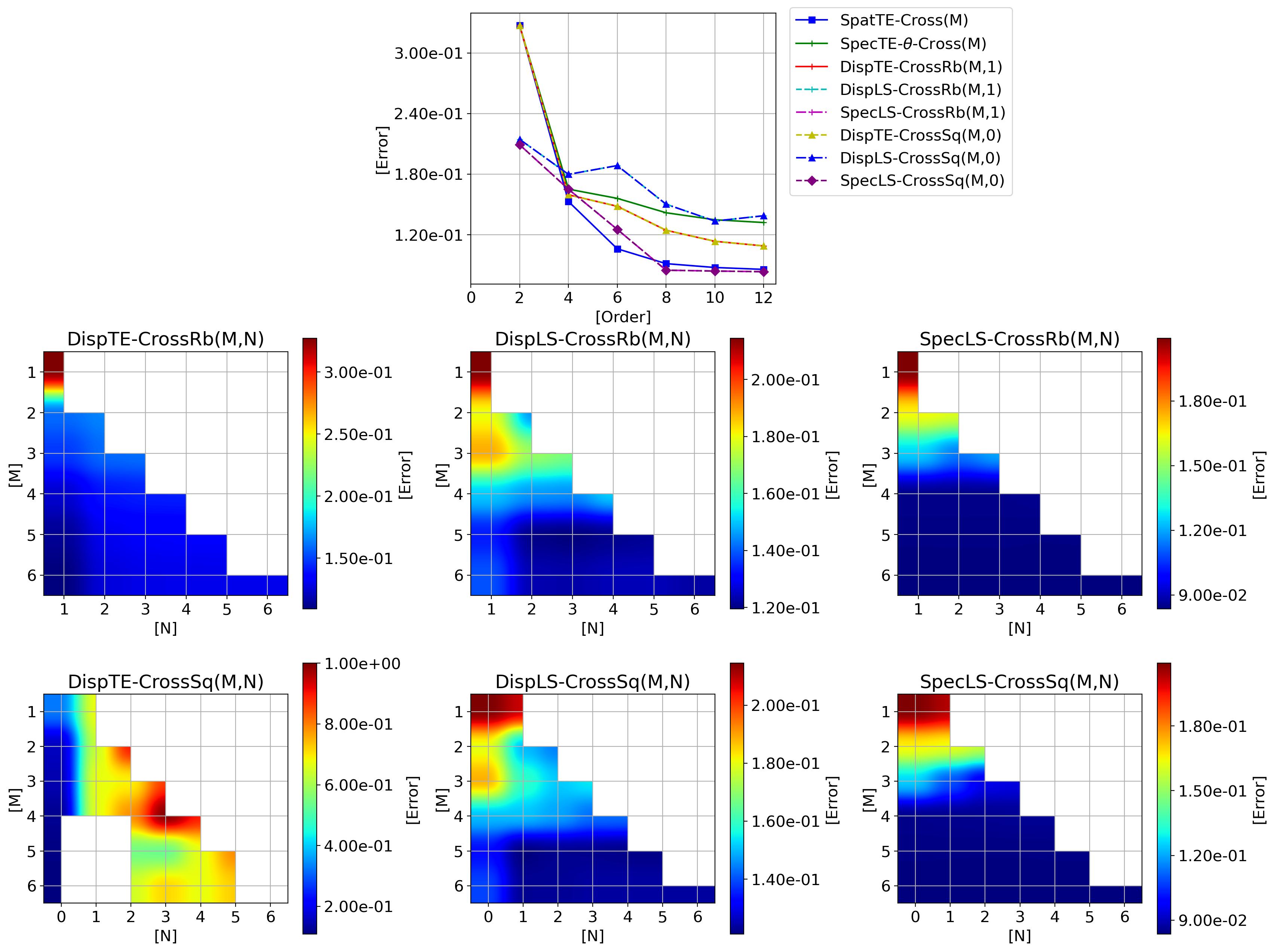} 
\caption{2D SEG EAGE Salt Velocity Model - Receiver - $\|\cdot\|_{2}$-norm at $T=1.95s$.}
\label{fig:segeage4}
\end{figure}

\begin{figure}[H]
\centering
\includegraphics[scale=.3]{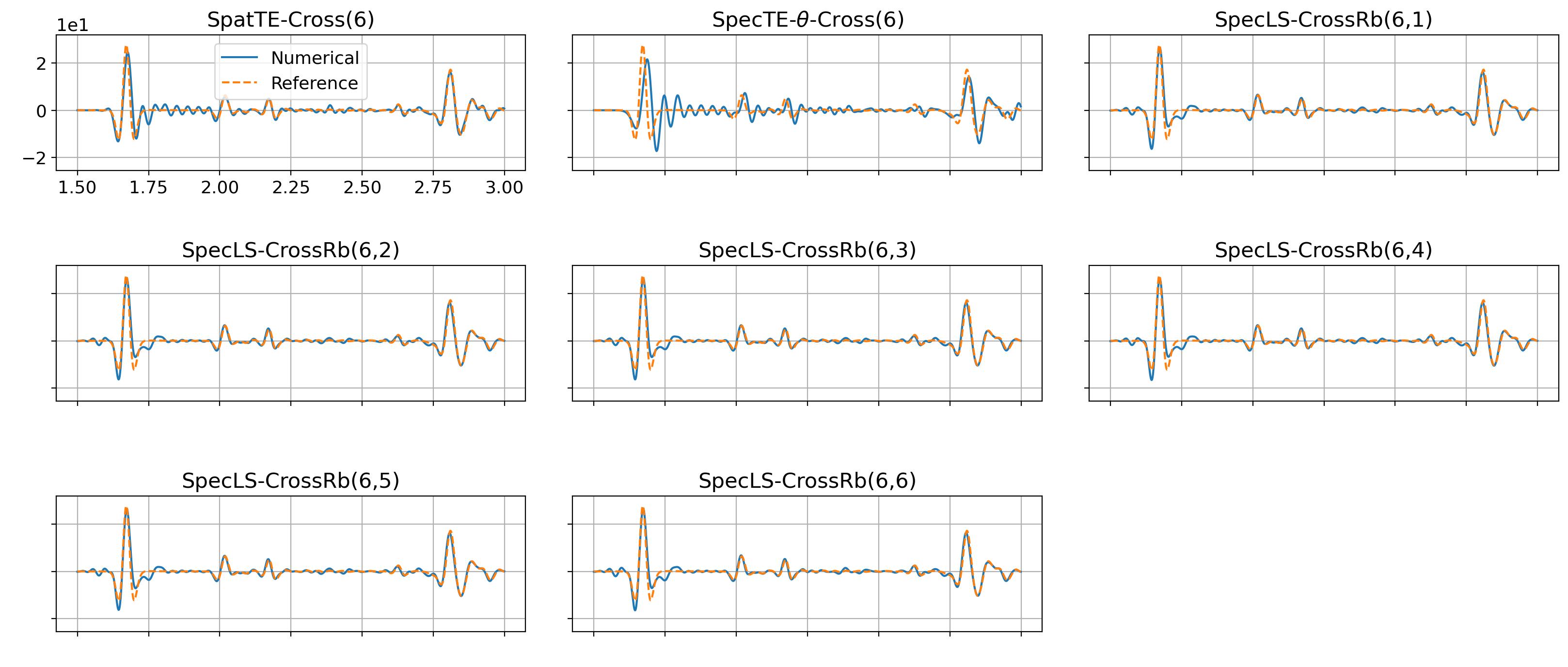} 
\caption{2D SEG EAGE Salt Velocity Model - Seismic Trace - $SpecLS-CrossRb(6,N)$ at $x=3600m$ and $z=20m$.}
\label{fig:segeage5}
\end{figure}

\begin{figure}[H]
\centering
\includegraphics[scale=.3]{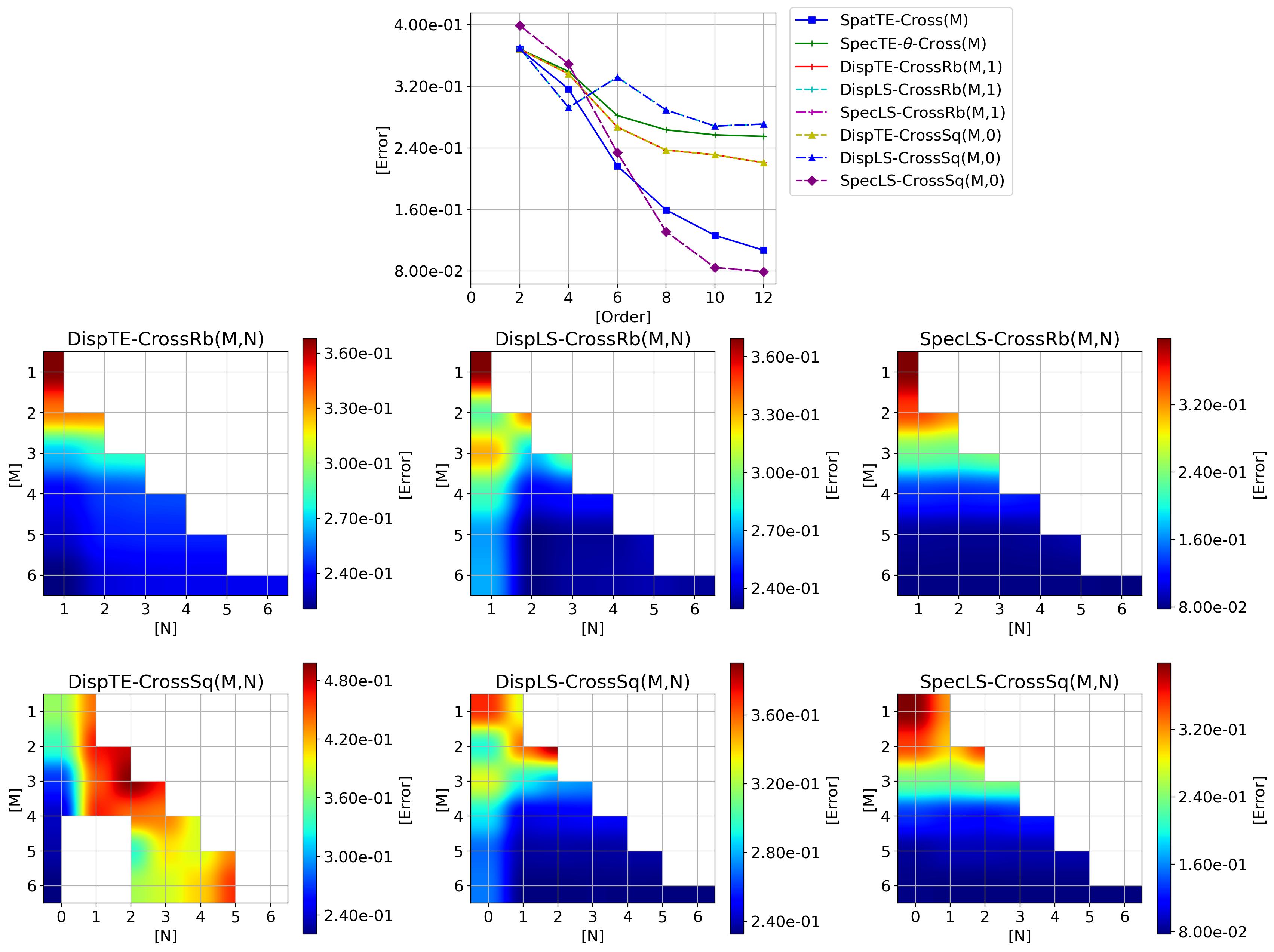} 
\caption{2D SEG EAGE Salt Velocity Model - Seismic Trace - $\|\cdot\|_{2}$-norm at $x=3600m$ and $z=20m$.}
\label{fig:segeage6}
\end{figure}

In Figures~\ref{fig:segeage3} and ~\ref{fig:segeage4}, dispersion effects over the seismograms are highlighted with red-bound boxes, demonstrating the reduced artifacts in $SpecLS$-$CrossRb(M,N)$ compared with $DispLS$-$CrossRb(M,N)$. Seismic trace comparisons in  Figures~\ref{fig:segeage5} and ~\ref{fig:segeage6} emphasize these findings: $SpecLS$-$CrossRb(M,1)$ aligns closely with the reference solution and clearly outperforms $SpatTE$-$Cross(6)$, providing a magnified view of these results, reinforcing the improved phase and amplitude accuracy of $SpecLS$-$CrossRb(M,N)$.

\begin{figure}[H]
\centering
\includegraphics[scale=.3]{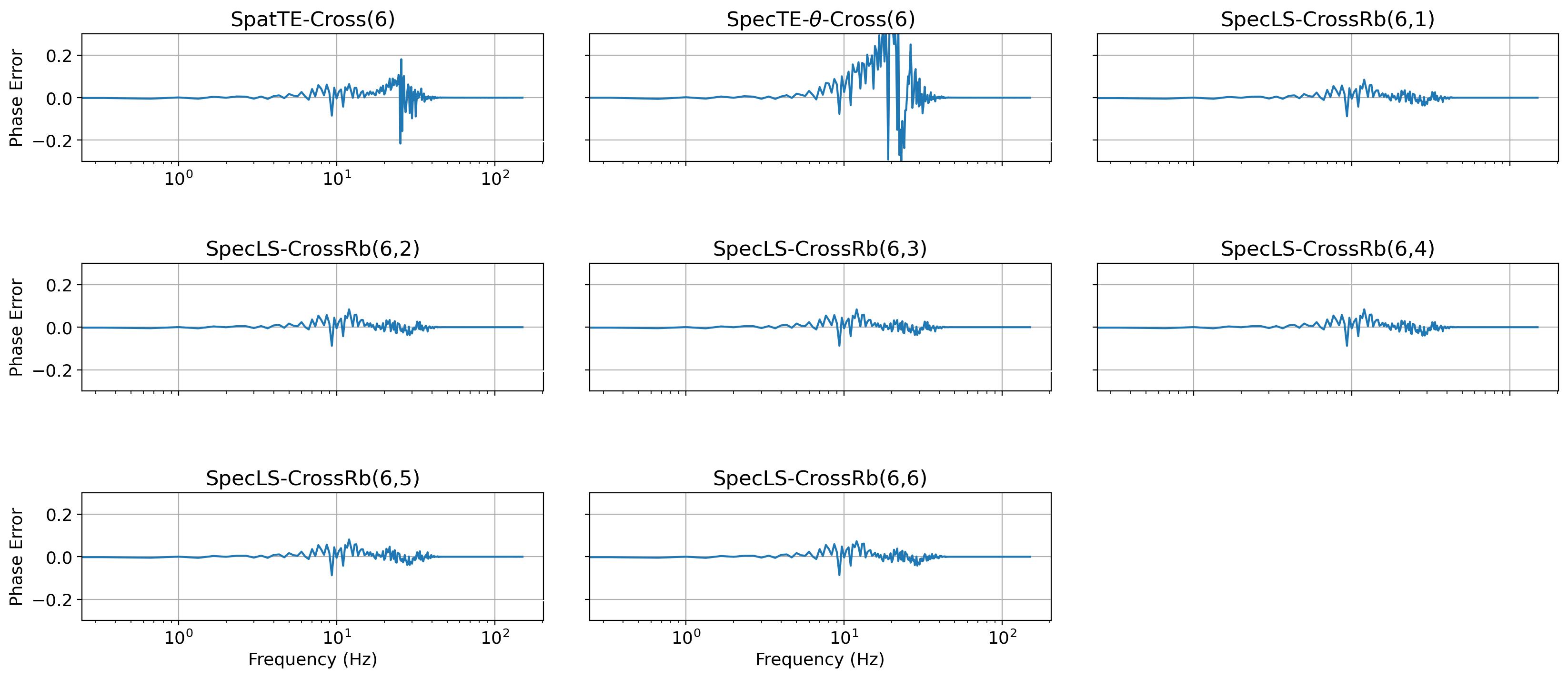} 
\caption{2D SEG EAGE Salt Velocity Model - FFT Analysis of Seismic Trace - $SpecLS-CrossRb(6,N)$ at $x=3600m$ and $z=20m$.}
\label{fig:segeage7}
\end{figure}

\begin{figure}[H]
\centering
\includegraphics[scale=.3]{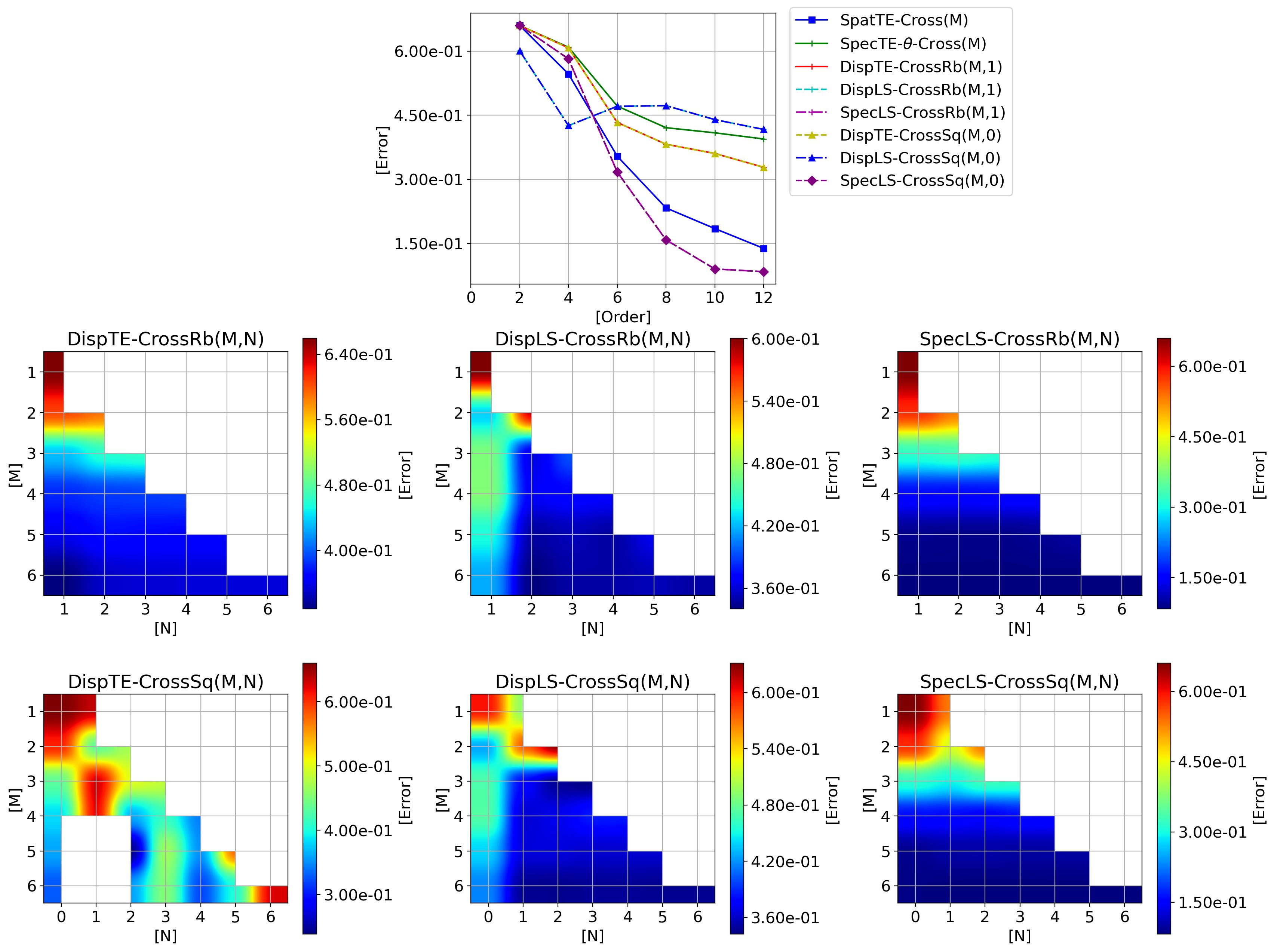}
\caption{2D SEG EAGE Salt Velocity Model - Seismic Trace - $\|\gamma^{\mathscr{F}}_{u_{st}(x_{*},z_{*})}\|_{2} 
$-norm at $x=3600m$ and $z=20m$.}
\label{fig:segeage8}
\end{figure}

Finally, the FFT dispersion curves in Figures~\ref{fig:segeage7} along with the corresponding norm quantification in Figure~\ref{fig:segeage8}, demonstrate how much more dispersive $DispLS$-$CrossRb(6,N)$ is compared to $SpecLS$-$CrossRb(6,N)$. These findings corroborate the previous analyses, where $DispLS$-$CrossRb(M,N)$ consistently produced inferior results relative to $SpecLS$-$CrossRb(M,N)$.

In summary, the \textit{SEG/EAGE 2D Salt} model illustrates the sensitivity of optimized FD schemes to velocity variations. When the velocity field is highly heterogeneous, schemes that rely on velocity-dependent weight computation tend to degrade in performance. Conversely, FD schemes that depend only on spatial parameters, such as $SpecLS$-$CrossRb(M,N)$ and $SpatTE$-$Cross(M)$, demonstrate greater robustness and reduced dispersion. These observations underscore the importance of carefully managing the dependence of FD optimization on both spatial–temporal parameters and velocity values to ensure stability and accuracy in complex media.

\subsection{Marmousi Velocity Model}
Another well-known benchmark is performed using a Marmousi Velocity Model \citep{versteeg1994marmousi}. As 2D SEG/EAGE Salt Model, this velocity model is widely used in several papers to evaluate the behavior of FD Schemes in a more realistic case, as it features a significant velocity contrast. We use a slice of the original velocity model, given in Figure \ref{fig:vel_seg}.
\begin{figure}[H]
\centering
\includegraphics[scale=.15]{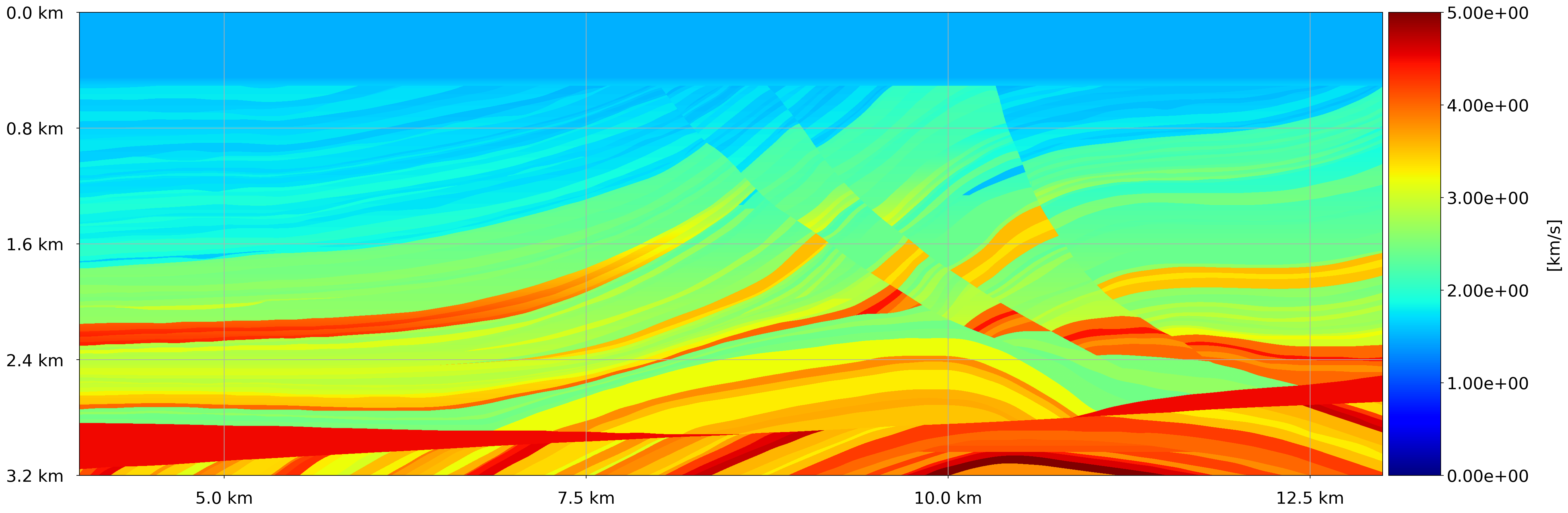} \caption{Marmousi Velocity Model.}
\label{fig:vel_marmousi}
\end{figure}
The set of parameters for the Marmousi Velocity Model reproduces similar properties that we found in Seismic Imaging Problems, as we can see in Table \ref{table1:marmousi}.
\begin{table}
\centering
\caption{Set of Parameters for Marmousi Velocity Model}
\label{table1:marmousi}
\begin{tabular}{|c|c|}
\hline
\textbf{Parameters} & \textbf{Values}    \\
\hline
[$x_{I}$,$x_{F}$]    & [$4000m$,$13000m$]   \\
\hline
[$z_{I}$,$z_{F}$]    & [$0m$,$3200m$]   \\
\hline
$t_{F}$    & $3000ms$   \\
\hline
$nx$       & $226$   \\
\hline
$nz$       & $81$   \\
\hline
$\Delta x = \Delta z$ & $40m$   \\
\hline
$\Delta t$ & $1.0ms$   \\
\hline
$nt$       & $3000$   \\
\hline
$f_{0}$    & $15Hz$   \\
\hline
\end{tabular}
\end{table}
Some information about the Marmousi Velocity Model is provided in Table \ref{table3:marmousi}, highlighting the source/receiver positions, as well as the number of receivers. The positions of the source/receivers replicate practical cases in which we seek to obtain information along the top of the oil basin.
\begin{table}
\centering
\caption{Marmousi Velocity Model Properties}
\label{table3:marmousi}
\begin{tabular}{|c|c|}
\hline
\textbf{Maximum Velocity} & $1.5Km/s$ \\
\hline
\textbf{Minimum Velocity} & $5.0Km/s$ \\
\hline
$\mathbf{x}$ \textbf{Source Position} & $8500m$ \\
\hline
$\mathbf{z}$ \textbf{Source Position} & $50m$ \\
\hline
$\mathbf{x}$ \textbf{Receiver Position} & from $4000m$ to $13000m$ spaced by $40m$ \\
\hline
$\mathbf{z}$ \textbf{Receiver Position} & $50m$ \\
\hline
\end{tabular}
\end{table}

\noindent

\begin{figure}[H]
\centering
\includegraphics[scale=.3]{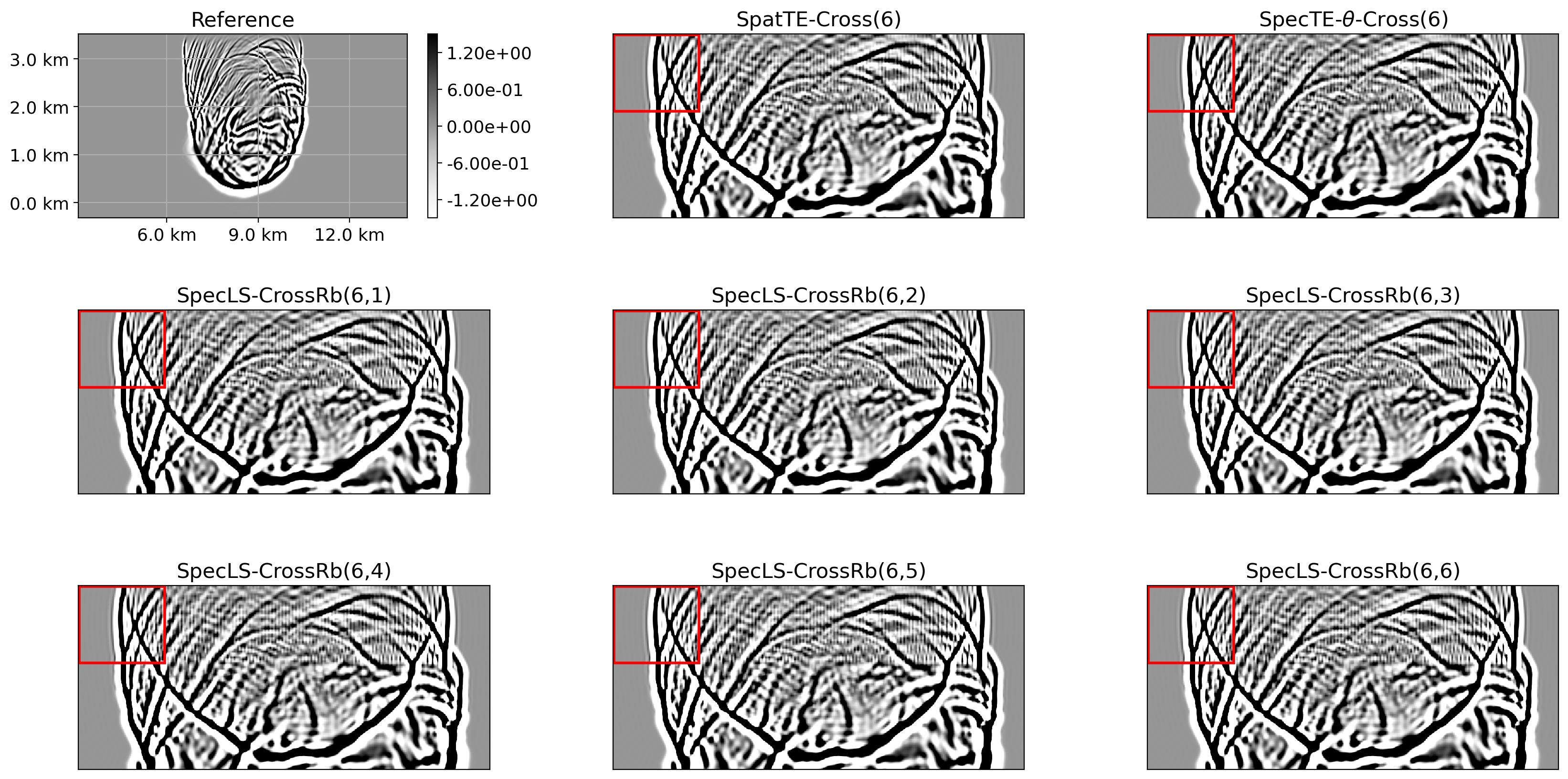} 
\caption{Marmousi Velocity Model - Displacement - $SpecLS-CrossRb(6,N)$ at $T=1.95s$.}
\label{fig:marmousi1}
\end{figure}

\begin{figure}[H]
\centering
\includegraphics[scale=.3]{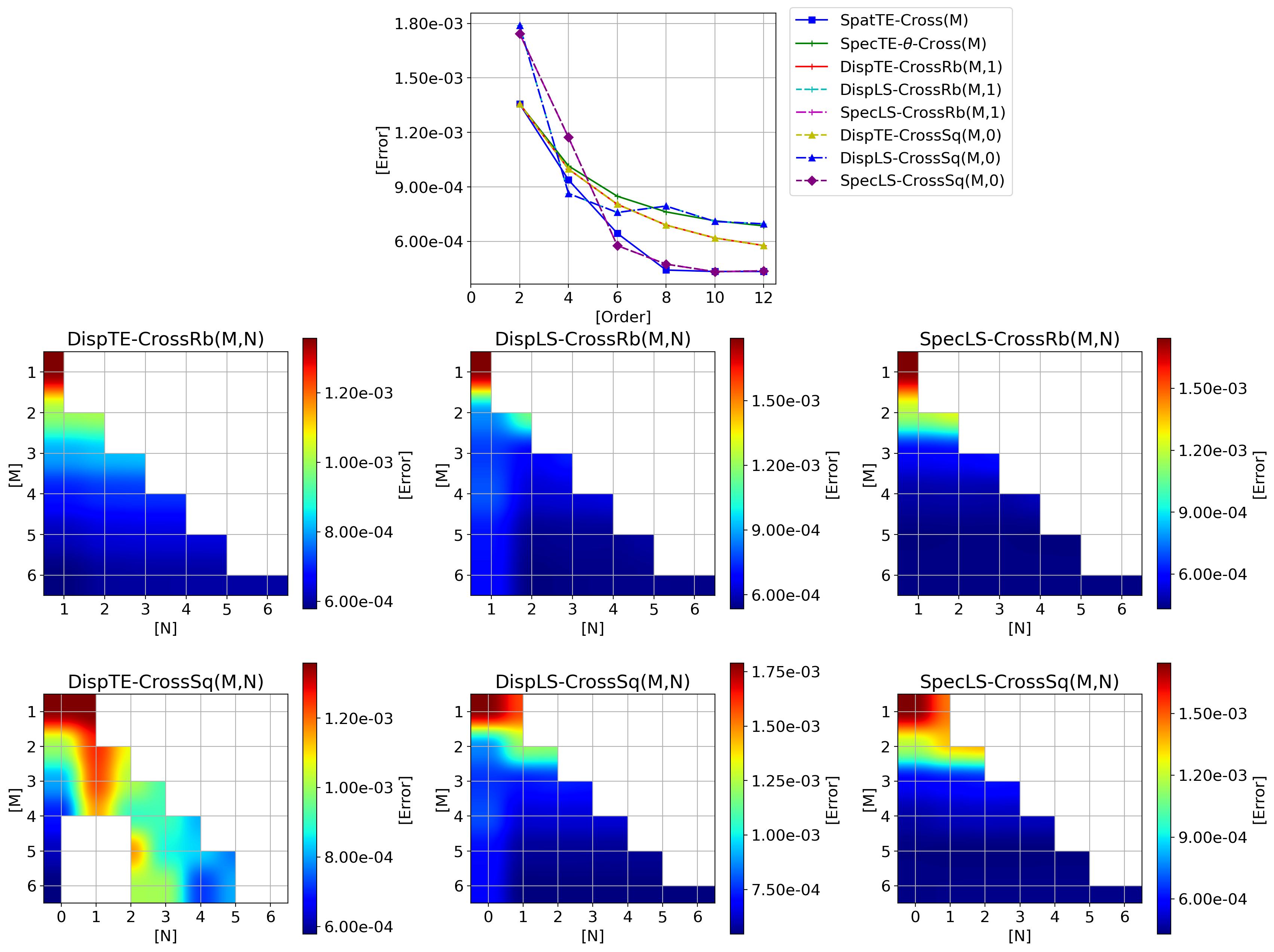} 
\caption{Marmousi Velocity Model - Displacement - $\|\cdot\|_{2}$-norm at $T=1.95s$.}
\label{fig:marmousi2}
\end{figure}

The \textit{Marmousi} model reinforces the strong sensitivity of optimized FD schemes to complex velocity variations, similarly to the \textit{2D SEG EAGE Salt} model (Figures~\ref{fig:marmousi1} and~\ref{fig:marmousi2}). In highly heterogeneous media, wave interactions become more complex, increasing dispersion effects and reducing the accuracy of velocity-dependent FD optimizations when a fixed scheme is applied across the grid. In contrast, schemes based mainly on spatial parameters, such as $SpecLS$-$CrossRb(M,N)$ and $SpatTE$-$Cross(M)$, generally present lower dispersion artifacts, as indicated by the norm evaluations in Figure~\ref{fig:marmousi2}.

\begin{figure}[H]
\centering
\includegraphics[scale=.3]{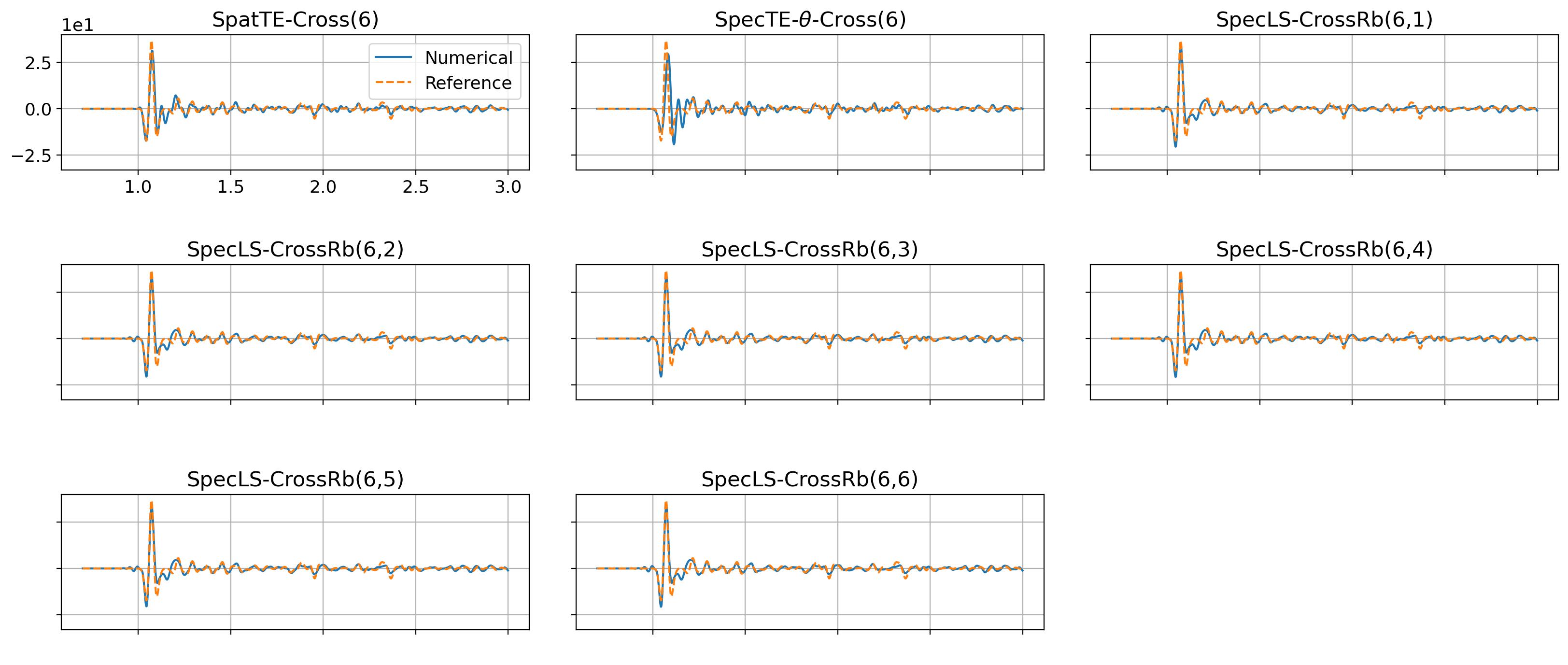} 
\caption{Marmousi Velocity Model - Seismic Trace - $SpecLS-CrossRb(6,N)$ at $x=50m$ and $z=7000m$.}
\label{fig:marmousi5}
\end{figure}

\begin{figure}[H]
\centering
\includegraphics[scale=.3]{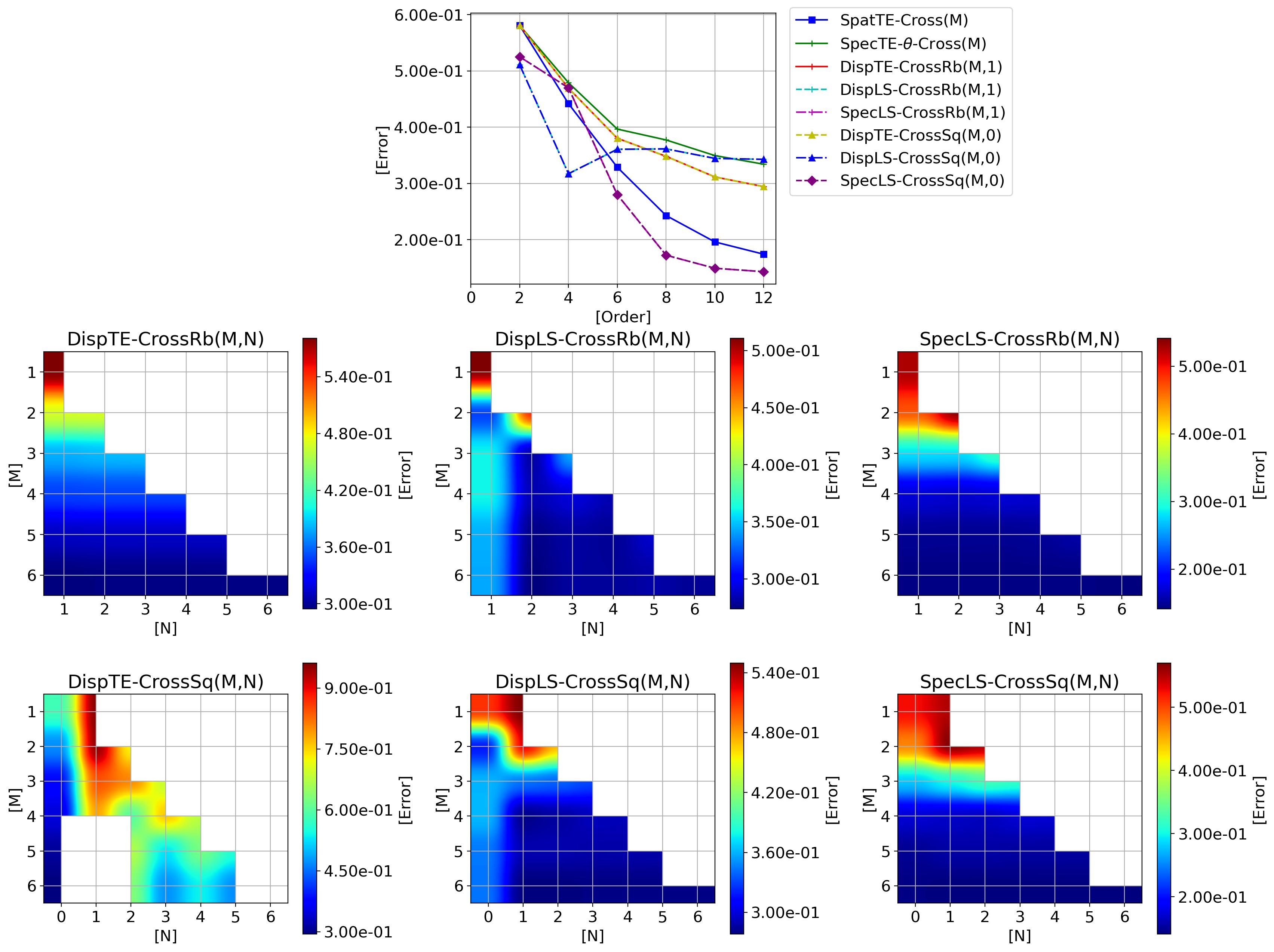} 
\caption{Marmousi Velocity Model - Seismic Trace - $\|\cdot\|_{2}$-norm at $x=50m$ and $z=7000m$.}
\label{fig:marmousi6}
\end{figure}

\begin{figure}[H]
\centering
\includegraphics[scale=.3]{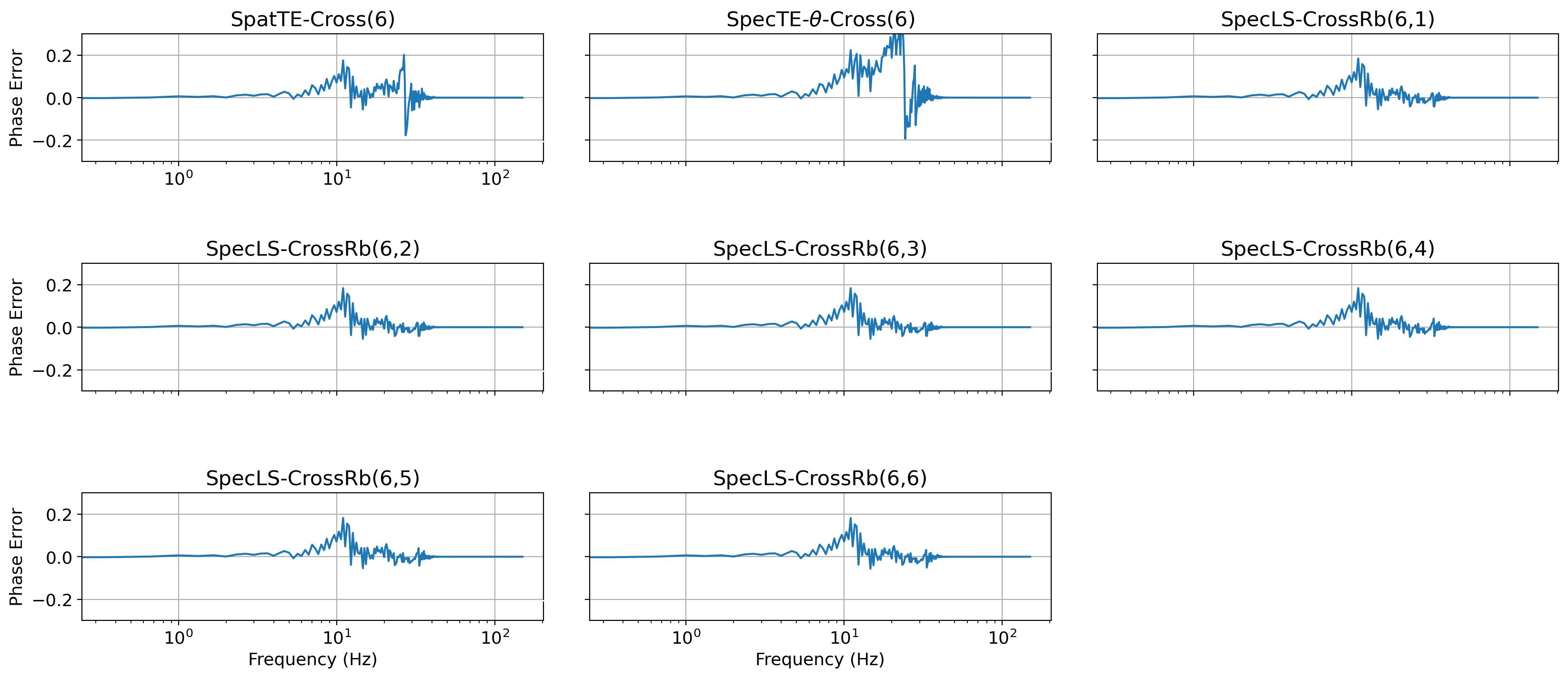} 
\caption{Marmousi Velocity Model - FFT Analysis of Seismic Trace - $SpecLS-CrossRb(6,N)$ at $x=50m$ and $z=7000m$.}
\label{fig:marmousi7}
\end{figure}

\begin{figure}[H]
\centering
\includegraphics[scale=.3]{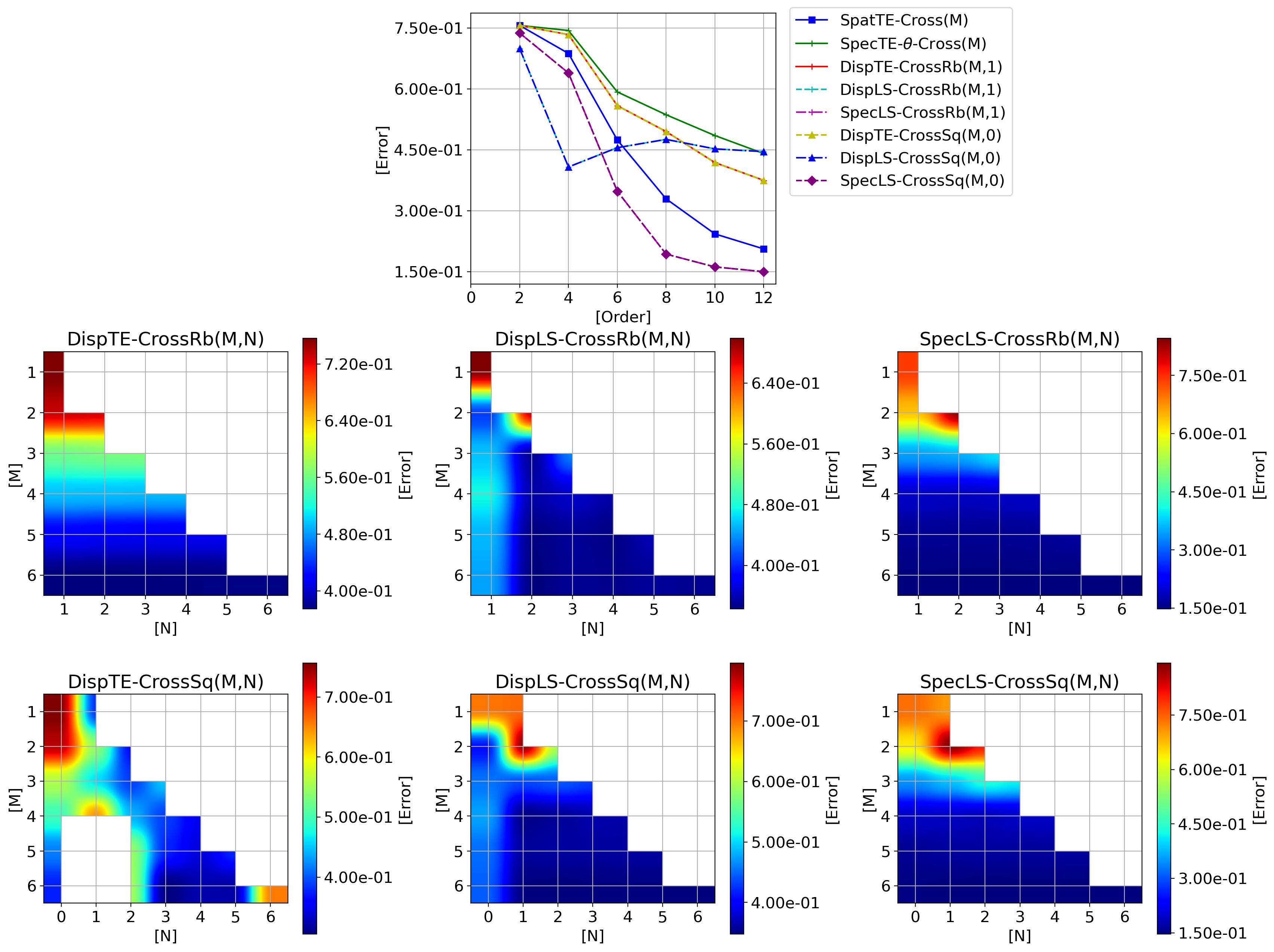}
\caption{Marmousi Velocity Model - Seismic Trace - $\|\gamma^{\mathscr{F}}_{u_{st}(x_{*},z_{*})}\|_{2} 
$-norm at $x=50m$ and $z=7000m$.}
\label{fig:marmousi8}
\end{figure}

Although $SpecLS$-$CrossRb(M,N)$ and $SpatTE$-$Cross(M)$ showed slightingly better results in some norm evaluations, for selected seismic traces (Figures~\ref{fig:marmousi5} and~\ref{fig:marmousi6}) and FFT analyses (Figures~\ref{fig:marmousi7} and~\ref{fig:marmousi8}) demonstrated their ability to reduce dispersion and remain closer to the reference solution compared to $DispLS$-$CrossRb(6,N)$.

Overall, the results emphasize that excessive reliance on heterogeneous velocity fields in FD optimization can compromise accuracy, whereas spatially optimized schemes tend to provide greater robustness in complex media. As discussed by \citet{wang2016effective}, combining optimized FD schemes across different velocity regions may be a more effective strategy to reduce dispersion effects in high-contrast models.

\pagebreak

\pagebreak

\subsection{Spatial Resolution and Space Order Trade-Off}
Based on the previous experiments, we observe that the Cross geometry yields better results than the Rhombus and Square geometries. Here, we examine the behavior of the Cross geometry as we simultaneously vary the spatial order and mesh refinement to achieve a trade-off between them.
For this purpose, we retain the same parameter settings used in the previous tests and vary only the spatial discretization parameter and the spatial order, evaluating in both the standard and the proposed FFT-based norms.

\begin{figure}[H]
\centering
\includegraphics[scale=.20]{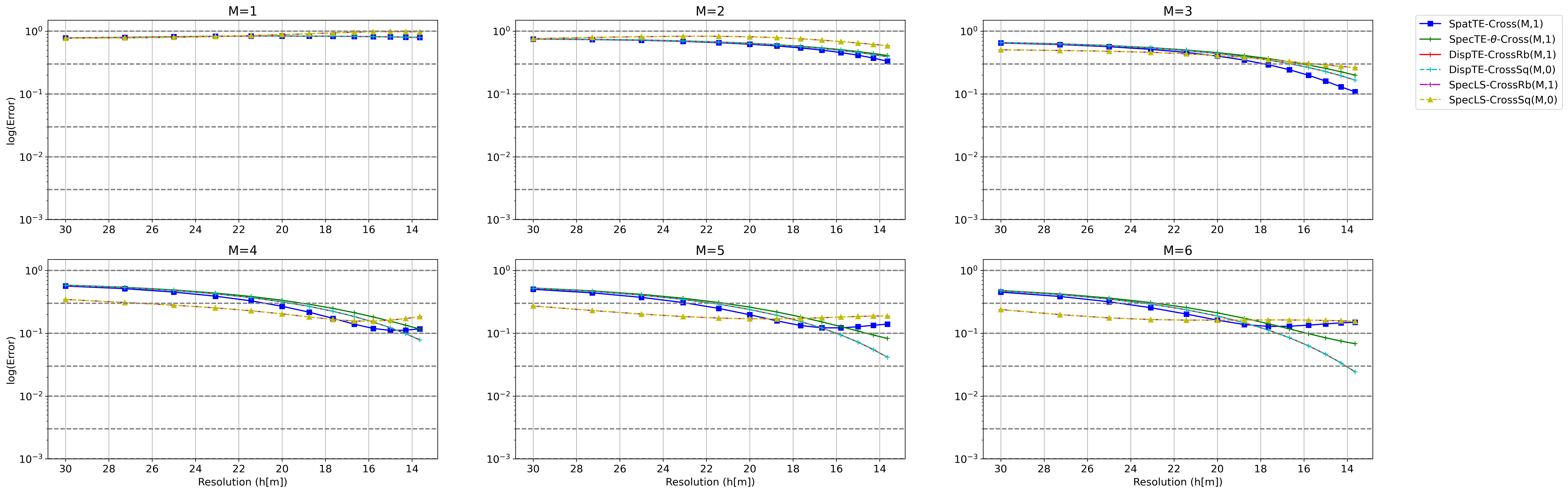}
\caption{Homogeneous Velocity Model - Seismic Trace - $\|\cdot\|_{2}$-norm at $x=1200m$ and $z=3000m$ for different spatial resolutions ($h$) and space order ($M$).}
\label{fig:spatial_tradeoff1}
\end{figure}

\begin{figure}[H]
\centering
\includegraphics[scale=.20]{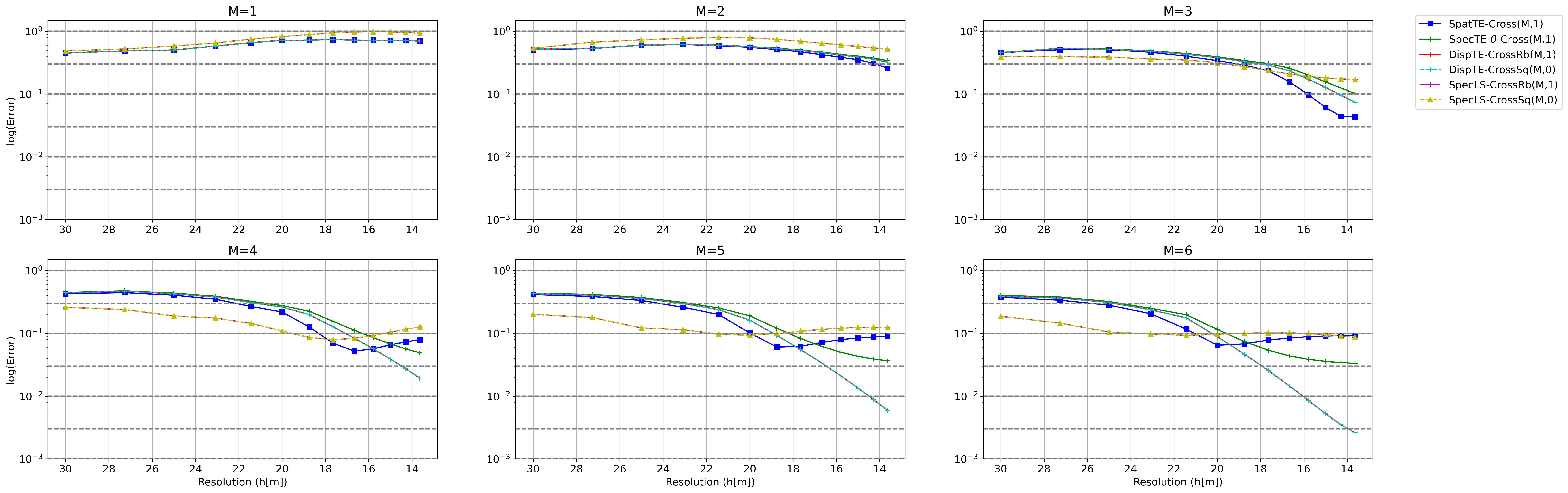}
\caption{Homogeneous Velocity Model - Seismic Trace - $\|\gamma^{\mathscr{F}}_{u_{st}(x_{*},z_{*})}\|_{2} $-norm at $x=1200m$ and $z=3000m$ for different spatial resolutions ($h$) and space order ($M$).}
\label{fig:spatial_tradeoff2}
\end{figure}

Figures \ref{fig:spatial_tradeoff1} and \ref{fig:spatial_tradeoff2} present the norm evaluations for the Homogeneous Velocity Model. 
In Figure \ref{fig:spatial_tradeoff1}, we evaluate the $\|\cdot\|{2}$ norm, where for order $M \leq 4$, all schemes exhibit similar behavior, with very small differences in norm. This indicates that, for lower-order FD methods, the Cross schemes exhibited similar behavior.
However, for order $M \geq 5$, the $DispTE$-$CrossRb(M,1)$ scheme shows a significant improvement, where the norm difference values increase significantly, showing that optimized FD schemes can hit best approximations in homogeneous velocity regions.  
In terms of the $\|\gamma^{\mathscr{F}}{u_{st}(x_{},z_{})}\|_{2}$-norm, we observe in Figure \ref{fig:spatial_tradeoff2} that for $M \geq 3$ $SpatTE$-$Cross(M)$ has less dispersion effects than other methods, when we refine the mesh. However, for $M \geq 4$, we observe that the $DispTE$-$CrossRb(M,1)$ scheme has an efficient dispersion reduction effect, considering a homogeneous velocity field.
Overall, both error measures exhibit the expected behavior: as the spatial resolution is refined and the spatial order increases, the standard $\|\cdot\|{2}$-norm decreases, and the same trend is observed for the dispersion error measured by the $\|\gamma^{\mathscr{F}}{u{st}(x_{},z_{})}\|_{2}$-norm. In some cases, the error curve became flat, indicating equilibrium in spatial order and refinement, accounting for dispersion error.
In general, for a homogeneous velocity model, low-order FD schemes ($M\leq3$) have a similar behavior in terms of standard error and dispersion error, showing that $SpatTE$-$Cross(M)$ can be chosen due to simplicity. 
On the other hand, looking for high-order FD schemes ($M \geq 4$), significant results can be reached using optimized schemes and mesh refinement, highlighting $DispTE$-$CrossRb(M,1)$. 

\begin{figure}[H]
\centering
\includegraphics[scale=.20]{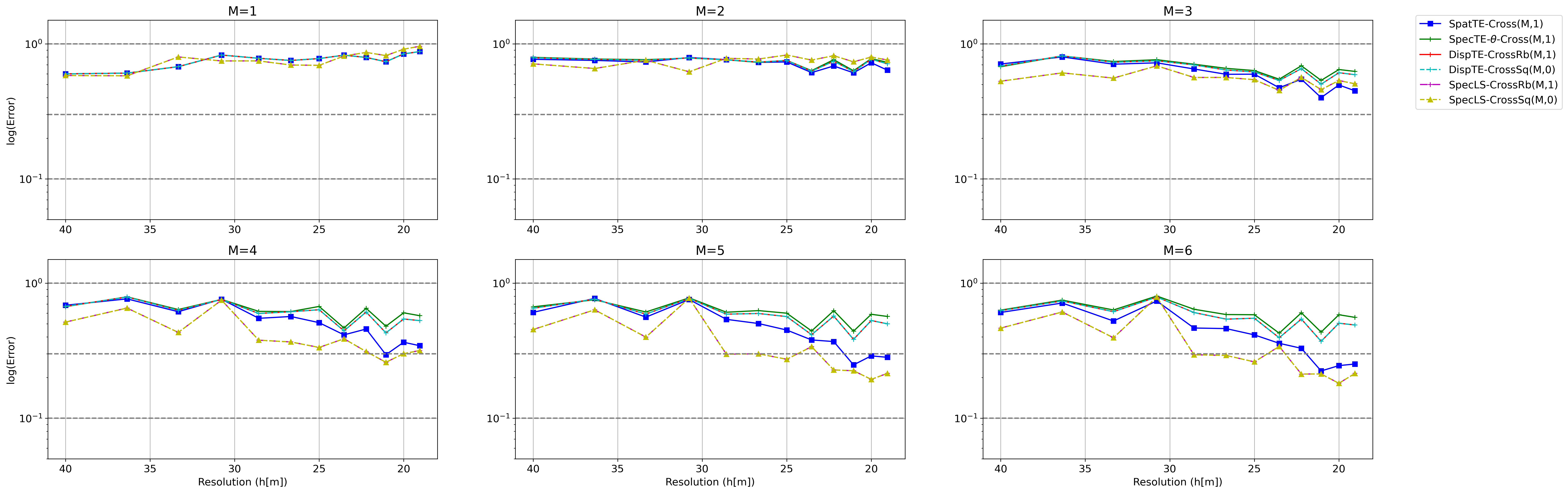}
\caption{2D SEG EAGE Salt Velocity Model  - Seismic Trace - $\|\cdot\|_{2}$-norm at $x=3600m$ and $z=20m$ for different spatial resolutions and space order.}
\label{fig:spatial_tradeoff3}
\end{figure}

\begin{figure}[H]
\centering
\includegraphics[scale=.20]{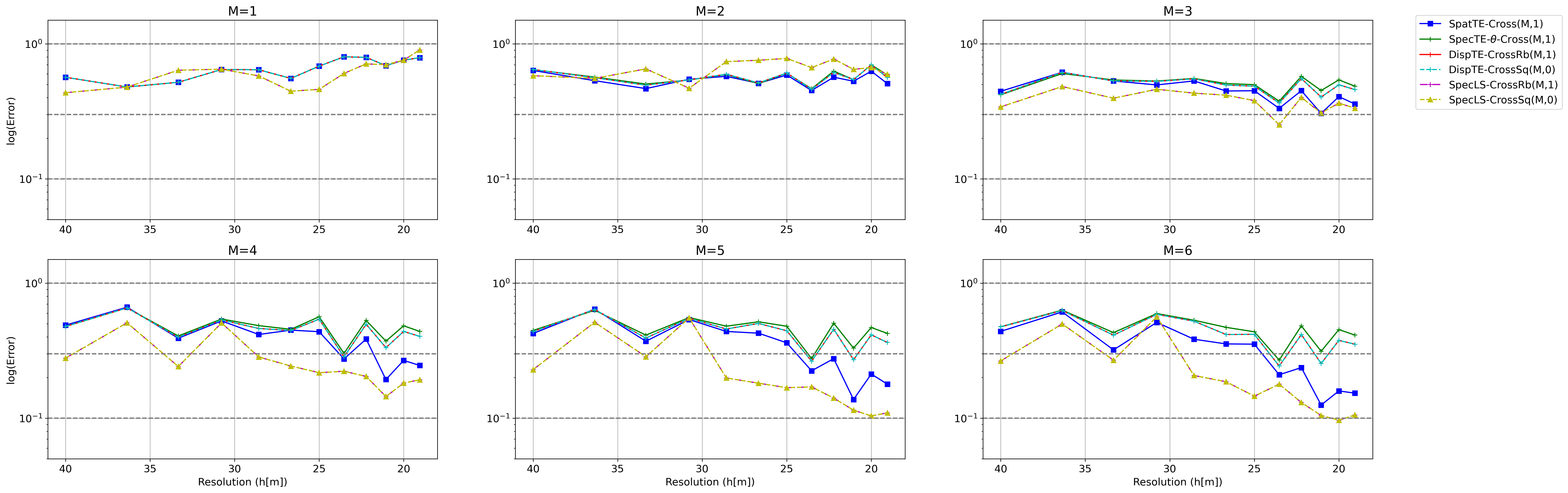}
\caption{2D SEG EAGE Salt Velocity Model - $\|\gamma^{\mathscr{F}}_{u_{st}(x_{*},z_{*})}\|_{2} $-norm at $x=3600m$ and $z=20m$ for different spatial resolutions and space order.}
\label{fig:spatial_tradeoff4}
\end{figure}

Following the test, Figures \ref{fig:spatial_tradeoff3} and \ref{fig:spatial_tradeoff4} show the results for the 2D SEG/EAGE Salt Velocity Model, where we have a complex velocity scenario. 
Since the proposed optimized schemes are constructed using the maximum velocity, this characteristic directly affects the error behavior, changing the usual convergence aspects, that is, refinement summing to high order leads to lower error values.  
In this test, the $SpatTE$-$Cross(M)$ and $SpecLS$-$CrossRb(M,1)$/$SpecLS$-$CrossSq(M,0)$ schemes provide the best results for both norms. Furthermore, as spatial resolution and spatial order increase, both norms tend to decrease, but not monotonically as expected; in some cases, we observe both growth and decrease in the norms.
In this heterogeneous test, increasing the FD order and/or mesh refinement can help reduce the norm values, but this is not the most important aspect. Here, changes in velocity affect the ability to choose the best option in a trade-off between order and refinement.

Overall, these experiments reinforce the conclusions drawn from the previous analyses. For the Homogeneous Velocity Model, the optimized schemes can outperform the standard $SpatTE$-$Cross(M)$ scheme for high-order FD schemes and proper mesh refinement. 
However, for velocity models with a broad velocity range, the optimized schemes are more sensitive because they depend on the maximum velocity; therefore, FD schemes that are independent of velocity offer the best trade-off between mesh refinement and FD order.

It is important highlight that these effects are consistently reflected in both error measures: while the $\|\cdot\|{2}$-norm evaluates the overall numerical error, in terms of a global value difference, the $\|\gamma^{\mathscr{F}}{u_{st}(x_{},z_{})}\|_{2}$-norm provides a more specific assessment of dispersion effects, where both measure can be used to define optimal choice of some FD scheme, considering the velocity field and discretization parameters.

\pagebreak 

\pagebreak

\pagebreak

\section{Conclusions}
In this work, we present a comprehensive mathematical construction for various types of FD schemes of arbitrary order and shape, highlighting the recent developments in this area through numerical results in Acoustic Wave propagation. We discuss the new construction of squared shaped stencils, as well as Rhombus and classic cross stencils, with and without dispersion optimization. Theoretical restrictions and possibilities are presented for each scheme. Overall, the schemes that require non-cross points (square and rhombus) require extra equations to be solved in the stencil coefficient computation, with the same order of accuracy as their cross counterpart. 

The potential gains of having extra non-cross points in the stencil are evaluated via an extensive numerical investigation. We propose analyses using Fourier and Wavelet Transform, which highlight key aspects of dispersion in each FD scheme, providing a more precise understanding of this effect than the numerical solution. Our numerical results, based on classical tests in the literature, provide a comprehensive understanding of the behavior of FD schemes across different idealized and realistic velocity models.

The implementation was performed in Devito \citep{devito-api}, therefore it is ready to be used in other production seismic modeling for anyone already using Devito.

Experiments with constant, or 2-layer media, show that optimized stencils can produce more accurate results for higher order at a similar cost compared to classic cross stencil. The use of extra points (non-cross points) is rarely advantageous.

In the presence of stronger velocity variations, such as in SEG EAGE and Marmousi, we observe that the choice of a unique time-space discretization-optimized FD scheme cannot achieve best results, unlike in cases where the velocity is constant or varies only slightly. Moreover, we conducted an extensive numerical test across different velocity fields and found that the choice of temporal and spatial parameters, as well as the velocity model, can influence the quality of the solution. Looking at the shapes, $Cross$ shapes have advantages over $Cross-Rb$ and $Cross-Sq$ shapes, due to their simplicity and computational cost, and because increasing the number of stencil weights does not significantly reduce the error. Regarding the methods, the use of $SpatTE$, $SpecTE-\theta$, $DispTE$, $DispLS$, or $SpecLS$ is all suitable; however, for complex velocity models, we see that the Cross $SpecLS$ (following \cite{liu2013time} and \cite{yang2006least}) schemes consistently achieve better accuracy at the same cost as other methods. The $SpecLS$ scheme does not depend on the velocity profile and time integration, only relying on the spatial computational setting ($h$), so it is robust and can be precomputed for efficient use in Acoustic RTM and FWI.

In the case of VTI/TTI acoustic wave equations, the FD schemes can be applied with some modifications. The schemes need to be adapted to look for unidirectional second-order derivatives. Our framework allows this adaptation when we select the $Cross$ shape, in which we need to select the unidirectional weights and recalculate the center weight. This adaptation is an important improvement to this model, in which numerical anisotropy and directional dispersion pose challenges that the optimized FD schemes can address satisfactorily. An open question is whether in the anisotropic case the non-cross stencils (rhombus, squares) can bring accuracy improvements that compensate its added cost.

%
%
%

In our main context, considering the Acoustic RTM and FWI pipeline, the use of optimized FD schemes can be directly applied to adjoint and gradient calculations, and their accurate approximation can be one of the keys to successful FWI convergence. From our results, schemes like $SpecLS$ have an advantage in this case due to the pre-calculation of the stencil and its application across the entire grid and over the overall step-times and FWI iterations, thereby reducing the computational effort, since the stencil can be stored in advance, independently of the timestep and velocity model.


\section*{Statements and Declarations}
\subsection*{Declaration of Competing Interest}
The authors declare that they have no known competing financial interests or personal relationships that could have appeared to influence the work reported in this paper.

\section*{Acknowledgments}
We gratefully acknowledge the technical guidance provided by TotalEnergies Research \& Technology (Houston, USA), with significant contributions of 
Dr. Christian Rivera, and the TotalEnergies E\&P Brazil Research Centre, especially the valuable and constructive discussions held throughout numerous technical meetings. This research was sponsored by TotalEnergies E\&P Brazil under the 1\% ANP obligation (AVENIR Project 24226‑3).

\begin{appendices}

\section{FFT and CWT based Norms}
Following \citep{torrence1998practical}, let be $u_{st}(x_{*},z_{*})$ a seismic trace in the spatial position $(x_{*},z_{*})$ over all the time times in $\Pi_{t}$. When we apply the DFT and/or CWT in $u_{st}(x_{*},z_{*})$ we will have:
\begin{equation}
\label{eqn:dft}
\mathscr{F}(u_{st}(x_{*},z_{*})) = \left\lbrace A_{k}\in\mathbb{C} / A_{k} =\left(\displaystyle\sum_{n=0}^{n_{t}-1}u_{*,*,n}e^{-2i\pi \frac{kn}{n_{t}}}\right)\mbox{, for } k=0,\ldots,n_{freq}.\right\rbrace,
\end{equation}
Where $i$ is the complex unit, note that $\mathscr{F}(u_{st}(x_{*},z_{*}))$ is a set of complex numbers, where each coefficient is associated with a one frequency $k$. To define the CWT, we need to determine the inverse of DFT, that is, the IDFT, where we will denote the IDFT by $\hat{\mathscr{F}}(\cdot)$, which recovers the original seismic trace $u_{st}(x_{*},z_{*})$. In this way, the IDFT is given by:
\begin{equation}
\label{eqn:idft}
\hat{\mathscr{F}}(u_{st}(x_{*},z_{*})) = \left\lbrace B_{n}\in\mathbb{R} / B_{n} =\left(\frac{1}{n_{freq}}\displaystyle\sum_{k=0}^{n_{freq}-1}A_{k}e^{-2i\pi\frac{kn}{n_{t}}}\right)\mbox{, for } n=0,\ldots,n_{t}-1.\right\rbrace,
\end{equation}
where $A_{k}$ is a coefficient obtained in \eqref{eqn:dft}. It is clear that $B_{n}=u_{*,*,n}$, for $n=0,\ldots,n_{t}$, that is, $\hat{\mathscr{F}}(u_{st}(x_{*},z_{*}))$ recovers the original seismic trace. To define de CWT we will introduce the value $W_{n}(s)$ given by:
\begin{equation}
\label{eqn:cwt1}
W_{n}(s) = \displaystyle\sum_{k=0}^{n_{freq}-1}A_{k}\hat{\Psi}^{*}(s\omega_{k})e^{i\omega_{k}n\tau}.   
\end{equation}
In the equation \eqref{eqn:cwt1} $s$ is a given scale value and $\omega_{k}$ is the angular frequency, given by:
\begin{equation}
\label{eqn:omegak}
\omega_{k} = \left\{
\begin{aligned}
    \displaystyle\frac{2\pi k}{n_{t}\tau} & ,& \hbox{if $k\leq\displaystyle\frac{n_{t}}{2}$,}\\
    \displaystyle-\frac{2\pi k}{n_{t}\tau} & ,& \hbox{if $k>\displaystyle\frac{n_{t}}{2}$.}
\end{aligned}
\right.
\end{equation}
Still in the equation \eqref{eqn:cwt1} we have $\hat{\Psi}^{*}(s\omega_{k})$ is the complex conjugate function given by:
\begin{equation}
\label{eqn:cwt2}
\hat{\Psi}^{*}(s\omega_{k}) = \displaystyle\left(\frac{2\pi s}{\tau}\right)^{1/2}\hat{\Psi}_{0}(s\omega_{k}),   
\end{equation}
where $\hat{\Psi}_{0}(s\omega_{k})$ is the Morlet continuous wavelet basis, with $\omega_{0}$ frequency, given by:
\begin{equation}
\label{eqn:cwt3}
\hat{\Psi}_{0}(s\omega_{k}) = \pi^{-1/4}H(\omega)e^{-(s\omega-\omega_{0})^{2}/2}, 
\end{equation}
Where $H(\omega)$ is the Heavisine function, given by:
\begin{equation}
\label{eqn:heavesidef}
H(\omega) = \left\{
\begin{aligned}
    1 & ,& \hbox{if $\omega>0$,}\\
    0 & ,& \hbox{otherwise.}
\end{aligned}
\right.
\end{equation}
Using the previous function definition, we have that the CWT is given by:
\begin{equation}
\label{eqn:cwt}
\mathscr{W}(u_{st}(x_{*},z_{*})) = \left\lbrace S_{n}\in\mathbb{C} / S_{n} = \hat{\mathscr{F}}(W_{n}) \mbox{, for } n=0,\ldots,n_{t}-1.\right\rbrace,
\end{equation}
There are other possibilities for the continuous wavelet basis, classified as real or complex, such as Mexican hat, Meyer, Shannon, Gaussian, Gabor, and others.  
Note that in \eqref{eqn:cwt} each $S_{n}$ is a set of complex numbers, for each step time $n=0,\ldots,n_{t}-1$. In this case, each $S_{n}$ is composed of $k$ coefficients and each coefficient is a complex number, which we will indicate by $S_{n,k}$.
As in the DFT, for each step time $n$, we have a set $S_{n,k}$ of complex numbers with $k$ coefficients, associated with each one a $ k$-frequency. 
Regarding our previous comment, the DFT works in the frequency domain, while CWT works in the time-frequency domain, as we can see mathematically by the earlier definitions. As DFT and CWT are complex numbers, we define the amplitude and phase for each transform as:
\begin{equation}
\label{eqn:ampldft}
\Gamma_{A_{k}} = \sqrt{\mathscr{Re}(A_{k}) + \mathscr{Im}(A_{k})},
\end{equation}
and
\begin{equation}
\label{eqn:amplcwt}
\Gamma_{S_{n,k}} = \sqrt{\mathscr{Re}(S_{n,k}) + \mathscr{Im}(S_{n,k})},
\end{equation}
where $\mathscr{Re}(\cdot)$ is the real part and $\mathscr{Im}(\cdot)$ the imaginary part of a complex number. The phase of each transformation is given by:
\begin{equation}
\label{eqn:phasedft}
\phi_{A_{k}} = \arctan\left(\frac{\mathscr{Im}(A_{k})}{\mathscr{Re}(A_{k})}\right),    
\end{equation}
and
\begin{equation}
\label{eqn:phasecwt}
\phi_{S_{n,k}} = \arctan\left(\frac{\mathscr{Im}(S_{n,k})}{\mathscr{Re}(S_{n,k})}\right),    
\end{equation}
where $\phi_{A_{k}}$ and $\phi_{S_{n,k}}$ are values in the interval $(-\pi,\pi]$. For the phases, it is easy to see that the maximum admissible value is $\pi$, so that we will define the maximum amplitude for each transformation as:
\begin{equation}
\label{eqn:amplmaxdft}
\Gamma_{A} = \max\limits_{0\leq k\leq n_{freq}-1}\left(\Gamma_{A_{k}}\right) ,
\end{equation}
and
\begin{equation}
\label{eqn:amplmaxcwt}
\Gamma_{S} = \max\limits_{0\leq k\leq n_{freq}-1,0\leq n\leq n_{t}-1}\left(\Gamma_{S_{n,k}}\right).
\end{equation}
In terms of units, we have that $[u_{st}(x_{*},z_{*})]=\displaystyle\frac{m}{s}$, so we have that $[\Gamma_{A_{k}}]=[\Gamma_{S_{n,k}}]=m$ and $[\phi_{A_{k}}]=[\phi_{S_{n,k}}]=\displaystyle\frac{rad}{s}$, where $m$ is meter, $s$ second and $rad$ radians. 

The previous construction is associated with $u_{st}(x_{*},z_{*})$, without considering whether this is a numerical or reference solution.  We will fix the spatial position $(x_{*},z_{*})$ and consider the associated numerical, $u_{st}^{num}$, and reference, $u_{st}^{ref}$, solutions without showing the position to reduce the complexity of the notation.
Each one of the previous values can be calculated for $u_{st}^{num}$ and $u_{st}^{ref}$, and to indicate the respective origin of the value, we will use the superscript \emph{num} or \emph{ref}. In terms of numerical and reference solutions, we have:
\begin{equation}
\label{eqn:dftnum}
\mathscr{F}(u_{st}^{num}(x_{*},z_{*})) = \left\lbrace A_{k}\in\mathbb{C} / A_{k} =\left(\displaystyle\sum_{n=0}^{n_{t}-1}u^{num}_{*,*,n}e^{-2i\pi \frac{kn}{n_{t}}}\right)\mbox{, for } k=0,\ldots,n_{freq}-1.\right\rbrace,
\end{equation}
\begin{equation}
\label{eqn:dftref}
\mathscr{F}(u_{st}^{ref}(x_{*},z_{*})) = \left\lbrace A'_{k}\in\mathbb{C} / A'_{k} =\left(\displaystyle\sum_{n=0}^{n_{t}-1}u^{ref}_{*,*,n}e^{-2i\pi \frac{kn}{n_{t}}}\right)\mbox{, for } k=0,\ldots,n_{freq}-1.\right\rbrace,
\end{equation}
\begin{equation}
\label{eqn:cwt1num}
W^{num}_{n}(s) = \displaystyle\sum_{k=0}^{n_{freq}-1}A_{k}\hat{\Psi}^{*}(s\omega_{k})e^{i\omega_{k}n\tau},   
\end{equation}
\begin{equation}
\label{eqn:cwt1ref}
W^{ref}_{n}(s) = \displaystyle\sum_{k=0}^{n_{freq}-1}A'_{k}\hat{\Psi}^{*}(s\omega_{k})e^{i\omega_{k}n\tau},   
\end{equation}
\begin{equation}
\label{eqn:cwtnum}
\mathscr{W}(u_{st}^{num}(x_{*},z_{*})) = \left\lbrace S_{n}\in\mathbb{C} / S_{n} = \hat{\mathscr{F}}(W_{n}^{num}) \mbox{, for } n=0,\ldots,n_{t}-1.\right\rbrace,
\end{equation}
\begin{equation}
\label{eqn:cwtref}
\mathscr{W}(u_{st}^{ref}(x_{*},z_{*})) = \left\lbrace S'_{n}\in\mathbb{C} / S'_{n} = \hat{\mathscr{F}}(W_{n}^{ref}) \mbox{, for } n=0,\ldots,n_{t}-1.\right\rbrace.
\end{equation}
The definitions \eqref{eqn:ampldft}-\eqref{eqn:amplmaxcwt} follow directly for the numerical and reference solutions using the previous notations, that is, we can obtain now $\Gamma_{A}^{num}$, $\Gamma_{A'}^{ref}$, $\Gamma_{S}^{num}$ and $\Gamma_{S'}^{ref}$ and their respectively components $\Gamma_{A_{k}}^{num}$, $\Gamma_{A'_{k}}^{ref}$, $\Gamma_{S_{n,k}}^{num}$ and $\Gamma_{S'_{n,k}}^{ref}$, also $\phi^{num}_{A_{k}}$, $\phi^{ref}_{A'_{k}}$, $\phi^{num}_{S_{n,k}}$ and $\phi^{ref}_{S'_{n,k}}$.
Considering the numerical and reference solution, we can define the phase difference as:
\begin{equation}
\label{eqn:phasesdifdft}
\phi_{dif}^{\mathscr{F}}(k) = \phi_{A_{k}}^{num} - \phi_{A'_{k}}^{ref}     
\end{equation}
\begin{equation}
\label{eqn:phasesdifcwt}
\phi_{dif}^{\mathscr{W}}(n,k) = \phi_{S_{n,k}}^{num} - \phi_{S'_{n,k}}^{ref}.    
\end{equation}
The phase difference can be rather than $\pi$ or less than $-\pi$. To resize the phase difference, we will use the following adjusted phase difference:
\begin{equation}
\label{eqn:phasesdifdftajusted}
\hat{\phi}_{dif}^{\mathscr{F}}(k) = \left\{
\begin{aligned}
    \phi_{dif}^{\mathscr{F}}(k) - 2\pi  & ,& \hbox{if $\phi_{dif}^{\mathscr{F}}(k) >+\pi$,}\\
    \phi_{dif}^{\mathscr{F}}(k)  + 2\pi & ,& \hbox{if $\phi_{dif}^{\mathscr{F}}(k)<-\pi$.}
\end{aligned}
\right. ,
\end{equation}
and
\begin{equation}
\label{eqn:phasesdifcwtajusted}
\hat{\phi}_{dif}^{\mathscr{W}}(n,k) = \left\{
\begin{aligned}
    \phi_{dif}^{\mathscr{W}}(n,k) - 2\pi  & ,& \hbox{if $\phi_{dif}^{\mathscr{W}}(n,k) >+\pi$,}\\
    \phi_{dif}^{\mathscr{W}}(n,k)  + 2\pi & ,& \hbox{if $\phi_{dif}^{\mathscr{W}}(n,k)<-\pi$.}
\end{aligned}
\right.
\end{equation}
The values $\phi_{dif}(\cdot)$ are resized to $(-\pi,\pi]$ using the definition of $\hat{\phi}_{dif}(\cdot)$. All the previous values are real numbers, even though the source of the values is complex numbers. So, we define the following dimensionless values:
\begin{equation}
\label{eqn:gammadft}
\gamma_{k}^{\mathscr{F}} = \displaystyle\left(\frac{\hat{\phi}_{dif}^{\mathscr{F}}(k)}{\pi}\right)\left(\frac{|\Gamma^{num}_{A_{k}}|}{\Gamma_{A}^{num}}\right),
\end{equation}
and
\begin{equation}
\label{eqn:gammacwt}
\gamma_{n,k}^{\mathscr{W}} = \displaystyle\left(\frac{\hat{\phi}_{dif}^{\mathscr{W}}(n,k)}{\pi}\right)\left(\frac{|\Gamma^{num}_{S_{n,k}}|}{\Gamma_{S}^{num}}\right).
\end{equation}
The values \eqref{eqn:gammadft} and \eqref{eqn:gammacwt} are dimensionless values and show an essential property of the numerical solution: as we said, the phase of the numerical and/or reference solution is connected to the dispersion properties of each solution. In this way, the phase difference shows us how delayed or advanced the phase of the reference solution is, that is, how dispersive the numerical solution is from the reference. Dividing by $\pi$ (that is, the maximum admissible phase difference value), we have a proportion of dispersion. 
On the other hand, the absolute value of the amplitude indicates to us how dissipative the solution is, and dividing by the maximum amplitude, we have a proportion of this dissipation. 
Moreover, each value $\gamma_{k}^{\mathscr{F}}$ or $\gamma_{n,k}^{\mathscr{W}}$ is calculated at one frequency, in the frequency or time frequency domain, respectively. Defining the following sets:
\begin{equation}
\label{eqn:gammadftset}
\gamma^{\mathscr{F}}_{u_{st}(x_{*},z_{*})} = \displaystyle\left\lbrace x\in\mathbb{R} / x = \displaystyle\left(\frac{\hat{\phi}_{dif}^{\mathscr{F}}(k)}{\pi}\right)\left(\frac{|\Gamma^{num}_{A_{k}}|}{\Gamma_{A}^{num}}\right)\mbox{, for } k=0,\ldots,n_{freq}-1.   \right\rbrace,
\end{equation}
and
\small{
\begin{equation}
\label{eqn:gammacwtset}
\begin{array}{ll}
\gamma^{\mathscr{W}}_{u_{st}(x_{*},z_{*})} = \\ \\ \displaystyle\left\lbrace x\in\mathbb{R} / x = \displaystyle\left(\frac{\hat{\phi}_{dif}^{\mathscr{W}}(k)} {\pi}\frac{|\Gamma^{num}_{S_{n,k}}|}{\Gamma_{S}^{num}}\right)\mbox{, for } n=0,\ldots,n_{t}-1 \mbox{ and } k=0,\ldots,n_{freq}-1.  \displaystyle\right\rbrace.    
\end{array}
\end{equation}
}
Note that \eqref{eqn:gammadftset} can be interpreted as a one-dimensional curve, where for each $k$ value we have a value $\gamma_{k}^{\mathscr{F}}$, that is, where we can evaluate the dispersion behavior for each $k$ frequency. Similarly, \eqref{eqn:gammacwtset} can be interpreted as a bi-dimensional curve, where for each $(n,k)$ pair we have a value $\gamma_{k}^{\mathscr{W}}$, that is, where we can evaluate the dispersion behavior for each $k$ frequency and each step time $n$.  

Our previous construction leads us to construct dimensionless curves for evaluating dispersion effects (for a suitable choice of FD scheme) on an arbitrary seismic trace $u_{st}(x_{*},z_{*})$. We can naturally extend this idea to a seismogram, that is, let $A_{xz}\subset\Pi_{xz}$ be an arbitrary set of selected points $(x_{p},z_{p})$, indexed by $p$. Regarding the definition of receiver, we can directly define the $\gamma^{\mathscr{F}}_{(\cdot)}$ and $\gamma^{\mathscr{W}}_{(\cdot)}$ by:
\begin{equation}
\label{eqn:gammadftrec}
\gamma^{\mathscr{F}}_{u_{rec}(A_{xz})} = \left\{\left. \gamma^{\mathscr{F}}_{u_{st}(x_{p},z_{p})}\right|_{\forall (x_{p},z_{p})\in A_{xz}}\right\},    
\end{equation}
and
\begin{equation}
\label{eqn:gammacwtrec}
\gamma^{\mathscr{W}}_{u_{rec}(A_{xz})} = \left\{\left. \gamma^{\mathscr{W}}_{u_{st}(x_{p},z_{p})}\right|_{\forall (x_{p},z_{p})\in A_{xz}}\right\}.  
\end{equation}
Note that now we can evaluate the behavior of a larger number of seismic traces, and looking at each curve can be unreasonable. Avoiding anTo avoid an individual evaluation, we can associate a number to each curve $\gamma^{\mathscr{F}}_{(\cdot)}$ and $\gamma^{\mathscr{W}}_{(\cdot)}$ e usual norms. We define a 2-norm based in $\gamma^{\mathscr{F}}_{(\cdot)}$ and $\gamma^{\mathscr{W}}_{(\cdot)}$ by: 
\begin{equation}
\label{eqn:norm2gammastdft}
\|\gamma^{\mathscr{F}}_{u_{st}(x_{*},z_{*})}\|_{2} = \displaystyle\sqrt{\frac{1}{n_{t}}\sum_{k=0}^{n_{t}}\left(\gamma_{k}^{\mathscr{F}}\right)^{2}}, 
\end{equation}
and
\begin{equation}
\label{eqn:norm2gammastcwt}
\|\gamma^{\mathscr{W}}_{u_{st}(x_{*},z_{*})}\|_{2} = \displaystyle\sqrt{\frac{1}{n_{t}}\frac{1}{n_{freq}}\sum_{n=0}^{n_{t}-1}\sum_{k=0}^{n_{freq}-1}\left(\gamma_{n,k}^{\mathscr{W}}\right)^{2}}. 
\end{equation}
The norms \eqref{eqn:norm2gammastdft} and \eqref{eqn:norm2gammastcwt} evaluate the square difference of phase difference in the frequency and time frequency domains, respectively, and the numbers \eqref{eqn:norm2gammastdft} and \eqref{eqn:norm2gammastcwt} are the average behavior of the dispersion effect. Regarding the FD schemes, the previous norms are indicators of how dispersive these schemes are, considering one seismic trace. Naturally, we can extend these norms for the receivers, using an average or maximum approach in the following way:
\begin{equation}
\label{eqn:norm2gammarecdftavg}
\|\gamma^{\mathscr{F}}_{u_{rec}(A_{xz})}||_{2}^{avg} = \displaystyle\frac{1}{n_{rec}}\sum_{p=1}^{n_{rec}} \| \gamma^{\mathscr{F}}_{u_{st}(x_{p},z_{p})} \|_{2} 
\end{equation}
\begin{equation}
\label{eqn:norm2gammarecdftmax}
\|\gamma^{\mathscr{F}}_{u_{rec}(A_{xz})}||_{2}^{max} = \max\limits_{1\leq p \leq n_{rec}}
\| \gamma^{\mathscr{F}}_{u_{st}(x_{p},z_{p})} \|_{2} 
\end{equation}
and
\begin{equation}
\label{eqn:norm2gammareccwtavg}
\|\gamma^{\mathscr{W}}_{u_{rec}(A_{xz})}||_{2}^{avg} = \displaystyle\frac{1}{n_{rec}}\sum_{p=1}^{n_{rec}} \| \gamma^{\mathscr{W}}_{u_{st}(x_{p},z_{p})} \|_{2} 
\end{equation}
\begin{equation}
\label{eqn:norm2gammareccwtmax}
\|\gamma^{\mathscr{W}}_{u_{rec}(A_{xz})}||_{2}^{max} = \max\limits_{1\leq p \leq n_{rec}}
\| \gamma^{\mathscr{W}}_{u_{st}(x_{p},z_{p})} \|_{2} 
\end{equation}

The norms \eqref{eqn:norm2gammarecdftavg} and \eqref{eqn:norm2gammareccwtavg} show us a global average behavior of the receivers associated with the set $A_{xz}$, while the norms \eqref{eqn:norm2gammarecdftmax} and \eqref{eqn:norm2gammareccwtmax} show us the maximum behavior of the receivers. The advantage of \eqref{eqn:norm2gammarecdftavg} and \eqref{eqn:norm2gammareccwtavg} is that we can compute a global behavior compensating regions where the dispersion effect is more substantial, that is, the numerical solution is out of phase (delay or advanced) and/or the wave peak has far value (increased or decreased) about the reference solution. The norms \eqref{eqn:norm2gammarecdftmax} and \eqref{eqn:norm2gammareccwtmax} can highlight a specific behavior of one seismic trace, without considering the contribution of the other traces.

\end{appendices}

\bibliography{silva_peixoto_refs}

\end{document}